\documentclass[11pt]{article}
\def\version{27.5.2016}\def\users{}  %
\def\users{final-layout}  
\usepackage{amsmath,amsfonts,amssymb,color}
\usepackage{mathrsfs} 
\usepackage{epsfig}
\usepackage{psfrag}
\usepackage{upgreek} 
\usepackage{cite}

\newtheorem{theorem}{Theorem}[section]

\newtheorem{definition}[theorem]{Definition}

\newtheorem{proposition}[theorem]{Proposition}

\newtheorem{remark}[theorem]{Remark}

\usepackage{ifthen}
\ifthenelse{\equal{\users}{final-layout}}{}{
\usepackage{fancyhdr}
\pagestyle{fancy}
\headheight=28pt\headwidth=16.5cm
\definecolor{gray}{gray}{0.5}
\rhead{\color{gray}Implicit time-discretisation in visco-elasto-dynamics\\
C.G.Panagiotopoulos, T.Roub\'\i\v cek}
\chead{}
\lhead{Version \version, file: \jobname.tex \\
compiled:
\number\day.\number\month.\number\year\ at
\the\hour:\ifnum\minute<10 0\fi\the\minute\ h }
\usepackage[notcite,notref,color]{showkeys}
\definecolor{labelkey}{rgb}{1.,.2,0.}

}
\newcount\hour \newcount\minute
\hour=\time
\divide \hour by 60
\minute=\time
\loop \ifnum \minute > 59 \advance \minute by -60 \repeat

\usepackage[normalem]{ulem}
\definecolor{brown}{rgb}{0.5,0,0}
\ifthenelse{\equal{\users}{final-layout}}{

    \newcommand{\DELETE}[1]{}
    \newcommand{\CHECK}[1]{#1}
    \newcommand{\COMMENT}[1]{}
   \newcommand{\COLOR}[1]{#1}
    \newcommand{\TINY}[1]{}
    \newcommand{\MARGINOTE}[1]{}
}
{

 \newcommand{\DELETE}[1]{{\color{brown}\sout{#1}\color{black}}}
  \newcommand{\CHECK}[1]{\color{brown}\uwave{#1}\color{black}}
 \newcommand{\COMMENT}[1]{{\color{red}\uuline{#1}\color{black}}}
 \newcommand{\COLOR}[1]{{\color{blue}{#1}}}
 \newcommand{\TINY}[1]{{\tiny#1}}
 \newcommand{\MARGINOTE}[1]{\marginpar{\color{red}\tiny\texttt{#1}}}
}

\renewcommand\dot[1]{\mathchoice
                 {{\buildrel{\hspace*{.1em}\text{\LARGE.}}\over{#1}}}
                 {{\buildrel{\hspace*{.1em}\text{\Large.}}\over{#1}}}
                 {{\buildrel{\hspace*{.1em}\text{\large.}}\over{#1}}}
                 {{\buildrel{\hspace*{.1em}\text{\large.}}\over{#1}}}}
\newcommand\DT{\dot}
\newcommand\DDT[1]{\mathchoice
   {{\buildrel{\hspace*{.1em}\text{\LARGE.\hspace*{-.1em}.}}\over{#1}}}
   {{\buildrel{\hspace*{.1em}\text{\Large.\hspace*{-.1em}.}}\over{#1}}}
   {{\buildrel{\hspace*{.1em}\text{\large.\hspace*{-.1em}.}}\over{#1}}}
   {{\buildrel{\hspace*{.1em}\text{\large.\hspace*{-.1em}.}}\over{#1}}}}

\def\R{{\mathbb R}}
\def\N{{\mathbb N}}
\newcommand\bbC{\mathbb C}
\newcommand\bbH{\mathbb H}
\newcommand\bbD{\mathbb D}

\newcommand\SYM{\R_\mathrm{sym}^{d\times d}}
\newcommand\GC{\mathchoice{\Gamma_{\hspace*{-.15em}\mbox{\tiny\rm C}}}
                          {\Gamma_{\hspace*{-.15em}\mbox{\tiny\rm C}}}
                          {\Gamma_{\hspace*{-.05em}\mbox{\tiny\rm C}}}
                          {\Gamma_{\hspace*{-.05em}\mbox{\tiny\rm C}}}}
\newcommand\GN{\mathchoice{\Gamma_{\hspace*{-.15em}\mbox{\tiny\rm N}}}
                          {\Gamma_{\hspace*{-.15em}\mbox{\tiny\rm N}}}
                          {\Gamma_{\hspace*{-.05em}\mbox{\tiny\rm N}}}
                          {\Gamma_{\hspace*{-.05em}\mbox{\tiny\rm N}}}}
\newcommand\GNN{\mathchoice{\Gamma_{\hspace*{-.15em}\mbox{\tiny\rm N1}}}
                          {\Gamma_{\hspace*{-.15em}\mbox{\tiny\rm N1}}}
                          {\Gamma_{\hspace*{-.05em}\mbox{\tiny\rm N1}}}
                          {\Gamma_{\hspace*{-.05em}\mbox{\tiny\rm N1}}}}
\newcommand\GNNN{\mathchoice{\Gamma_{\hspace*{-.15em}\mbox{\tiny\rm N0}}}
                          {\Gamma_{\hspace*{-.15em}\mbox{\tiny\rm N0}}}
                          {\Gamma_{\hspace*{-.05em}\mbox{\tiny\rm N0}}}
                          {\Gamma_{\hspace*{-.05em}\mbox{\tiny\rm N0}}}}
\newcommand{\Sdir}{\Sigma_{\mbox{\tiny\rm D}}}
\newcommand{\Gdir}{\Gamma_{\mbox{\tiny\rm D}}}
\newcommand{\GD}{\Gdir}

\newcommand{\Snew}{\Sigma_{\mbox{\tiny\rm N}}}

\newcommand{\nablaS}{\nabla_{\scriptscriptstyle\textrm{\hspace*{-.3em}S}}^{}}
\newcommand{\divS}{\mathrm{div}_{\scriptscriptstyle\textrm{\hspace*{-.1em}S}}^{}}

\renewcommand\d{\mathrm d}

\begin{document}
\begin{sloppypar}
\def\rmn{{\rm n}}
\def\rmt{{\rm t}}
\def\pl{{\partial}}
\def\In{{\in}}
\newcommand{\eq}[1]{\eqref{#1}}
\def\uhalf{{u_\tau^{k-1/2}}}
\def\vhalf{{v_\tau^{k-1/2}}}
\def\zhalf{{z_\tau^{k-1/2}}}
\def\pihalf{{\pi_\tau^{k-1/2}}}
\newcommand{\QED}{\ \hfill$\ \hfill\Box$\medskip}
\newcommand{\overf}{\hspace*{.08em}\overline{\hspace*{-.08em}f}}

\MARGINOTE{POSSIBLE~JOURNALS:
\\SIAM~J.Sci.Comp.
\\Numer.Funct.Anal.Optim.
\\SIAM~J.Numer.Anal.,
\\DCDS-S,
\\.......???}

\vspace*{2cm}

%

\baselineskip=24pt\noindent
{\LARGE\sf 
\COLOR{ENERGY-CONSERVING} TIME-DISCRETISATION OF\\
ABSTRACT DYNAMIC PROBLEMS WITH APPLICATIONS\\
IN 
CONTINUUM MECHANICS OF SOLIDS }\\[10mm]


\bigskip

\noindent
{\large\sf Tom{\'a}{\v s} Roub{\'\i}{\v c}ek}$\,^{1,2}$, {\large\sf Christos G. Panagiotopoulos}$\,^{3}$

\bigskip
\bigskip
\baselineskip=10pt

\noindent{\footnotesize\it
$^1$ Mathematical Institute, Charles University,
Sokolovsk{\'a}~83, CZ-186~75~Praha~8, Czech Republic.\\
$^2$ Institute of Thermomechanics, Czech Acad.\ Sci., Dolej\v skova~5,
CZ--182~00 Praha 8, Czech Republic\\
$^3$ Institute of Applied and Computational Mathematics, 
Foundation for Research and Technology - Hellas,\\\hspace*{.5em} 
Nikolaou Plastira 100, Vassilika Vouton, GR-700 13 Heraklion, Crete, Greece 
}


\bigskip
\bigskip

\noindent\hspace*{1cm}
\begin{minipage}[t]{14.3cm}

\baselineskip=12pt

{\small
\noindent{\bf Abstract.} 
An abstract 2nd-order evolution equation or inclusion is discretised in time
in such a way that the energy is conserved at least in qualified cases, 
typically in the cases when the governing energy is component-wise quadratic
\COLOR{or ``slightly-perturbed'' quadratic}. Specific 
applications in continuum mechanics of solids possibly with various internal 
variables cover vibrations or waves in 
linear viscoelastic materials 
at small strains, coupled with some inelastic processes as plasticity, 
damage, or phase transformations, and also some surface variants related 
to contact mechanics. The applicability is illustrated by numerical simulations 
of vibrations interacting with a frictional contact or waves emitted by an 
adhesive contact of a 2-dimensional viscoelastic body.  

\medskip


\medskip

} 
\end{minipage}


\vspace*{1em}
\baselineskip=14pt
\section{\sf INTRODUCTION}\label{sect-intro}

Dynamical effects due to inertia play a prominent role in a lot of application 
of continuum mechanics, and their interaction with various other processes 
is of a particular interest. Various nonlinear (possibly activated) processes 
as plasticity, damage, or phase transformations may thus mutually interact 
with these dynamical forces.

In computational mechanics, so-called transient versus 
wave propagation problems (i.e.\ low-frequency vibrations versus 
high-frequency waves, respectively) are distinguished
and different numerical methods are used especially as far as time 
discretisation concerns. Here we focus on implicit
methods, which are also well fitted with mathematical analysis.
It is well known that the backward Euler (also called backward-difference) 
scheme serves well for theoretical purposes as the so-called Rothe 
method (see e.g.\ \cite[Chap.5]{MieRou15RIST})
but brings unacceptable artificial numerical 
attenuation by the inertial term which practically
prevents any usage for wave computations even if an extremely
small time stepping is used. Henceforth, we focus on discretisation schemes
that conserves energy at least at some occasions.

For clarity, we make the exposition of the main ideas on an abstract level by 
considering the initial-value problem for an evolution system of two equations 
(or, if $\Psi$'s or $\Phi$ are not smooth, inclusions):
\begin{subequations}\label{abstract-hyper-system}
\begin{align}\label{abstract-hyper-system-a}
\mathscr{T}'\DDT u+&\pl\Psi_1(\DT u)+\pl^{}_u\Phi(u,z)\ni f(t),
\ \ \ \ \ \ u|_{t=0}^{}=u_0,\ \ \ \DT u|_{t=0}^{}=v_0,
\\\label{abstract-hyper-system-b}
&\pl\Psi_2(\DT z)+\pl^{}_z\Phi(u,z)\ni g(t),\ \ \ \ \ \ z|_{t=0}^{}=z_0,
\end{align}\end{subequations}
where $\mathscr{T}:H\to\R$, $\Phi:U\times Z\to\R$, 
$\Psi_1:U\to\R$, and $\Psi_2:Z\to\R$ are the 
(quadratic) kinetic-energy functional, 
a stored energy functional, and two dissipation functionals,
respectively, and $f:[0,T]\to U^*$ and $g:[0,T]\to Z^*$ are
the time-dependent loadings, while $U\subset H=H^*$ and $Z$ are
Banach spaces and ``$\partial$'' denotes the convex subdifferential
of the possibly nonsmooth functionals $\Psi_i$ or a partial G\^ateaux
derivative. If the functional is smooth, then the G\^ateaux
derivative will be briefly denoted by $(\cdot)'$; it is used for 
$\mathscr{T}'$ in \eqref{abstract-hyper-system-a} \COLOR{which is 
linear so that we write $\mathscr{T}'\DDT u$ instead of $\mathscr{T}'(\DDT u)$}.

For efficient computations of 
dynamical processes, 
\COLOR{various} integration methods, more 
sophisticated in comparison with the mentioned 
backward-Euler one,  have been devised by Newmark \cite{Newm59MCSD} 
and further generalized by Hilber, Hughes, and Taylor \cite{HiHuTa77INDT} 
and then used widely in engineering and computational physics. In fact, 
for a special choice of parameters, the latter method gives the classical 
Crank-Nicolson scheme \cite{CraNic47PMNE} if applied to a transformed system of 
three 1st-order equations (inclusions)
\begin{subequations}\label{abstract-2nd}\begin{align}\label{abstract-2nd-v}
&\DT u=v\,,&&u|^{}_{t=0}\!=u_0,&&
\\\label{abstract-2nd-u}
&\mathscr{T}'\DT v+\pl\Psi_1(v)+\pl^{}_u\Phi(u,z)\ni f\,,
&&v|^{}_{t=0}\!=v_0,
\\\label{abstract-2nd-b}
&\pl\Psi_2(\DT z)+\pl^{}_z\Phi(u,z)\ni g\,,&&z|_{t=0}^{}=z_0,
\end{align}\end{subequations}
cf.\ also e.g.\ \cite{ZWHT84USSS}. The Crank-Nicolson scheme
was originally devised for heat equation and later used for 
2nd-order problems in the form \eqref{abstract-hyper-system}, see 
e.g.\ \cite[Ch.6, Sect.9]{GlLiTr81NAVI}. It is different if applied to 
the dynamical equations transformed into the form \eq{abstract-2nd}; then
it is sometimes called just a central-difference scheme or 
generalized midpoint
scheme, cf.\ e.g.\ \cite[Sect.\,12.2]{Wang07FSW} 
or \cite[Sect.\,1.6]{SimHug98CI}, respectively.

If $\Phi$ is quadratic, this method conserves energy even
after time discretisation. It can be used for $\Psi_2$ 1-homogeneous 
\CHECK{when the limit passage in the discrete semi-stability 
developed for the backward Euler scheme in the modern
theory of rate-independent processes is suitably modified}. 
\COMMENT{THIS MAY NOT WORK - JUST A WEAK SOLUTION - I WILL FIX IT}
This will be done in Section~\ref{sec-quadratic}.
This already allows e.g.\ for application to visco-elasto-dynamic 
problems coupled with plasticity at small strains like 
in \cite[Sect.\,1.6.1]{SimHug98CI} or \cite{Wang07FSW},
as outlined in Section~\ref{sec-cont-mech}.

If $\Phi$ is only separately quadratic,
one can still combine the above outlined Crank-Nicolson
type discretisation with 
the fractional-step method to obtain a 
suitably decoupled scheme using two Crank-Nicolson
formulas.
This will be done in Section~\ref{sec-nonconvex}.
In special case, even generalization for  potentials 
which are nonquadratic in $u$ or in $z$ can be devised by using 
suitably defined differential quotients, which will be done on
an abstract level in Section~\ref{sec-nonquadratic}.

The resting plan of the paper is to discuss convergence of the devised 
discretisations under suitable data qualification in Section~\ref{sec-conv}
to outline application in continuum mechanics of solids under small strains
in Section~\ref{sec-cont-mech}, and then, in Section~\ref{sec-comput}, 
to demonstrate computational 
efficiency of such discretisations on a 2-dimensional dynamic adhesive 
Mode-II contact of a visco-elastic body with a rigid obstacle.

\TINY{
\bigskip\hrule{\tt
Simo ....
Numerical analysis and simulation of plasticity
JC Simo - Handbook of numerical analysis, 1998 - Elsevier

\cite[Sect.\,1.6.1]{SimHug98CI}

N.D. Cristescu: Dynamic plasticity - 2007 - World Scientific

M. L. Wilkins: Computer Simulation of Dynamic Phenomena

the Crank-Nicolson also in:

O. C. Zienkiewicz, W. L. Wood, N. W. Hine, and R. L. Taylor
A unified set of single step algorithms. Part 1: General formulation and applications.
International Journal for Numerical Methods in Engineering
1984
Volume 20, Issue 8, pages 1529-1552


}\hrule}


\vspace*{1em}
\section{\sf A CRANK-NICOLSON SCHEME FOR $\Phi$ QUADRATIC}\label{sec-quadratic}

Rather for notational simplicity, we consider a time step $\tau>0$
which do not vary in particular time levels, leading 
to an equidistant partition of the considered time interval.
Let us emphasize that, however, a varying time-step and non-equidistant 
partitions are easily possible to implement because we will always 
consider only first-order time differences and one-step formulas.
In fact, such a varying time-step can be advantageously used for a certain 
adaptivity to optimize computational costs.

The backward Euler formula applied to \eqref{abstract-hyper-system} serves 
well for theoretical purposes even for 2nd-order systems under 
the name ``Rothe method'' (cf.\ e.g.\ \cite{MieRou15RIST,Roub13NPDE})
but it is well known that, due to the artificial numerical attenuation, 
it practically
cannot be used for realistic calculation of wave propagation unless $\tau>0$
is made extremely small. Many other methods have been devised for 
dynamical problems. A simple method
consists in application of the \emph{Crank-Nicolson formula} 
\cite{CraNic47PMNE} to the transformed 1st-order system \eq{abstract-2nd}.
This results to the system for the triple $(u_\tau^k,v_\tau^k,z_\tau^k)$:
\begin{subequations}\label{C-N-2nd}
\begin{align}\label{C-N-2nd+}
&\frac{u_\tau^k{-}u_\tau^{k-1}\hspace*{-.5em}}\tau=\vhalf,&&
u_\tau^0=u_0,
\\&\label{C-N-2nd++}
\mathscr{T}'\frac{v_\tau^k{-}v_\tau^{k-1}\hspace*{-.5em}}\tau
+\pl\Psi_1\big(v_\tau^{k-1/2}\big)
+\pl^{}_u\Phi\big(\uhalf,\zhalf\big)\ni
f_\tau^k\ ,
&&v_\tau^0=v_0,&&
\\\label{C-N-2nd+++}
&\pl\Psi_2\Big(\frac{z_\tau^k{-}z_\tau^{k-1}\hspace*{-.5em}}\tau\hspace*{.5em}\Big)
+\pl^{}_z\Phi\big(u_\tau^{k-1/2},z_\tau^{k-1/2}\big)\ni g_\tau^k\ ,
&&z_\tau^0=z_0,&&
\\\label{half-notation}
&\text{with the abbreviation } 
\uhalf:=\frac{u_\tau^k{+}u_\tau^{k-1}\hspace*{-.5em}}2,
\ \,\vhalf:=\frac{v_\tau^k{+}v_\tau^{k-1}\hspace*{-.5em}}2,
\ \,\zhalf:=\frac{z_\tau^k{+}z_\tau^{k-1}\hspace*{-.5em}}2\hspace*{-7em}
\end{align}\end{subequations}
to be solved recursively for $k=1,...,T/\tau$. 
The 
right-hand sides can be taken e.g.\ as 
$f_\tau^k:=f(k\tau)$, or $f_\tau^k:=f((k{-}\frac12)\tau)$,
or $f_\tau^k:=\frac12f(k\tau)+\frac12f((k{-}1)\tau)$,
or $f_\tau^k:=\frac1\tau\int_{(k{-}1)\tau}^{k\tau}f(t)\,\d t$,
and similarly for $g_\tau^k$, the particular choice not being essential
for our considerations below. This approximation was suggested in engineering 
literature for the system without $z$-variable, i.e.\ for 
(\ref{C-N-2nd}a,b),
e.g.\ in \cite{PaPaTa09RDCT,WelHam81ASCP,ZWHT84USSS}, possibly even in a 
nonlinear variant 
using a non-centered adaptively tuned formula. 
Actually, it falls into a broader class of so-called 
Hilber-Hughes-Taylor formulas widening the popular Newmark method 
\cite{Newm59MCSD} as a special choice of parameters (namely $\alpha=\beta=1/2$ 
and $\gamma=1$ in the usual notation, cf.\ \cite{HiHuTa77INDT}). To be mentioned here that an extension of the Hilber-Hughes-Taylor method has been presented known as the generalized-$\alpha$ method \cite{ChuHu93INDGA} allowing high frequency energy dissipation with second order accuracy.
\COMMENT{the paper Hilber-Hughes-Taylor better to check for this notation}

This system does not satisfy the usual symmetry condition 
and thus does not have any potential, but eliminating $v_\tau^k$ by 
substituting $v_\tau^k=\frac2\tau(u_\tau^k{-}u_\tau^{k-1})-v_\tau^{k-1}$ into 
\eq{C-N-2nd++}, one again obtains a potential problem for 
the couple $(u_\tau^k,z_\tau^k)$. To be more specific, 
$(u_\tau^k,z_\tau^k)\in U\times Z$ is a minimizer of the functional
\begin{align}\nonumber
\!\!(u,z)\mapsto 
\frac2\tau\Phi\Big(\frac{u{+}u_\tau^{k-1}\hspace*{-.5em}}2\hspace*{.5em},\frac{z{+}z_\tau^{k-1}\hspace*{-.5em}}2\hspace*{.5em}\Big)
+\Psi_1\Big(\frac{u{-}u_\tau^{k-1}\hspace*{-.5em}}\tau\hspace*{.5em}\Big)
+\Psi_2\Big(\frac{z{-}z_\tau^{k-1}\hspace*{-.5em}}\tau\hspace*{.5em}\Big)\
\\+\,2\tau\mathscr{T}\Big(\frac{u{-}\tau v_\tau^{k-1}{-}u_\tau^{k-1}\hspace*{-.5em}}{\tau^2}\hspace*{.5em}\Big)
-\Big\langle f_\tau^k,\frac u\tau\Big\rangle
-\Big\langle g_\tau^k,\frac z\tau\Big\rangle\;,
\label{C-R-formula-2nd}
\end{align}
and then one simply calculate
$v_\tau^k=\frac2\tau(u_\tau^k{-}u_\tau^{k-1})-v_\tau^{k-1}\in H$.

Existence of a potential is thus also advantageous to ensure existence of 
a solution $(u_\tau^k,v_\tau^k,z_\tau^k)\in U\times H\times Z$ to \eq{C-N-2nd} just
by the direct method. In this section, we will assume:
\begin{subequations}\label{struct-ass}\begin{align}
&\Phi:U\times Z\to\R\ 
\text{ coercive, 
quadratic, 
}
\label{struct-ass-a}
\\\label{struct-ass-b}
&\mathscr{T}:H\to\R\ \text{ quadratic, coercive,
}\ \ 
\\\nonumber
&\Psi_1:U\to\R\cup\infty,\ \Psi_2:Z\to\R\cup\infty\ \ \text{ convex,} 
\\\label{struct-ass-c}
&\hspace*{5em}\text{lower semicontinuous,
$p_i$-coercive (i.e.\ $\Psi_i(\cdot)\ge\epsilon\|\cdot\|^{p_i}$)\,.}
\end{align}\end{subequations}
The  coercivity of $\Phi$ means that 
$\lim_{\|u\|\to\infty,\,\|z\|\to\infty}\Phi(u,z)/(\|u\|+\|z\|)=\infty$
and, in fact, can be weakened if combined with a coercivity of 
$\Psi_1$ or $\Psi_2$. 
Altogether, the potential in \eq{C-R-formula-2nd} is 
convex and weakly
lower semicontinuous, which ensures the mentioned existence of its  
minimizer which solves also \eq{C-N-2nd} 
provided $U\times Z$ is 
reflexive and 
\begin{align}\label{qualif-f-g-IC}
\text{$f_\tau^k\in U^*\ $ and $\ g_\tau^k\in Z^*\ $ and 
$\ (u_0,v_0,z_0)\in U\times H\times Z$.}
\end{align}

By testing \eq{C-N-2nd++} by $v_\tau^k{+}v_\tau^{k-1}$ and substituting 
also $v_\tau^k{+}v_\tau^{k-1}=\frac2\tau(u_\tau^k{-}u_\tau^{k-1})$ due to 
\eq{C-N-2nd+} and by testing \eq{C-N-2nd+++} by
$z_\tau^k{-}z_\tau^{k-1}$, after summation we obtain the equality
\begin{align}\nonumber
\frac{\mathscr{T}(v_\tau^k)-\mathscr{T}(v_\tau^{k-1})}\tau
+\Xi_1\big(\vhalf\big)+\Xi_2\Big(\frac{z_\tau^k{-}z_\tau^{k-1}\hspace*{-.5em}}\tau\hspace*{.5em}\Big)
+\frac{\Phi(u_\tau^k,z_\tau^k)-\Phi(u_\tau^{k-1},z_\tau^{k-1})}\tau
\\=\big\langle f_\tau^k,\vhalf\big\rangle+\Big\langle g_\tau^k,
\frac{z_\tau^k{-}z_\tau^{k-1}\hspace*{-.5em}}\tau\hspace*{.5em}\Big\rangle
\label{C-N-2nd-conserv}
\end{align}
with the dissipation rates defined by 
\begin{align}
\Xi_1(v):=\big\langle\pl\Psi_1(v),v\big\rangle\quad\text{ and }\quad
\Xi_2(\DT z):=\big\langle\pl\Psi_1(\DT z),\DT z\big\rangle,
\end{align}
and where we used the structural assumption that both $\mathscr{T}$
and $\Phi$ are quadratic.
More specifically, we used the two following binomial formulas:
\begin{subequations}\label{binomial}\begin{align}\label{binomial-T}
&\mathscr{T}'\frac{v_\tau^k{-}v_\tau^{k-1}\!\!}\tau\:\cdot
\frac{v_\tau^k{+}v_\tau^{k-1}}2
=\frac{\mathscr{T}(v_\tau^k)-\mathscr{T}(v_\tau^{k-1})}\tau,
\\\nonumber
&\Phi'_z\Big(\frac{u_\tau^k{+}u_\tau^{k-1}\!\!}2\ ,\frac{z_\tau^k{+}z_\tau^{k-1}\!\!}2\ \Big)
\cdot\Big(\frac{v_\tau^k{+}v_\tau^{k-1}\!\!}2\ ,\frac{z_\tau^k{-}z_\tau^{k-1}\!\!}\tau\ \Big)
=\Phi'_z\Big(\frac{u_\tau^k{+}u_\tau^{k-1}\!\!}2\ ,\frac{z_\tau^k{+}z_\tau^{k-1}\!\!}2\ \Big)\cdot
\\&\hspace*{11em}\cdot\Big(\frac{u_\tau^k{-}u_\tau^{k-1}\!\!}\tau\ ,\frac{z_\tau^k{-}z_\tau^{k-1}\!\!}\tau\ \Big)=
\frac{\Phi(u_\tau^k,z_\tau^k)-\Phi(u_\tau^{k-1},z_\tau^{k-1})}\tau.
\label{binomial-Phi}\end{align}\end{subequations}
In particular as a special case if $\Psi_1=\Psi_2=0$, $f=0$, and $h=0$, the
equality \eqref{C-N-2nd-conserv} shows that the discrete scheme 
\eq{C-N-2nd} conserves the kinetic and stored energy:
$\mathscr{T}(v_\tau^k)+\Phi(u_\tau^k,z_\tau^k)=\,$constant. 

The scheme investigated in this Section~\ref{sec-quadratic}
covers various linear rheological models (as Kelvin-Voigt's, Maxwell's,
Jeffreys', Burgers', etc.) which use the dissipation potentials $\Psi$'s 
quadratic, i.e.\ $p_1=p_2=2$, 
and could be easily implemented numerically \cite{PaMaRo14BVE}, 
cf.\ also Remark~\ref{rem-Max} below. 

The non-quadratic potentials $\Psi$'s, considered
above too, allow e.g.\ for modeling of certain inelastic processes. 
Even they can be nonsmooth at 0, i.e.\ these processes may be
activated in the sense that their evolution needs the corresponding
driving force to achieve a certain threshold. Moreover, 
some of them can be 1-homogeneous, so that $p_1=1$ or $p_2=1$, 
which means that these processes can
be rate independent, although the whole system remains rate dependent
due to the inertia and possibly also due to the other $\Psi$-potential. 

Although \eqref{struct-ass-a} still brings substantial restriction on 
generality, various variants of rate-independent 
linearized plasticity, with kinematic or isotropic hardening, or 
without hardening (perfect plasticity) or with plastic-strain gradient
are thus covered.
\TINY{, cf.\ also Remark~\ref{rem-elasto-plasto} below.}


\begin{remark}\label{rem-Psi}
\upshape
In fact, in this section the additive splitting 
$\Psi(\DT u,\DT z)=\Psi_1(\DT u)+\Psi_2(\DT z)$ can easily be avoided
and a general coupling of dissipative forces can be considered.
\end{remark}

\TINY{
\begin{remark}[{\it Set constraints:
combination with the backward-Euler scheme.}]\label{rem-nonquadratic-}
\upshape
A lot of applications involve various constraints in the stored 
energy. This augments the structure \eq{struct-ass} as follows:
\begin{subequations}\label{struct-ass+K}\begin{align}\nonumber
&\Phi=\Phi_0+\delta_K\ \ 
\text{ coercive, with }\ 
\delta_K(u,z)=\begin{cases}0&\text{if }(u,z)\in K,
\\\infty&\text{if }(u,z)\not\in K,\end{cases}
\\[-.4em]\label{struct-ass-K}&\qquad\text{ with $K\subset U\times Z$ convex, 
closed,
}
\\&\qquad\text{ and $\Phi_0:U\times Z\to\R$ quadratic, convex, 
}
\label{struct-ass-a+K}
\\\label{struct-ass-b+K}
&\mathscr{T}\text{ quadratic, convex},\ \ 
\Psi_i\text{ convex, lower semicontinuous, 
$p_i$-coercive (i.e.\ $\Psi_i\ge\epsilon\|\cdot\|^{p_i}$)
}.
\end{align}\end{subequations}
The Crank-Nicolson scheme \eqref{C-N-2nd} is combined with 
the backward-Euler scheme for the set constraints:
\begin{subequations}\label{C-N-2nd+}
\begin{align}\label{C-N-2nd+}
&\frac{u_\tau^k{-}u_\tau^{k-1}\hspace*{-.5em}}\tau=\vhalf,&&
u_\tau^0=u_0,
\\&\label{C-N-2nd++}
\mathscr{T}'\frac{v_\tau^k{-}v_\tau^{k-1}\hspace*{-.5em}}\tau
+\pl\Psi_1\big(v_\tau^{k-1/2}\big)
+\pl^{}_u\Phi_0\big(\uhalf,\zhalf\big)\ni
f_\tau^k-I_1^*N_{K
}(u_\tau^k,z_\tau^k),
&&v_\tau^0=v_0,&&
\\\label{C-N-2nd+++}
&\pl\Psi_2\Big(\frac{z_\tau^k{-}z_\tau^{k-1}\hspace*{-.5em}}\tau\hspace*{.5em}\Big)
+\pl^{}_z\Phi_0\big(u_\tau^{k-1/2},z_\tau^{k-1/2}\big)\ni g_\tau^k-I_2^*N_{K
}(u_\tau^k,z_\tau^k)\ ,
&&z_\tau^0=z_0,&&
\end{align}\end{subequations}
to be solved recursively for $k=1,...,T/\tau$, where we use the notation for
the embeddings $I_1:U\mapsto U\times Z$ and $I_2:Z\mapsto U\times Z$.
The potential \eq{C-R-formula-2nd} now involves $\delta_K(u,z)$.
Unfortunately, the energy balance \eq{C-N-2nd-conserv} 
now holds only as an inequality ``$\,\le\,$''. This is still enough for 
numerical stability of this discretisation and for its convergence in qualified
cases. The artificial attenuation caused by the backward-Euler discretisation of
these constraints is effective only if these constraints are ``in action'' and
presumably might not substantially destroy ability of this scheme 
to calculate vibrations.
\end{remark}
}


\vspace*{1em}
\section{\sf FRACTIONAL-STEP SPLITTING OF CRANK-NICOLSON SCHEME 
}\label{sec-nonconvex}

The quadratic (and in particular convex) structure of $\Phi$ considered in 
Section~\ref{sec-quadratic} is a severe restriction and excludes interesting 
applications. In particular, suddenly triggered processes (like rupture)
are hard to model. Thus relaxing this structural restriction is highly 
desirable. This can be achieved by decoupling the time-discretised system 
suitably, namely ``componentwise''. This allows to qualify $\Phi$ 
only ``componentwise'' and works successfully if the dissipation potentials 
$\Psi$'s are separated, as indeed the case of our system 
\eq{abstract-hyper-system}. It is called a fractional-step method or 
sometimes also a Lie-Trotter (or sequential) splitting, and there is 
an extensive literature about it, cf.\ \cite{Marc90SADM,Yane71MFS}.
Actually, the Crank-Nicolson scheme itself can be understood as a splitting, 
cf.\ \cite{Fara05SMAA}.

Let us first relax \eqref{struct-ass-a} by assuming that 
\begin{subequations}\label{structural-ass}\begin{align}
\Phi:U\times Z\to\R
\ \ \text{ coercive},\quad
\label{structural-ass-a}
&\forall z\in Z:\quad\Phi(\cdot,z):U\to\R\ \text{ quadratic, convex, 
}
\\\label{structural-ass-b}
&\forall u\in U:\quad\Phi(u,\cdot):Z\to\R\ \text{ quadratic, convex, 
}
\end{align}\end{subequations}
while (\ref{struct-ass}b,c) remains unchanged. 
We now modify \eqref{C-N-2nd} as follows:
\begin{subequations}\label{abstract-hyper-disc}
\begin{align}\label{abstract-hyper-disc-a}
&\frac{u_\tau^k{-}u_\tau^{k-1}\hspace*{-.5em}}\tau=\vhalf,&&
u_\tau^0=u_0,
\\\label{abstract-hyper-disc-b}
&\mathscr{T}'\frac{v_\tau^k{-}v_\tau^{k-1}\hspace*{-.5em}}\tau+\pl\Psi_1\big(\vhalf\big)
+\pl^{}_u\Phi\big(\uhalf,z_\tau^{k-1}\big)\ni f_\tau^k,&&v_\tau^0=v_0,
\\\label{abstract-hyper-disc-c}
&\pl\Psi_2\Big(\frac{z_\tau^k{-}z_\tau^{k-1}\hspace*{-.5em}}\tau\hspace*{.5em}\Big)
+\pl^{}_z\Phi(u_\tau^k,z_\tau^{k-1/2})\ni g_\tau^k,&&z_\tau^0=z_0,
\end{align}\end{subequations}
where we again use the notation \eqref{half-notation}. Note that the system 
\eq{abstract-hyper-disc} is indeed decoupled: first (\ref{abstract-hyper-disc}a,b) 
is to be solved for $(u_\tau^k,v_\tau^k)$ and 
then \eq{abstract-hyper-disc-c} is to be solved for $z_\tau^k$. Sometimes, this 
componentwise-split Crank-Nicolson method is also called the second-order
Yanenko method \cite{Fara05SMAA}.

Likewise \eqref{C-R-formula-2nd}, an algorithmically useful observation is that
these problems possess potentials but, in contrast to 
\eqref{C-N-2nd}, now two potentials are to be identified, namely
\begin{subequations}\label{C-R-system-potential}\begin{align}
&\label{C-R-system-potential-u}
\!\!\!\!u\mapsto\frac2\tau\Phi\Big(\frac{u{+}u_\tau^{k-1}\!\!}2,z_\tau^{k-1}\Big)
+\Psi_1\Big(\frac{u{-}u_\tau^{k-1}\hspace*{-.5em}}\tau\hspace*{.5em}\Big)
+2\tau\mathscr{T}\Big(\frac{u{-}\tau v_\tau^{k-1}{-}u_\tau^{k-1}\!\!}{\tau^2}\ \Big)
-\Big\langle f_\tau^k,\frac u\tau\Big\rangle,\!\!
\\&\label{C-R-system-potential-z}
\!\!\!\!z\mapsto
\frac2\tau\Phi\Big(u^k,\frac{z{+}z_\tau^{k-1}}2\Big)
+\Psi_2\Big(\frac{z{-}z_\tau^{k-1}\hspace*{-.5em}}\tau\hspace*{.5em}\Big)
-\Big\langle g_\tau^k,\frac z\tau\Big\rangle\:.
\end{align}\end{subequations}
Note that, under the assumptions \eqref{struct-ass-b} and 
\eqref{structural-ass}, both these potentials
are convex and coercive, which may algorithmically 
facilitate numerical solution of (\ref{abstract-hyper-disc}a,b) and 
\eqref{abstract-hyper-disc-c}.
And, like in Section~\ref{sec-quadratic}, existence of 
a solution to \eq{abstract-hyper-disc} is guaranteed by
the direct-method arguments provided again \eq{qualif-f-g-IC} holds.

To show energy conservation even in the discrete scheme,
we use the same test as we made for \eqref{C-N-2nd}, namely
$\vhalf$ for \eqref{abstract-hyper-disc-b}
and $\frac{z_\tau^k{-}z_\tau^{k-1}}\tau$ for \eqref{abstract-hyper-disc-c}.
Using \eqref{structural-ass},
it gives 
\begin{subequations}\label{C-N-2nd-system}\begin{align}
&\frac{\mathscr{T}(v_\tau^k)-\mathscr{T}(v_\tau^{k-1})}\tau
+\Xi_1\big(\vhalf\big)
+\frac{\Phi(u_\tau^k,z_\tau^{k-1})-\Phi(u_\tau^{k-1},z_\tau^{k-1})}\tau
=\big\langle f_\tau^k,\vhalf\big\rangle,
\label{C-N-2nd-system-a}
\\&
\Xi_2\Big(\frac{z_\tau^k{-}z_\tau^{k-1}\hspace*{-.5em}}\tau\hspace*{.5em}\Big)
+\frac{\Phi(u_\tau^k,z_\tau^k)-\Phi(u_\tau^k,z_\tau^{k-1})\hspace*{-.5em}}\tau
=\Big\langle g_\tau^k,\frac{z_\tau^k{-}z_\tau^{k-1}\!\!}\tau\ \Big\rangle.
\label{C-N-2nd-system-b}
\end{align}\end{subequations}
It is important that the scheme is carefully decoupled in such a way that,
when summing \eqref{C-N-2nd-system} up, we benefit from the 
cancellation of the terms $\pm\Phi(u_\tau^k,z_\tau^{k-1})$ and obtain again the 
energy equality \eqref{C-N-2nd-conserv}. Here we used
together with the binomial formula  \eqref{binomial-T}, other two binomial 
formulas instead of only one in \eqref{binomial-Phi}, namely
\begin{subequations}\label{2x-binomial}\begin{align}\nonumber
\pl_u\Phi\Big(\frac{u_\tau^k{+}u_\tau^{k-1}\!\!}2\ ,z_\tau^{k-1}\Big)
\cdot\frac{v_\tau^k{+}v_\tau^{k-1}\!\!}2 
&=\pl_u\Phi\Big(\frac{u_\tau^k{+}u_\tau^{k-1}\!\!}2\ ,z_\tau^{k-1}\Big)\cdot
\frac{u_\tau^k{-}u_\tau^{k-1}\!\!}\tau\ 
\\\hspace*{11em}&=
\frac{\Phi(u_\tau^k,z_\tau^{k-1})-\Phi(u_\tau^{k-1}z_\tau^{k-1})}\tau,\ \ \text{ and}
\label{2x-binomial-u}\\\label{2x-binomial-z}
\pl_z\Phi\Big(u_\tau^k,\frac{z_\tau^k{+}z_\tau^{k-1}\!\!}2\ \Big)
\,\cdot\,\frac{z_\tau^k{-}z_\tau^{k-1}\!\!}\tau\ &=
\frac{\Phi(u_\tau^k,z_\tau^k)-\Phi(u_\tau^k,z_\tau^{k-1})}\tau.
\end{align}\end{subequations}
%


\begin{remark}[{\it More general dissipation.}]\label{rem-more-dissip}
\upshape
Making the dissipation potentials $\Psi_1(\DT u)$ and $\Psi_2(\DT z)$ 
dependent also on the state $(u,z)$ is easy and widens applications.
Then the subdifferentials in \eq{abstract-hyper-system} and 
\eq{abstract-2nd} should be only
partial with respect to $\DT u$ and $\DT z$ respectively, and 
the discrete scheme \eq{C-N-2nd++} and \eq{abstract-hyper-disc} 
should use 
$\pl_{\DT u}\Psi_1(u_\tau^{k-1},z_\tau^{k-1},v_\tau^{k-1/2})$ and 
$\pl_{\DT z}\Psi_2(u_\tau^{k-1},z_\tau^{k-1},\frac{z_\tau^k{-}z_\tau^{k-1}}\tau)$
instead of $\pl\Psi_1(v_\tau^{k-1/2})$  and 
$\pl\Psi_2(\frac{z_\tau^k{-}z_\tau^{k-1}}\tau)$, respectively.
\end{remark}

\vspace*{1em}
\section{\sf SPECIAL NONQUADRATIC POTENTIALS $\Phi$
}\label{sec-nonquadratic}

\COLOR{
A generalization for nonquadratic cases is very desirable for some applications.
Still holding energy-conservation, it can sometimes be realized by a 
modification of (\ref{abstract-hyper-disc}b,c) in the spirit of 
\cite[Sec.\,3.1]{CoMeSu11SAPF} where a specific gradient-flow problem
or \cite{ArmPet98FACA} where a specific compliance contact 
have been considered. In contrast to e.g.\ \cite{WelHam81ASCP}, we confine 
ourselves on such special cases where we will not need any adaptively tuned 
formula needed iterative implementation.
In the abstract case we can assume existence of differential quotients
\begin{align}\mathfrak{D}_{u}^{}\Phi:U\times U\times Z\to U^*
\qquad\text{ and }\qquad
\mathfrak{D}_{z}^{}\Phi:U\times Z\times Z\to Z^*
\end{align}
approximating respectively $\pl_u\Phi$ and $\pl_z\Phi$ in the sense that 
\begin{subequations}\label{ass-for-Phi}\begin{align}\label{ass-for-Phi-a}
&\forall\,u\in U,\ z\in Z:\quad\ \
\mathfrak{D}_{u}^{}\Phi(u,u,z)=\pl_u\Phi(u,z)
\ \ \text{ and }\ \ \mathfrak{D}_{z}^{}\Phi(u,z,z)=\pl_z\Phi(u,z)\,,
\\&\nonumber\forall\,u,\tilde z\!\in\!U,\ z,\tilde z\!\in\!Z:\
\big\langle\mathfrak{D}_{u}^{}\Phi(u,\tilde u,z),u{-}\tilde u\big\rangle
=\Phi(u,z)-\Phi(\tilde u,z)\ \ \text{ and }
\\&\label{ass-for-Phi-b}\hspace*{9em}
\big\langle\mathfrak{D}_{z}^{}\Phi(u,z,\tilde z),z{-}\tilde z\big\rangle
=\Phi(u,z)-\Phi(u,\tilde z),
\\&\label{ass-for-Phi-cont}
\mathfrak{D}_{u}^{}\Phi,\mathfrak{D}_{z}^{}\Phi\ \text{are continuous, and possibly also }
\\&\label{ass-for-Phi-pot-u}
\exists\,\mathfrak{F}_{\tilde u,z}:U\to\R:\qquad
\mathfrak{D}_{u}^{}\Phi(u,\tilde u,z)=\pl\mathfrak{F}_{\tilde u,z}(u),
\ \text{ and }
\\&\label{ass-for-Phi-pot-z}
\exists\,\mathfrak{G}_{u,\tilde z}:Z\to\R:\qquad
\mathfrak{D}_{z}^{}\Phi(u,z,\tilde z)=\pl\mathfrak{G}_{u,\tilde z}(z),
\end{align}\end{subequations}
where the continuity assumption \eqref{ass-for-Phi-cont} refers to suitable 
topologies depending on particular situations. Then, omitting now the notation 
$\vhalf$, \eqref{abstract-hyper-disc} is to be modified for $k=1,...,T/\tau$ as
\begin{subequations}\label{abstract-hyper-disc+}
\begin{align}\label{abstract-hyper-disc-a+}
&\frac{u_\tau^k{-}u_\tau^{k-1}\hspace*{-.5em}}\tau=
\frac{v_\tau^k{+}v_\tau^{k-1}\hspace*{-.5em}}2,
&&u_\tau^0=u_0,
\\\label{abstract-hyper-disc-b+}
&\mathscr{T}'\frac{v_\tau^k{-}v_\tau^{k-1}\hspace*{-.5em}}\tau+\pl\Psi_1
\Big(\frac{u_\tau^k{-}u_\tau^{k-1}\hspace*{-.5em}}\tau\hspace*{.5em}\Big)
+
\mathfrak{D}_{z}^{}(u_\tau^k,u_\tau^{k-1},z_\tau^{k-1})\ni f_\tau^k,&&v_\tau^0=v_0,
\\\label{abstract-hyper-disc-c+}
&\pl\Psi_2\Big(\frac{z_\tau^k{-}z_\tau^{k-1}\hspace*{-.5em}}\tau\hspace*{.5em}\Big)
+
\mathfrak{D}_{z}^{}\Phi(u_\tau^k,
z_\tau^k,z_\tau^{k-1})\ni g_\tau^k,&& z_\tau^0=z_0.
\end{align}\end{subequations}
Obviously, \eqref{ass-for-Phi-b} ensures that the test of 
\eqref{abstract-hyper-disc-b+} and \eqref{abstract-hyper-disc-c+} successively by
$\frac{u_\tau^k-u_\tau^{k-1}}\tau$ and 
$\frac{z_\tau^k-z_\tau^{k-1}}\tau$ gives again \eqref{C-N-2nd-system}
and one can again benefit from cancellation of the ``mixed'' terms 
$\pm\Phi(u_\tau^k,z_\tau^{k-1})$ when summing  \eqref{C-N-2nd-system} 
up. In this way, we again obtain the discrete energy conservation 
\eq{C-N-2nd-conserv}. 


If also (\ref{ass-for-Phi}d,e) holds, 
two potentials \eqref{C-R-system-potential}
should be modified as
\begin{subequations}\label{C-R-system-potential+}\begin{align}
&\label{C-R-system-potential-u+}
\!\!\!\!u\mapsto
\frac1\tau\mathfrak{F}_{u_\tau^{k-1},z_\tau^{k-1}}(u)
+\Psi_1\Big(\frac{u{-}u_\tau^{k-1}\hspace*{-.5em}}\tau\hspace*{.5em}\Big)
+2\tau\mathscr{T}\Big(\frac{u{-}\tau v_\tau^{k-1}{-}u_\tau^{k-1}\!\!}{\tau^2}\ \Big)
-\Big\langle f_\tau^k,\frac u\tau\Big\rangle,\!\!
\\&\label{C-R-system-potential-z+}
\!\!\!\!z\mapsto
\frac1\tau\mathfrak{G}_{u_\tau^k,z_\tau^{k-1}}(z)
+\Psi_2\Big(\frac{z{-}z_\tau^{k-1}\hspace*{-.5em}}\tau\hspace*{.5em}\Big)
-\Big\langle g_\tau^k,\frac z\tau\Big\rangle\:.
\end{align}\end{subequations}
Existence of such potentials underlying the scheme \eqref{abstract-hyper-disc+}
can in particular cases facilitate numerical solution when appropriate
algorithms are used.

This scheme indeed generalizes Sect.~\ref{sec-nonconvex}. Actually, if 
$\Phi(\cdot,z)$ or $\Phi(u,\cdot)$ are quadratic, one can take simply 
\begin{align}\label{DPhi-if-quadratic}
\mathfrak{D}_{u}^{}\Phi(u,\tilde u,z)
:=\pl_z\Phi\Big(\frac{u{+}\tilde u}2,z\Big)
\quad\text{ or }\quad 
\mathfrak{D}_{z}^{}\Phi(u,z,\tilde z)
:=\pl_z\Phi\Big(u,\frac{z{+}\tilde z}2\Big)
\end{align} 
so that \eqref{abstract-hyper-disc-b+} or \eqref{abstract-hyper-disc-c+} 
coincide with \eqref{abstract-hyper-disc-c}, respectively. 
Then also $\mathfrak{F}_{\tilde u,z}$ and $\mathfrak{G}_{u,\tilde z}$ 
from (\ref{ass-for-Phi}d,e) exist and 
can be taken as
\begin{align}
\mathfrak{F}_{\tilde u,z}(u)=2\Phi\Big(\frac{u{+}\tilde u}2,z\Big)
\qquad\text{ and }\qquad 
\mathfrak{G}_{u,\tilde z}(z)=2\Phi\Big(u,\frac{z{+}\tilde z}2\Big).
\end{align} 
In a general 
non-quadratic $\Phi(\cdot,z)$ or $\Phi(u,\cdot)$,
existence of $\mathfrak{D}_{u}^{}\Phi$ or $\mathfrak{D}_{z}^{}\Phi$ satisfying 
\eqref{ass-for-Phi} however relies rather on a special local scalar character of
involved operators, cf.\ the example 
\eqref{quotient} below.
}

\TINY{
\begin{remark}[{\it Non-quadratic $\Phi(u,\cdot)$: 
combination with a backward-Euler scheme.}]\label{rem-nonquadratic}
\upshape
Some\-times, even the assumption \eqref{structural-ass-b}\COMMENT{LINK TO BELOW!!!} is unrealistically
strong and we must still relax by considering only   \COMMENT{BETTER ALREADY IN REM.~\ref{rem-nonquadratic-} -- BUT CONVEXITY NOT NEEDED THERE!}
\begin{align}\label{structural-ass-b-weak}
\forall u\in K_1:
\quad\Phi(u,\cdot):Z\to\R\cup\{\infty\}\ \text{ convex, lower semicontinuous}.
\end{align}
Then the binomial formula \eqref{2x-binomial-z} does not work any more
and, instead of using the fractional-step strategy to decoupling 
the system to two Crank-Nicolson formulas \eqref{abstract-hyper-disc},
we now use it to combine one Crank-Nicolson scheme with the backward-Euler formula.
More specifically, we now modify \eqref{abstract-hyper-disc} in the way
that $z_\tau^{k-1/2}$ in \eqref{abstract-hyper-disc-c} is replaced by $z_\tau^k$
\begin{align}\label{abstract-hyper-disc-Euler}
&&\pl\Psi_2\Big(\frac{z_\tau^k{-}z_\tau^{k-1}\hspace*{-.5em}}\tau\hspace*{.5em}\Big)
+\pl^{}_z\Phi(u_\tau^k,z_\tau^k)\ni g_\tau^k\ \ \ \ \text{ for }k=1,2...,
&&z_\tau^0=z_0,&&&&
\end{align}
while (\ref{abstract-hyper-disc}a,b) is kept unchanged. This problem has, 
instead of
\eqref{C-R-system-potential-z}, the convex potential:
\begin{align}\label{backward-Euler-potential}
z\mapsto\frac{\Phi(u_\tau^k,z)}\tau
+\Psi_2\Big(\frac{z{-}z_\tau^{k-1}\hspace*{-.5em}}\tau\hspace*{.5em}\Big)
-\Big\langle g_\tau^k,\frac z\tau\Big\rangle\:.
\end{align}
%
%
Its coercivity and lower semicontinuity inherited from \eq{structural-ass-b}
with \eq{struct-ass-b} ensures existence of some $z_\tau^k$ solving 
\eq{abstract-hyper-disc-Euler}.
To reveal the energetics, we use the same test as we made for 
\eqref{C-N-2nd}, namely $\vhalf$ for \eqref{abstract-hyper-disc-b}
and $\frac{z_\tau^k{-}z_\tau^{k-1}}\tau$ for \eqref{abstract-hyper-disc-c}.
Using \eqref{structural-ass},\COMMENT{CHECK LABEL}
it gives now only one equality but another one is only an upper estimate,
namely
\begin{subequations}\label{C-N-2nd-system+}\begin{align}\nonumber
&\frac{\mathscr{T}(v_\tau^k)-\mathscr{T}(v_\tau^{k-1})}\tau
+\Xi_1\big(\vhalf\big)
\\[-.2em]&\hspace*{7.5em}
+\frac{\Phi(u_\tau^k,z_\tau^{k-1})-\Phi(u_\tau^{k-1},z_\tau^{k-1})}\tau
=\big\langle f_\tau^k,\vhalf\big\rangle,
\label{C-N-2nd-system-a}
\\&
\Xi_2\Big(\frac{z_\tau^k{-}z_\tau^{k-1}\hspace*{-.5em}}\tau\hspace*{.5em}\Big)
+\frac{\Phi(u_\tau^k,z_\tau^k)-\Phi(u_\tau^k,z_\tau^{k-1})}\tau
\le\Big\langle g_\tau^k,\frac{z_\tau^k{-}z_\tau^{k-1}\!\!}\tau\ \Big\rangle.
\label{C-N-2nd-system-b}
\end{align}\end{subequations}
The inequality in \eqref{C-N-2nd-system-b} relies
on the convexity \eq{structural-ass-b-weak}.
Like previously in \eqref{C-N-2nd-system},
we again benefit from the 
cancellation of the terms $\pm\Phi(u_\tau^k,z_\tau^{k-1})$ when summing 
\eqref{C-N-2nd-system+} and obtain the 
energy-type estimate like \eqref{C-N-2nd-conserv} but now
only as an inequality:
 \begin{align}\nonumber
&\frac{\mathscr{T}(v_\tau^k)-\mathscr{T}(v_\tau^{k-1})}\tau
+\Xi_1\big(\vhalf\big)
+\Xi_2\Big(\frac{z_\tau^k{-}z_\tau^{k-1}\hspace*{-.5em}}\tau\hspace*{.5em}\Big)
\\&\hspace*{8em}
+\frac{\Phi(u_\tau^k,z_\tau^k)-\Phi(u_\tau^{k-1},z_\tau^{k-1})}\tau
\le\big\langle f_\tau^k,\vhalf\big\rangle
+\Big\langle g_\tau^k,\frac{z_\tau^k{-}z_\tau^{k-1}\!\!}\tau\ \Big\rangle.
\label{C-N-2nd-system+}\end{align}
This is a certain reasonable compromise and a certain price one must pay for 
a weakening of \eqref{structural-ass-b} towards \eqref{structural-ass-b-weak}:
one has always a guaranteed upper estimate which makes the decoupled 
scheme (\ref{abstract-hyper-disc}a,b)--\eqref{abstract-hyper-disc-Euler} 
numerical stable and, in qualified cases, convergent and, at least if the 
``inelastic'' variable $z$ does not evolve, the energy is still conserved 
even in the discrete case.
Actually, if $\Phi$ is non-quadratic in terms of $u$, one 
can use the Crank-Nicolson formula for \eqref{abstract-2nd-v}
while the backward-Euler formula is used already for \eqref{abstract-2nd-u},
cf.\ \cite{CarDol04TSDN}.
\end{remark}
}

\TINY{
\begin{remark}[{\it Set constraints.}]
\upshape
A combination of the Crank-Nicolson scheme with the fractional-step
decoupling together with the backward-Euler discretisation of possible 
constraints can be made like in Remark~\ref{rem-nonquadratic-} would allow
for relaxation of the qualification \eqref{structural-ass} as
\begin{subequations}\label{structural-ass+}\begin{align}
&\Phi=\Phi_0+\delta_K
\ \ \text{ coercive
with 
$K\subset U{\times}Z$ convex, closed,}
\label{structural-ass-0}
\\\label{structural-ass-a}
&\forall z\in Z:\quad\Phi_0(\cdot,z):U\to\R\ \text{ quadratic, convex, 
}
\\\label{structural-ass-b}
&\forall u\in U:\quad\Phi_0(u,\cdot):Z\to\R\ \text{ quadratic, convex. 
}
\end{align}\end{subequations}
Of course, as in Remark~\ref{rem-nonquadratic-}, only upper energy 
estimate is at disposal. \COMMENT{HERE I WILL STILL CHECK
WHETHER THE DECOUPLING FOR $I_1^*N_K$ and $I_2^*N_K$ WORKS?
..... $\delta_K(u,z_\tau^{k-1})$ and $\delta_K(u_\tau^k,z)$ to be added
in \eqref{C-R-system-potential-u} and \eqref{C-R-system-potential-z},
respectively}
\end{remark}
}

\TINY{
\begin{remark}[{\it Weakening convexity of $\Phi(u,\cdot)$: semiconvexity.}]
\label{rem-semi-convex}
\upshape
If $\Psi_2$ is uniformly convex, then \eqref{structural-ass-b-weak} 
can still be weakened by requiring only a semi-convex of $\Phi(u,\cdot)$ 
uniform in $u$, i.e.\ $\Phi(u,\cdot)+L\|\cdot\|^2$ convex for enough 
big $L$. Then, for $\tau>0$ small enough, the functional 
\eqref{backward-Euler-potential} still remains convex and 
the estimate \eqref{C-N-2nd-system+} holds but only with a 
factor in $\Psi_2$ which only asymptotically approaches 1 from below, see 
e.g.\ \cite[Rem.\,8.24]{Roub13NPDE} or also 
\cite[Estimate~(5.1.29)]{MieRou15RIST}.
\end{remark}
}

\TINY{
\begin{remark}[{\it Nonconvex $\Phi(u,\cdot)$.}]\label{rem-nonconvex}
\upshape
Some application needs $\Phi(u,\cdot)$ weakly lower-semicontinuous but 
nonconvex in a way that 
Remarks~\ref{rem-semi-convex} or \ref{rem-more-steps} do not apply.
One can anyway rely on \eqref{backward-Euler-potential} which is then 
a global minimization problem. When $\Psi_2$ is 1-homogeneous, 
a rich theory was developed in \cite[Chap.\,5]{MieRou15RIST} 
by using the backward-Euler formula also for \eqref{abstract-hyper-system-a}.
It can be modified when combined with our Crank-Nicolson scheme for 
(\ref{abstract-2nd}a,b).
The upper-energy estimate is an obvious combination of 
the technique from the theory of rate-independent processes (comparing
values of the incremental minimization at two subsequent time levels
with using the 1-homogeneity of $\Psi_2$ and the cancellation effect 
when combined with the Crank-Nicolson formula). 
\end{remark}
}

\TINY{
\begin{remark}[{\it Failure of a lower energy estimate.}]\label{rem-lower-est}
\upshape
If $\Psi_2$ is 1-homogeneous and \eqref{structural-ass-b-weak} holds, 
it is well known that the backward Euler scheme (which then realizes 
the global minimizing of the potential \eq{backward-Euler-potential} of the
incremental formula \eqref{abstract-hyper-disc-Euler})
allows for the two-sided energy estimate, cf.\ \cite{MieRou15RIST}. 
Also, it is known (and we showed 
also here) that the Crank-Nicolson formula in situations of the 
quadratic stored energy \eqref{abstract-hyper-system-a} provides energy
equality. One can then hope that merging these two formulas
\eqref{abstract-hyper-disc} might lead not only to the upper estimate 
\eqref{C-N-2nd-system+} but also to a lower one. Unfortunately, 
this is not true. In special cases, one can at least get a 
locally-asymptotic lower estimate. Let us demonstrate it
under the qualification \eqref{struct-ass-a}. Using that $\Psi_2$ is 
now 1-homogeneous and $\Phi(u,\cdot)$ is convex, 
minimization of \eq{C-R-system-potential-z} yields the semistability 
$\Phi(u_\tau^k,z_\tau^k)\le\Phi(u_\tau^k,\tilde z)
+\Psi_2(\tilde z-z_\tau^k)+\langle h_\tau^k,\tilde z-z_\tau^k\rangle$
for any $\tilde z$ which, written for $k-1$ and tested by 
$\tilde z= z_\tau^k$, gives after some simple manipulation 
\begin{align}\nonumber
\Psi_2\Big(\frac{z_\tau^k{-}z_\tau^{k-1}\hspace*{-.5em}}\tau\hspace*{.5em}\Big)
+\frac{\Phi(u_\tau^{k-1},z_\tau^k)-\Phi(u_\tau^{k-1},z_\tau^{k-1})}\tau
\ge\Big\langle g_\tau^{k-1},\frac{z_\tau^k{-}z_\tau^{k-1}\!\!}\tau\ \Big\rangle.
\end{align}
By adding $\frac1\tau\Phi(u_\tau^k,z_\tau^k)-\frac1\tau\Phi(u_\tau^k,z_\tau^{k-1})$ 
to both sides, we obtain 
\begin{align}\nonumber
\Psi_2\Big(\frac{z_\tau^k{-}z_\tau^{k-1}\hspace*{-.5em}}\tau\hspace*{.5em}\Big)
&+\frac{\Phi(u_\tau^k,z_\tau^k)-\Phi(u_\tau^k,z_\tau^{k-1})}\tau
\ge\Big\langle g_\tau^{k-1},\frac{z_\tau^k{-}z_\tau^{k-1}\!\!}\tau\ \Big\rangle
+r_\tau^k
\\&\hspace*{0em}
\text{ where }\ r_\tau^k:=
\frac{\Phi(u_\tau^k,z_\tau^k)-\Phi(u_\tau^{k-1},z_\tau^k)}\tau
-\frac{\Phi(u_\tau^k,z_\tau^{k-1})-\Phi(u_\tau^{k-1},z_\tau^{k-1})}\tau
\end{align}
and, when sum it with the equality \eq{C-N-2nd-system-a}
we again enjoy the cancellation of the terms $\pm\Phi(u_\tau^k,z_\tau^{k-1})$
and now obtain the lower energy-type estimate
\begin{align}\nonumber
&\frac{\mathscr{T}(v_\tau^k)-\mathscr{T}(v_\tau^{k-1})}\tau
+\Xi_1\big(\vhalf\big)
+\Psi_2\Big(\frac{z_\tau^k{-}z_\tau^{k-1}\hspace*{-.5em}}\tau\hspace*{.5em}\Big)
\\&\hspace*{5em}\nonumber
+\frac{\Phi(u_\tau^k,z_\tau^k)-\Phi(u_\tau^{k-1},z_\tau^{k-1})}\tau
\ge\big\langle f_\tau^k,\vhalf\big\rangle+\Big\langle g_\tau^{k-1},\frac{z_\tau^k{-}z_\tau^{k-1}\!\!}\tau\ \Big\rangle+r_\tau^k.
\label{C-N-2nd-system+lower}\end{align}
Using again the notation $\langle\cdot,\cdot\rangle_\Phi$ for
the bilinear form related to the quadratic functional $\Phi$ in the sense
$\Phi(u,z)=\frac12\langle(u,z),(u,z)\rangle_\Phi$, we can write the reminder 
more explicitly as
$r_\tau^k=\langle(u_\tau^k{-}u_\tau^{k-1},0),(0,z_\tau^k{-}z_\tau^{k-1})\rangle_\Phi/\tau$. Obviously, it has no specific sign and only allows for an estimate
\begin{align}
\big|r_\tau^k\big|\le2\tau C_\Phi^{}\Big|\frac{u_\tau^k{-}u_\tau^{k-1}\hspace*{-.5em}}\tau\hspace*{.5em}\Big|
\Big|\frac{z_\tau^k{-}z_\tau^{k-1}\hspace*{-.5em}}\tau\hspace*{.5em}\Big|
\end{align}
with a constant $C_\Phi:=\sup_{\|(u,z)\|\le1}\Phi(u,z)$. In a finite-dimensional
case, one can assume $|\DT z_\tau|$ bounded in $L^1(0,T)$ while
$|\DT u_\tau|$ bounded in $L^\infty(0,T)$ due to the kinetic term.
Thus $|r_\tau^k|=\mathscr{O}(\tau)$. Yet, these argumentation breaks in 
infinite-dimensional case where the kinetic term usually does not 
control $u$ in the norm needed for $\Phi$. On top of it, 
the estimate \eqref{C-N-2nd-system+lower} is only local and,
when sum it up for time-levels $k=1,2,...$ whose number is 
$\mathscr{O}(1/\tau)$, the overall reminder in the lower-energy balance
does not converge to 0 for $\tau\to0$ even in the finite-dimensional case.
\end{remark}
}

\begin{remark}[{\it More fractional steps.}]\label{rem-more-steps}
\upshape
The fractional-step splitting can easily be iterated and the canceling 
effect is then ``telescopic'', cf.\ \cite[Rem.\,8.25]{Roub13NPDE}. Thus 
some other equations can be easily added. In particular, the heat-transfer 
problem can be considered to get a 3-step decoupled scheme, 
cf.\ \cite{Roub??ECTD}.
\end{remark}


\vspace*{1em}
\section{\sf CONVERGENCE ANALYSIS IN PARTICULAR CASES}\label{sec-conv}

Beside the standard notation for the Lebesgue $L^p$-spaces, we will use 
$W^{k,p}$ for Sobolev spaces whose $k$-th derivatives 
are in $L^p$-spaces. We abbreviate $H^k=W^{k,2}$. We consider a 
fixed time interval $I=[0,T]$ and, for a Banach space $X$, we denote
by $L^p(I;X)$ the standard Bochner space of Bochner-measurable mappings 
$I\to X$ with whose norm in $X$ is in $L^p(I)$. 
Also, $W^{k,p}(I;X)$ denotes the Banach space of mappings from $L^p(I;X)$ 
whose $k$-th distributional derivative in time is also in $L^p(I;X)$.
Also, $C(I;X)$ and $C_\text{weak}(I;X)$ will denote the Banach space of 
continuous and weakly continuous mappings $I\to X$, respectively.
Moreover, we denote by ${\rm BV}(I;X)$ the Banach space
of the mappings $I\to X$ that have 
bounded variation on $I$, and by ${\rm B}(I;X)$ the space of 
Bochner measurable, everywhere defined, and bounded mappings $I\to X$.
By ${\rm Meas}(I;X)$ we denote the space of $X$-valued measures on $I$.

First, we define suitably a weak solution to \eq{abstract-2nd}.
Assuming $\Phi$ smooth and using just the definition of the convex
subdifferentials $\partial\Psi_1$ and $\partial\Psi_2$, the inclusions in 
(\ref{abstract-2nd}b,c) read as 
\begin{subequations}\label{abstract-2nd++}\begin{align}
\label{abstract-2nd-u++}
&\forall\,\tilde v\in U:\qquad\Psi_1(v)\le\Psi_1(\tilde v)
+\big\langle f-\pl^{}_u\Phi(u,z)-\mathscr{T}'\DT v\,,\,\tilde v-v\big\rangle,
\\\label{abstract-2nd-b++}
&\forall\,\tilde z\in Z:\qquad
\Psi_2(\DT z)\le\Psi_2(\tilde z)
+\big\langle g-\pl^{}_z\Phi(u,z)\,,\,\tilde z-\DT z\big\rangle\,.
\end{align}\end{subequations}
Summing them and integrating over the time interval $I=[0,T]$, 
using \eq{abstract-2nd-v} while making also
the calculus $\langle\mathscr{T}'\DT v,\DT u\rangle
=\langle\mathscr{T}'\DT v,v\rangle=
\langle\mathscr{T}'v,\DT v\rangle=\frac{\d}{\d t}\mathscr{T}(v)$
and $\langle\partial_u\Phi(u,z),\DT u\rangle
+\langle\partial_z\Phi(u,z),\DT z\rangle
=\frac{\d}{\d t}\Phi(u,z)$, we eventually obtain:

\begin{definition}[Weak solution.]\label{def-weak-sln}
We call the triple $u\in W^{1,p_1}(I;U)$, 
$v\in L^{p_1}(I;U)\cap C_{\rm weak}(I;H)\cap W^{1,1}(I;U^*)$, and 
$z\in W^{1,p_2}(I;Z)$ a weak solution to \eq{abstract-2nd} if 
\begin{subequations}\label{weak-form}\begin{align}\label{weak-form-a}
&\DT u=v\ \text{ in the sense of distributions, and }
\\\nonumber
&\!\!\int_0^T\!\!
\Psi_1(\tilde v)+\Psi_2(\tilde z)
+\big\langle
\mathscr{T}'\DT v{+}\partial_u\Phi(u,z){-}f,\tilde v\big\rangle
+\big\langle\partial_z\Phi(u,z){-}g,\tilde z\big\rangle
+\big\langle f,\DT u\big\rangle
+\big\langle g,\DT z\big\rangle
\,\d t
\\[-.3em]&\hspace{4em}
+\mathscr{T}(v_0)+\Phi(u_0,z_0)\ge\mathscr{T}(v(T))+\Phi(u(T),z(T))
+\int_0^T\!\!\Psi_1(\DT u)+\Psi_2(\DT z)\,\d t
\label{weak-form-b}
\end{align}\end{subequations}
holds for all $\tilde v\in L^\infty(I;U)$ and $\tilde z\in L^\infty(I;Z)$.
\end{definition}

Considering a fixed time step $\tau>0$ as in the previous sections such that $T/\tau$ is integer, and $\{u_\tau^k\}_{k=0,...,K}$ with $K=T/\tau$,
we introduce a notation for the piecewise-constant and the piecewise affine 
interpolants defined respectively by
\begin{subequations}\label{def-of-interpolants}
\begin{align}\label{def-of-interpolants-}
&&&
\overline{u}_\tau(t)= u_\tau^k,\qquad\ \
\underline u_\tau(t)= u_\tau^{k-1},\qquad\ \
\underline{\overline u}_\tau(t)=
u_\tau^{k-1/2},
&&\text{and}
&&
\\&&&\label{def-of-interpolants+}
u_\tau(t)=\frac{t-(k{-}1)\tau}\tau u_\tau^k
+\frac{k\tau-t}\tau u_\tau^{k-1}
&&\hspace*{-7em}\text{for }(k{-}1)\tau<t\le k\tau.
\end{align}\end{subequations}
Similar meaning has also $v_\tau$, $\overline v_\tau$, etc.

First, let us investigate the situation in Section~\ref{sec-quadratic},
relying on the structural assumption \eqref{struct-ass}
except that the $p$-homogeneity will not be particularly exploited.
In terms of the above introduced interpolants, one can write 
the scheme \eq{C-N-2nd} analogously to \eq{abstract-2nd} as
\begin{subequations}\label{C-N-2nd+interp}
\begin{align}\label{C-N-2nd+a}
&\DT u_\tau=\underline{\overline v}_\tau,&&u_\tau\big|^{}_{t=0}\!=u_0,&&
\\&\mathscr{T}'\DT v_\tau+\partial\Psi_1(\DT u_\tau)
+\partial_u\Phi(\underline{\overline u}_\tau,\underline{\overline z}_\tau)
=\overf_\tau,&&v_\tau\big|^{}_{t=0}\!=v_0,
\\&\partial\Psi_2(\DT z_\tau)
+\partial_z\Phi(\underline{\overline u}_\tau,\underline{\overline z}_\tau)=\overline g_\tau,
&&z_\tau\big|^{}_{t=0}\!=z_0,
\end{align}\end{subequations}
considered on the time interval $I=[0,T]$. Assuming $\Psi_i$ are convex, 
cf.\ \eqref{struct-ass-b}, and $\Phi$ smooth, the definition of the convex
subdifferentials $\partial\Psi_1$ and $\partial\Psi_2$ allows for writing 
\eqref{C-N-2nd+interp} in sum as the variational inequality 
\begin{align}\nonumber
\big\langle\mathscr{T}'\DT v_\tau+\partial_u\Phi(\underline{\overline u}_\tau,\underline{\overline z}_\tau)-\overf_\tau,\tilde u-\DT u_\tau\big\rangle+\Psi_1(\tilde u)
&
\\+\,\big\langle\partial_z\Phi(\underline{\overline u}_\tau,\underline{\overline z}_\tau)-\overline g_\tau,\tilde z-\DT z_\tau\big\rangle+\Psi_2(\tilde z)&
\ge\Psi_1(\DT u_\tau)+\Psi_2(\DT z_\tau)
\label{discrete-VI-t}\end{align}
to hold for any $(\tilde u,\tilde z)\in U\times Z$ and for a.a.\ $t\in I$. 
We further use the binomial formulas 
\eqref{binomial} after summation for $k=1,...,T/\tau$ written in the form
\begin{align}\nonumber
&\int_0^T\!\!
\big\langle\mathscr{T}'\DT v_\tau+\partial_u\Phi(\underline{\overline u}_\tau,\underline{\overline z}_\tau)
,\DT u_\tau\big\rangle
+\big\langle\partial_z\Phi(\underline{\overline u}_\tau,\underline{\overline z}_\tau)
,\DT z_\tau\big\rangle\,\d t
\\[-.3em]&\qquad\qquad
=
\mathscr{T}(v_\tau(T))+\Phi(u_\tau(T),z_\tau(T))
-\mathscr{T}(v_0)-\Phi(u_0,z_0).
\label{C-N-2nd-conserv+}\end{align}
Substituting it into \eq{discrete-VI-t} integrated over $I$,
we obtain the discrete analog of \eq{weak-form-b}, namely
\begin{align}\nonumber
&\int_0^T\!\!
\bigg(
\Psi_1(\tilde v)+\Psi_2(\tilde z)
+\big\langle
\mathscr{T}'\DT v_\tau{+}\partial_u\Phi(\underline{\overline u}_\tau,\underline{\overline z}_\tau){-}\overf_\tau,\tilde v\big\rangle
+\big\langle\partial_z\Phi(\underline{\overline u}_\tau,\underline{\overline z}_\tau)
{-}\overline g_\tau,\tilde z\big\rangle
\\[-.3em]\nonumber&\hspace{4.5em}
+\big\langle\overf_\tau,\DT u_\tau\big\rangle
+\big\langle\overline g_\tau,\DT z_\tau\big\rangle
\bigg)
\,\d t
+\mathscr{T}(v_0)+\Phi(u_0,z_0)
\\[-.5em]&\hspace{10em}\ge\mathscr{T}(v_\tau(T))+\Phi(u_\tau(T),z_\tau(T))+\!\int_0^T\!\!\Psi_1(\DT u_\tau)+\Psi_2(\DT z_\tau)
\,\d t,
\label{weak-form-disc}
\end{align}
while an analog of \eqref{weak-form-a} is just \eqref{C-N-2nd+a}.

In terms of these interpolants, the discrete energy conservation 
\eqref{C-N-2nd-conserv} summed for $k=1,...,T/\tau$ can be written 
as:
\begin{align}\nonumber
&\mathscr{T}(v_\tau(T))+\Phi(u_\tau(T),z_\tau(T))+\int_0^T\!\!\Xi_1(\underline{\overline v}_\tau)+\Xi_2(\DT z_\tau)\,\d t
\\[-.5em]&\qquad\qquad\qquad\qquad\qquad\qquad\qquad\ =
\mathscr{T}(v_0)+\Phi(u_0,z_0)+\int_0^T\!\!\big\langle\overf_\tau,\DT u_\tau\big\rangle
+\big\langle\overline g_\tau,\DT z_\tau\big\rangle\,\d t\,.
\label{C-N-2nd-conserv++}\end{align}

Note however that neither $\overline v_\tau$ not $\underline v_\tau$ is  
the velocity corresponding to $u_\tau$, i.e.\ $\overline v_\tau\ne\DT u_\tau$
and $\underline v_\tau\ne\DT u_\tau$ in 
general, although in the limit both $\overline v_\tau-\DT u_\tau\to0$ and
$\underline v_\tau-\DT u_\tau\to0$ for $\tau\to0$. If 
$v_0\!\in\!H{\setminus}U$, then even both $\overline v_\tau$ and 
$\underline v_\tau$ are not valued in $U$, although 
$\underline{\overline v}_\tau\!\in\!L^\infty(I;U)$.
In parti\-cular, the relation $\DT u_\tau=\underline{\overline v}_\tau$ 
must be taken into account as a vital ingredient accompanied the
variational inequality \eq{C-N-2nd-conserv+}.


\begin{proposition}[Numerical stability of \eq{C-N-2nd} and 
convergence.]\label{prop-conv-1}
Let $U$ and $Z$ be reflexive Banach spaces, $U$ densely embedded into a 
a Hilbert space $H$, \eq{struct-ass} hold, 
and furthermore $\|\pl\Psi_1(\cdot)\|_{U^*}\le C(1+\|\cdot\|_U^{p_1-1})$ with some $C\in\R$, and 
let  $u_0\in U$, $v_0\in H$, $z_0\in Z$,
$f\in L^{p_1'}(I;U^*)+L^1(I;H)$, and $g\in L^{p_2'}(I;Z^*)$ with 
$p_i'=p_i^{}/(p_i^{}{-}1)$ with $p_i^{}>1$ from \eqref{struct-ass-c}.
Then the following a-priori estimates hold: 
\begin{subequations}\label{C-N-2nd-est}\begin{align}
\label{C-N-2nd-est1}
&\!\!\!\|u_\tau\|^{}_{L^\infty(I;U)\,\cap\,W^{1,\infty}(I;H)\,\cap\,W^{1,p_1}(I;U)}\!\le C,\!&&
\|z_\tau\|^{}_{W^{1,p_2}(I;Z)}\!\le C,\!\!
\\&\!\!\!\|\underline{\overline v}_\tau\|^{}_{L^\infty(I;H)\,\cap\,L^{p_1}(I;U)}\!
\le C,&& 
\|v_\tau\|^{}_{W^{1,\max(p_1',2)}(I;U^*)+W^{1,1}(I;H)}\le C\,.
\label{C-N-2nd-est2}
\end{align}\end{subequations}
Moreover, the sequence $(u_\tau,\underline{\overline v}_\tau,z_\tau)$ 
converges weakly* in the topologies
indicated by the first three estimates \eq{C-N-2nd-est} to 
the unique weak solution $(u,v,z)$ due to Definition~\ref{def-weak-sln}.
Moreover, if the Banach space $U$ and the potential $\Psi_1$ 
(resp.\ $Z$ and $\Psi_2$) are uniformly 
convex and if $p_1\le2$, we have even the strong convergences
\begin{subequations}\label{strong-conv}\begin{align}
&&&&&u_\tau\to u&&\text{ in }\ W^{1,p_1}(I;U),\ \ \text{ resp.}&&&&&&
\\&&&&&z_\tau\to z&&\text{ in }
W^{1,p_2}(I;Z).
\end{align}\end{subequations}
\end{proposition}

\noindent{\it Proof}.
The first three a-priori estimates in
\eq{C-N-2nd-est}
can be obtained by usual estimation from \eq{C-N-2nd-conserv}, using 
the Young and the discrete Gronwall inequalities.
Then, from \eq{C-N-2nd+interp} by comparison, we obtain also the last 
estimate in \eq{C-N-2nd-est2}.

Next, by the Banach selection principle, we choose a weakly* convergent 
subsequence with respect to the weak* topologies in the spaces 
indicated in the  first three a-priori estimates in
\eq{C-N-2nd-est}. Having relevant $L^\infty$-estimates at disposal together
with corresponding time derivatives estimated, we can
also rely on 
$$
u_\tau(T)\to u(T)\ \text{ in }\ U,\qquad
v_\tau(T)\to v(T)\ \text{ in }\ H,\qquad
z_\tau(T)\to z(T)\ \text{ in }\ Z\ \text{ weakly}.
$$
Then the limit passage from \eq{weak-form-disc} to \eq{weak-form}
is easy when exploiting the assumed convexity which makes the 
right-hand side of \eq{weak-form-disc} weakly lower-semicontinuous and 
when realizing the assumption \eqref{struct-ass-a}
which makes both $\pl_u\Phi$ and $\pl_z\Phi$ linear so that 
the right-hand side of \eq{weak-form-disc} is weakly continuous.
Moreover, passing to the limit in the relation 
$\DT u_\tau=\underline{\overline v}_\tau$ yields $\DT u=v$.

In fact, not only the selected subsequence but even the whole sequence converges
to the weak solution $(u,v,z)$ because it is unique.
To see this uniqueness, we subtract the equation (inequality) for two
solutions and test it by the difference of time derivatives of them. 
The monotonicity of $\pl\Psi_1$ and $\pl\Psi_2$ is then to be used 
together with linearity of $\Phi'$ and $\mathscr{T}'$;  cf.\ 
e.g.\ \cite[Prop.\,11.35]{Roub13NPDE}. 

The strong convergence \eq{strong-conv} can then be seen by
estimation:
\begin{align}\nonumber
&\!
\int_0^T\!\!\Xi_1(v)+\Xi_2(\DT z)\,\d t
\le\liminf_{\tau\to0}
\int_0^T\!\Xi_1(\underline{\overline v}_\tau)
+\Xi_2(\DT z_\tau)\,\d t
\le\limsup_{\tau\to0}
\int_0^T\!\Xi_1(\underline{\overline v}_\tau)
+\Xi_2(\DT z_\tau)\,\d t
\\&\nonumber\quad
=\lim_{\tau\to0}\bigg(\mathscr{T}(v_0)+\Phi(u_0,z_0)+\int_0^T\!\!
\big\langle\overf_\tau,\underline{\overline v}_\tau\big\rangle
+\big\langle\overline g_\tau,\DT z_\tau\big\rangle\,\d t
-\mathscr{T}(v_\tau(T))-\Phi(u_\tau(T),z_\tau(T))\bigg)
\\[-.3em]&\nonumber\quad\le\mathscr{T}(v_0)+\Phi(u_0,z_0)+\int_0^T
\big\langle f,v\big\rangle+\big\langle g,\DT z\big\rangle\,\d t
-\mathscr{T}(v(T))-\Phi(u(T),z(T))
\\[-.3em]&\quad=
\int_0^T\!\!\Xi_1(v)+\Xi_2(\DT z)\,\d t
\label{trick}\end{align}
where the first equality has used \eq{C-N-2nd-conserv++}
while the last equality is based on that \eq{abstract-2nd} is already proved 
together with the fact that $\DT v\in L^{\max(p_1',2)}(I;U^*)+L^1(I;H)$ 
is in duality with $v\in L^{p_1}(I;U)\cap L^\infty(I;H)$ 
so that the by-part integration of the $\mathscr{T}$-term can legally
be executed; here we needed $p_1\le2$.
Note that $\DT v$ indeed remains in the nonreflexive space 
$L^{\max(p_1',2)}(I;U^*)+L^1(I;H)$ by a comparison argument because $f-\pl\Psi_1(v)-\pl\Phi(u,z)$ is in 
this space. Similarly, we use the estimates on $\DT u$ and $\DT z$ for the 
by-part integration of the $\Phi$-term;
note that $\DT u\in L^\infty(I;H)\cap L^{p_1}(I;U)$ is certainly in duality
with $\partial_u\Phi(u,z)\in L^\infty(I;U^*)$
and $\DT z\in L^{p_2}(I;Z)$ in duality with 
$\partial_z\Phi(u,z)\in L^\infty(I;Z^*)$. Therefore, \eq{trick} implies that
\begin{align}
\lim_{\tau\to0}
\int_0^T\!\!\Xi_1(\underline{\overline v}_\tau)
+\Xi_2(\DT z_\tau)\,\d t
=\int_0^T\!\!\Xi_1(v)+\Xi_2(\DT z)\,\d t.
\end{align}
The assumed uniform convexity of 
$\Psi$'s
and of the underlying Banach spaces $U$ and $Z$ (and thus of $L^{p_1}(I;U)$ and $L^{p_2}(I;Z)$ too)
together with the weak convergence then yields the strong convergence \eqref{strong-conv}
by the Fan-Glicksberg theorem.
\QED

%
\MARGINOTE{references?? check D. Richtmyer and K. W. Morton}

Let us now come to the decoupled scheme from Section~\ref{sec-nonconvex},
relying on the weakened structural assumption \eqref{structural-ass}
instead of \eqref{struct-ass-a}. 
Having in mind the Yanenko-type 
time-discrete scheme \eqref{abstract-hyper-disc},
in terms of the interpolant-notation \eqref{def-of-interpolants}, 
the discrete variational inequality like \eqref{discrete-VI-t}
must reads as:
\begin{align}\nonumber
\big\langle\mathscr{T}'\DT v_\tau+\partial_u\Phi(\underline{\overline u}_\tau,\underline{z}_\tau)-\overf_\tau,\tilde v-\DT u_\tau\big\rangle+\Psi_1(\tilde v)
&
\\+\,\big\langle\partial_z\Phi(\overline u_\tau,\underline{\overline z}_\tau)
-\overline g_\tau,\tilde z-\DT z_\tau\big\rangle+\Psi_2(\tilde z)&
\ge\Psi_1(\DT u_\tau)+\Psi_2(\DT z_\tau).
\label{discrete-VI-t+}\end{align}
We further use the binomial formulas
\eqref{2x-binomial}
and the 
cancellation effect in \eq{C-N-2nd-system} and,
after summation for $k=1,...,T/\tau$,
we obtain the slightly modified equality \eq{C-N-2nd-conserv+},
namely
\begin{align}\nonumber
&\int_0^T\!\!
\big\langle\mathscr{T}'\DT v_\tau+\partial_u\Phi(\underline{\overline u}_\tau,\underline{z}_\tau)
,\DT u_\tau\big\rangle
+\big\langle\partial_z\Phi(\overline u_\tau,\underline{\overline z}_\tau)
,\DT z_\tau\big\rangle\,\d t
\\[-.3em]&\qquad\qquad\qquad
=
\mathscr{T}(v_\tau(T))+\Phi(u_\tau(T),z_\tau(T))
-\mathscr{T}(v_0)-\Phi(u_0,z_0).
\label{C-N-2nd-conserv+++}
\end{align}
Substituting it in into \eq{discrete-VI-t+} integrated over $I$, we obtain the 
discrete analog of \eq{weak-form} like \eq{weak-form-disc} but now modified as
\begin{align}\nonumber
&\int_0^T\!\!
\Psi_1(\tilde v)+\Psi_2(\tilde z)
+\big\langle
\mathscr{T}'\DT v_\tau{+}\partial_u\Phi(\underline{\overline u}_\tau,\underline{z}_\tau){-}\overf_\tau,\tilde v\big\rangle
+\big\langle\partial_z\Phi(\underline{\overline u}_\tau,\overline z_\tau)
{-}\overline g_\tau,\tilde z\big\rangle
\\[-.3em]\nonumber&\hspace{4.5em}
+\big\langle\overf_\tau,\DT u_\tau\big\rangle
+\big\langle\overline g_\tau,\DT z_\tau\big\rangle
\,\d t
+\mathscr{T}(v_0)+\Phi(u_0,z_0)
\\[-.3em]&\hspace{9em}\ge\mathscr{T}(v_\tau(T))+\Phi(u_\tau(T),z_\tau(T))+\!\int_0^T\!\!\Psi_1(\DT u_\tau)+\Psi_2(\DT z_\tau)
\,\d t.
\label{weak-form-disc+}
\end{align}
The a-priori estimates \eq{C-N-2nd-est1} can be derived as before because 
\eq{C-N-2nd-conserv} is at disposal, while the estimates \eq{C-N-2nd-est2}
can be obtained again by comparison from \eq{C-N-2nd+interp} but
modified by replacing 
$\pl_u\Phi(\underline{\overline u}_\tau,\underline{\overline z}_\tau)$ and 
$\pl_z\Phi(\underline{\overline u}_\tau,\underline{\overline z}_\tau)$ 
respectively by $\pl_u\Phi(\underline{\overline u}_\tau,\underline{z}_\tau)$ and 
$\pl_z\Phi({\overline u}_\tau,\underline{\overline z}_\tau)$. As these terms are 
now nonlinear even if $\Phi$ is smooth but only component-wise quadratic, the 
limit passage for $\tau\to0$ is more difficult in comparison with 
Proposition~\ref{prop-conv-1}. We impose general assumptions:\hspace*{-.1em}
\begin{subequations}\label{qualif}\begin{align}\label{qualif-a}
&\Phi,\Psi_1,\Psi_2\ \text{ are weakly lower-semicontinuous},\ \ 
\|\pl\Psi_1(\cdot)\|_{U^*}\le C(1+\|\cdot\|_U^{p_1-1}),
\\\label{qualif-b}&\partial_u\Phi\!:\!U\!\times\!Z\to U^*,\ 
\partial_z\Phi\!:\!U\!\times\!Z\to Z^*
\ \text{ are (weak,weak)-continuous}.
\end{align}\end{subequations}

\begin{proposition}[Numerical stability of \eqref{abstract-hyper-disc} and 
convergence.]
Let $U$, $H$, $Z$, $f$, and $g$ be qualified as in 
Proposition~\ref{prop-conv-1}, let
 $u_0\in U$, $v_0\in H$, $z_0\in Z$, let further the structural as\-sum\-p\-tions \eq{structural-ass} with 
{\rm(\ref{struct-ass}b,c)} and the qualification \eqref{qualif} hold.
Then all the statements of Proposition~\ref{prop-conv-1} hold.
\end{proposition}

\noindent{\it Sketch of the Proof}. 
The a-priori estimates \eq{C-N-2nd-est} work by the same way as 
in Proposition~\ref{prop-conv-1} by exploiting \eq{C-N-2nd-conserv}.
Then, after selection of weakly convergent subsequences, 
the limit passage of \eq{weak-form-disc+} towards \eq{weak-form-b} 
is immediate due to \eq{qualif}. Having proved that the limit is 
the weak solution, \eq{trick} works in an unchanged way and yields
the strong convergence \eq{strong-conv}.
%
\QED



An interesting special situations occur if some processes are
much faster than the external loading or the wave speed,
and can be well considered as arbitrarily fast and thus rate independent. 
Typically is concerns 
 the internal variable $z$, and 
then one is to consider $\Psi_2$ homogeneous of degree-1 and coercive on some 
Banach space $X\supset Z$, cf.\ also \cite[Sect.5.1-5.2]{MieRou15RIST}.
Then $z\in L^\infty(I;Z)\cap{\rm BV}(I;X)$ and Definition~\ref{def-weak-sln} 
is to be modified by replacing $\int_0^T\Psi_2(\DT z)\,\d t$ by the total 
variation 
\begin{align}\label{variation}
\mathrm{Diss}_{\Psi_2}(I;z)=\sup_{0\le 
t_0<t_1<...<t_N\le T,\ N\in\N}\sum_{i=1}^N
\Psi_2\big(z(t_i){-}z(t_{i-1})\big).
\end{align}
In particular, $\DT z$ is a $X$-valued measure in general.  
If $\DT z$ is regular, such modified definition holds also on 
any subinterval $[t_1,t_2]\subset I$, which is the concept of 
weak solution used which is, under mild qualification, equivalent to
a so-called local-solution concept used in the theory of rate-independent
processes, cf.\ \cite[Proposition 3.3.5]{MieRou15RIST}.
In case of the convexity of $\Phi$ as in Proposition~\ref{prop-conv-1} but
with $\Psi_2$ 1-homogeneous, even the uniqueness of the solution 
is again at disposal, cf.\ \cite[Prop.\,5.1.11]{MieRou15RIST}.

\TINY{
Let us illustrate it
on The time derivative of $z$
is then controlled only in the sense of measures and special 
mathematical techniques have to be used in particular if both $\Phi$ and 
$\Psi_2$ contain constraints and Definition~\ref{def-weak-sln} cannot
be simply used, cf.\ especially 
\cite{Miel05ERIS,MieRou15RIST,MieThe04RIHM}. In particular, assuming $\Psi_1$ 
smooth, Definition~\ref{def-weak-sln} modifies as follows, 
cf.\ also \cite[Def.\,5.1.1 \& Rem.\,5.1.8]{MieRou15RIST}:

\begin{definition}\label{def-ES}
We call $u\in W^{1,p_1}(I;U)$, 
$v\in L^{p_1}(I;U)\cap C_{\rm weak}(I;H)\cap W^{1,p_1'}(I;U^*)$ 
and 
$z\in\mathrm{BV}(I;Z)$  
an a.e.-local 
solution to \eq{abstract-2nd} if \eqref{weak-form-a} and \COMMENT{\eq{ES-b} BETTER AS AN EQUALITY??}
\begin{subequations}\label{ES}\begin{align}\nonumber
&\!\!\int_0^T\!\!
\big\langle
\mathscr{T}'\DT v{+}\Psi_1'(v)+\partial_u\Phi(u,z){-}f,\tilde u-u\big\rangle
-\big\langle\mathscr{T}'v,\DT{\tilde u}-\DT u\big\rangle
\,\d t
\\[-.3em]&\hspace{14em}
+\big\langle\mathscr{T}'v(T),\tilde u(T)-u(T)\big\rangle
\ge
\big\langle\mathscr{T}'v_0,\tilde u(0)-u_0\big\rangle
\label{ES-b},
\\\nonumber
&
\mathscr{T}(v(t_2))+\Phi(u(t_2),z(t_2))
+\int_{t_1}^{t_2}\!\!\Psi_1(\DT u)\,\d t
+\mathrm{Diss}_{\Psi_2}([t_1,t_2];z)-\big\langle g(t_2),z(t_2)\big\rangle
\\[-.7em]&\hspace{4em}
\le\mathscr{T}(v(t_1))+\Phi(u(t_1),z(t_1))-\big\langle g(t_1),z(t_1)\big\rangle
+\int_{t_1}^{t_2}\!\!\big\langle f,\DT u\big\rangle-\big\langle \DT g,z\big\rangle\,\d t,
\label{ES-c}
\\\label{ES-semistability}
&
\Phi(u(t),z(t))-\big\langle g(t),z(t)\big\rangle
\le\Phi(u(t),\tilde z)-\big\langle g(t),\tilde z\big\rangle
+\Psi_2(\tilde z{-}z(t))
\end{align}\end{subequations}
holds for all $\tilde u\!\in\!W^{1,\infty}(I;U)$,
for 
all $\tilde z\!\in\!Z$,
and for a.a.\ $t_1,t_2,t\!\in\![0,T]$, where 
\begin{align}
\mathrm{Diss}_{\Psi_2}([t_1,t_2];z)=\sup\sum_{i=1}^N
\Psi_2\big(z(t^{(i)}){-}z(t^{(i-1)})\big),
\end{align}
where the supremum is taken over all partitions of the type 
$t_1=t^{(0)}<t^{(1)}<...<t^{(N)}=t_2$ with $N\!\in\!\N$.
\end{definition}

In view of 
\eqref{C-R-system-potential-z}
and the 1-homogeneity of $\Psi_2$, $z_\tau^k$ minimizes the functional 
\begin{align}\nonumber
z&\mapsto\Psi_2(z{-}z_\tau^{k-1})+
2\Phi\Big(u_\tau^k,\frac{z{+}z_\tau^{k-1}}2\Big)
-\big\langle g_\tau^k,z\big\rangle.
\end{align}
Using 1-homogeneity of $\Psi_2$, we obtain the discrete semi-stability 
\begin{align}\label{disc-semi-stab}
\forall\tilde z\in Z:\quad2\Phi\Big(u_\tau^k,\frac{z_\tau^k{+}z_\tau^{k-1}}2\Big)
\le2\Phi\Big(u_\tau^k,\frac{\tilde z{+}z_\tau^{k-1}}2\Big)
+\Psi_2(\tilde z{-}z_\tau^k)-\big\langle g_\tau^k,\tilde z{-}z_\tau^k\big\rangle.
\end{align}
This can be written in terms of the interpolants \eq{def-of-interpolants} as
\begin{align}\label{disc-semi-stab+}
\forall_{\rm a.a.}t\!\in\!I\ \forall\tilde z\!\in\!Z{:}
\ \ 2\Phi\big(\overline u_\tau(t),\overline{\underline z}_\tau(t)\big)
\le2\Phi\Big(\overline u_\tau(t),\frac{\tilde z{+}{\underline z}_\tau(t)}2\Big)
+\Psi_2\big(\tilde z{-}\overline z_\tau(t)\big)
-\big\langle\overline g_\tau(t),\tilde z{-}\overline z_\tau(t)\big\rangle.
\end{align}


The difficult parts of the limit passage towards such solutions are the terms
$\mathscr{T}(v(t_1))+\Phi(u(t_1),z(t_1))$ in 
\eq{ES-c} which requires strong convergence the approximate solutions and 
the semistability \eq{ES-semistability} which requires an explicit construction of
a so-called mutual recovery sequence \cite{MiRoSt08GLRR}. A lot of examples of such 
constructions are known for the backward-Euler formula \eq{abstract-hyper-disc-Euler}, 
cf.\ e.g.\ \cite{MieRou15RIST}. For the Crank-Nicolson formula \eq{abstract-hyper-disc-c},
it deserves a suitable modification of the construction itself or 
of the proof because \eq{disc-semi-stab+} is not exactly in the form of
\eq{ES-semistability}. Let us demonstrate it on the case when $\Phi(u,\cdot)$ is quadratic.

Denoting the bilinear form underlying $\Phi(u,\cdot)$ by 
$\langle\cdot,\cdot\rangle_{\Phi,u}:Z\times Z\to\R$ the bilinear form defined by 
$$
\langle z,\tilde z\rangle_{\Phi,u}:=\big\langle\pl_z\Phi(u)z,\tilde z\big\rangle
$$
where the last duality is between $Z^*$ and $Z$, we have assumed the following joint continuity 
which is to be proved in particular cases:
\begin{align}\label{joint-cont}
\forall \tilde z\in Z:\qquad(u,z)\mapsto\big\langle z,\tilde z\big\rangle_{\Phi,u}
\ \text{ is (strong$\times$weak)-continuous}.
\end{align}

\begin{proposition}\label{prop-RI-quadratic}
Let ......, $\pl_u\Phi:U\times Z\to U^*$ continuous with 
at most linear growth,\COMMENT{CAN BE WEAKENED} 
$\Psi_1$ be quadratic and coercive on $U$, 
\end{proposition}

\noindent{\it Proof}. ....................
The coercivity of $\Psi_1$ allows us to show the 
strong convergence $\overline{u}_\tau\to u$ in $L^2(I;U)$ and thus
also $\overline{u}_\tau(t)\to u(t)$ in $U$ for a.a.\ $t\in I$, cf.\ 
\cite[Sect.\,5.1]{MieRou15RIST}.  

By Helly's selection principle, we can select a subsequence such that 
$\overline{z}_\tau(t)\to z(t)$ weakly if $\tau\to0$ 
even for any $t\in I$.
Simultaneously, we have also $\underline{z}_\tau(t)\to z(t)$ 
weakly for a.a.\ $t$, in particular for all points of continuity of $z$;
recall that functions with bounded variations are continuous with 
the exception at most countable number of time instances. We now take $t$ fixed.
Exploiting that $\Phi$ is quadratic and choosing $\tilde z$ arbitrary, 
we use the binomial trick with
the so-called mutual recovery sequence $\tilde z_\tau=\overline{z}_\tau(t)
+\tilde z-z(t)$, cf.\ \cite{MieRou15RIST} for details about this technique for
the backward Euler formula. The main motivation of this choice is 
to make $\Psi_2(\tilde z_\tau{-}\overline{z}_\tau(t))=\Psi_2(\tilde z{-}z(t))$
simply constant while $\tilde z_\tau\to\tilde z$ weakly, which makes the limit 
passage through the nonlinearities in \eqref{disc-semi-stab+} possible. Here, 
using elementary algebra,
this choice gives
\begin{align}\nonumber
&\Phi\Big(\overline{u}_\tau(t),\frac{\overline{z}_\tau(t){+}\underline{z}_\tau(t)}2\Big)
-\Phi\Big(\overline{u}_\tau(t),\frac{\tilde z_\tau{+}\underline{z}_\tau(t)}2\Big)
\\&\qquad\qquad\nonumber=\Phi\Big(\overline{u}_\tau(t),\frac{\overline{z}_\tau(t)}2\Big)
+\frac14\big\langle \overline{z}_\tau(t),\underline{z}_\tau(t)\big\rangle_{\Phi,\overline{u}_\tau(t)}
-\Phi\Big(\overline{u}_\tau(t),\frac{\tilde z_\tau}2\Big)
-\frac14\big\langle\tilde z_\tau,\underline{z}_\tau(t)\big\rangle_{\Phi,\overline{u}_\tau(t)}
\\&\qquad\qquad\nonumber
=\frac18\big\langle \overline{z}_\tau(t){+}\tilde z_\tau,\overline{z}_\tau(t){-}\tilde z_\tau\big\rangle_{\Phi,\overline{u}_\tau(t)}
+\frac14\big\langle \overline{z}_\tau(t){-}\tilde z_\tau,\underline{z}_\tau(t)\big\rangle_{\Phi,\overline{u}_\tau(t)}
\\&\qquad\qquad\nonumber
=\frac18\big\langle \overline{z}_\tau(t){+}\tilde z_\tau,z(t){-}\tilde z\big\rangle_{\Phi,\overline{u}_\tau(t)}
+\frac14\big\langle z(t){-}\tilde z,\underline{z}_\tau(t)\big\rangle_{\Phi,\overline{u}_\tau(t)}
\\&\qquad\qquad\nonumber
\to\frac18\big\langle z(t){+}\tilde z,z(t){-}\tilde z\big\rangle_{\Phi,u(t)}
+\frac14\big\langle z(t){-}\tilde z,z(t)\big\rangle_{\Phi,u(t)}
\\&\qquad\qquad
=\Phi(u(t),z(t))-\Phi\Big(u(t),\frac{\tilde z{+}z(t)}2\Big).
\label{binomial-MRS}
\end{align}
For the convergence, we used $\overline{u}_\tau(t)\to u(t)$ strongly in $U$ 
with the assumption \eq{joint-cont}.
Then we can easily perform the limit passage in \eqref{disc-semi-stab},
obtaining 
\begin{align}\label{semi-stab}
\forall\tilde z\in Z:\quad\Phi(u(t),z(t))\le\Phi\Big(u(t),\frac{\tilde z{+}z(t)}2\Big)
+\Psi_2\Big(\frac{\tilde z{-}z(t)}2\Big)
-\Big\langle g(t),\frac{\tilde z{-}z(t)}2\Big\rangle.
\end{align}
Now substituting $\tilde z=2\hat z-z(t)$ so that
$\frac{\tilde z{+}z(t)}2=\hat z$ and $\frac{\tilde z{-}z(t)}2=\hat z-z(t)$,
we obtain
the desired 
semistability \eqref{ES-semistability} only with $\hat z$ instead of 
$\tilde z$ .
\QED

}


\TINY{\hrule
{The convergence in
the rate-independent part then relies on an explicit construction of a
so-called mutual-recovery sequence \cite{MiRoSt08GLRR} in particular cases, see
also \cite[Chap.\,5]{MieRou15RIST}, but here modified .........}
\hrule}

\TINY{
\begin{remark}[{\it Nonquadratic $\Phi(u,\cdot)$.}]
\upshape
Let us now proceed to the weaker qualification \eq{structural-ass-b-weak}.
The decoupled Crank-Nicolson/backward-Euler scheme 
(\ref{abstract-hyper-disc}a,b)--\eqref{abstract-hyper-disc-Euler}
gives, by using the convexity of $\Phi(u,\cdot)$,
%
%
instead of the equality \eq{C-N-2nd-conserv+}
or \eq{C-N-2nd-conserv++} OR \eq{C-N-2nd-conserv+++}?????????????????????, the inequality
\begin{align}\nonumber
&\int_0^T\!\!
\big\langle\mathscr{T}'\DT v_\tau+\partial_u\Phi(\underline{\overline u}_\tau,\underline{z}_\tau)
,\DT u_\tau\big\rangle
+\big\langle\partial_z\Phi(\overline u_\tau,\overline z_\tau)
,\DT z_\tau\big\rangle\,\d t
\\[-.3em]&\qquad\qquad\qquad
\ge
\mathscr{T}(v_\tau(T))+\Phi(u_\tau(T),z_\tau(T))
-\mathscr{T}(v_0)-\Phi(u_0,z_0).
\label{C-N-2nd-conserv++++}\end{align}
Substituting it into \eq{discrete-VI-t} integrated over $I$,
we obtain \eq{weak-form-disc+} but modified by replacing 
$\partial_z\Phi(\underline{\overline u}_\tau,\overline z_\tau)$
by $\partial_z\Phi(\overline u_\tau,\overline z_\tau)$.
%
The a-priori estimates \eqref{C-N-2nd-est} still holds, but
the limit passage in \eqref{weak-form-disc+} towards \eqref{weak-form}
is now more involved because $\partial_u\Phi$ and $\partial_z\Phi$ are no 
longer linear. 
In the rate-independent cases like we mentioned above, the standard 
constructions of mutual recovery sequences can be used essentially 
by the same way as in \cite[Chap.\,5]{MieRou15RIST} where the backward-Euler 
formula was used also for \eq{abstract-hyper-system-a}.
This technique allows for a generalization for $\Phi$ non-quadratic and even 
nonconvex in $z$. 
\end{remark}
}

\COLOR{
\begin{remark}[{\it Convergence of the scheme from 
Sect.\,\ref{sec-nonquadratic}.}]\label{rem-conv-modified}
\upshape
The convergence of the modified scheme 
\eq{abstract-hyper-disc+}, i.e.
\begin{subequations}\label{C-N-2nd+interp++}
\begin{align}\label{C-N-2nd+a++}
&\DT u_\tau=\underline{\overline v}_\tau,&&u_\tau\big|^{}_{t=0}\!=u_0,&&
\\\label{C-N-2nd+b++}
&\mathscr{T}'\DT v_\tau+\partial\Psi_1(\DT u_\tau)+
\mathfrak{D}_{u}^{}\Phi(\overline u_\tau,\underline u_\tau,\underline z_\tau)
=\overf_\tau,&&v_\tau\big|^{}_{t=0}\!=v_0,
\\&\partial\Psi_2(\DT z_\tau)+
\mathfrak{D}_{z}^{}\Phi(\overline u_\tau,\overline z_\tau,\underline z_\tau)
=\overline g_\tau,
&&z_\tau\big|^{}_{t=0}\!=z_0,
\end{align}\end{subequations}
is to exploit again \eqref{discrete-VI-t+}--\eqref{weak-form-disc+} but with 
$\mathfrak{D}_{u}^{}\Phi(\overline u_\tau,\underline u_\tau,\underline z_\tau)$ and 
$\mathfrak{D}_{z}^{}\Phi(\overline u_\tau,\overline z_\tau,\underline z_\tau)$
in place of $\partial_u\Phi(\overline{\underline u}_\tau,\underline z_\tau)$
and $\partial_z\Phi(\overline u_\tau,\underline{\overline z}_\tau)$, respectively.
It needs still the continuity assumption \eqref{ass-for-Phi-cont}. The simplest 
option is to require \eqref{ass-for-Phi-cont} in the (weak,weak)-topology.
The weak continuity of \eqref{ass-for-Phi-cont} may sometimes be inadequately 
strong, cf.\ e.g.\ the example \eq{quotient} where the weakening of 
\eqref{ass-for-Phi-cont} by using strong topology as far as $(u,\pi)$-variables 
in \eq{quotient} concerns (or here 
as far as $u$-variable in $\mathfrak{D}_{z}^{}\Phi$ concerns). 
This needs to prove the strong convergence of $\overline u_\tau$ still
before executing possibly \eq{trick}. Standardly, one can make it 
by assuming $\Phi(\cdot,z)$ quadratic so that 
$\mathfrak{D}_{u}^{}\Phi(\overline u_\tau,\underline u_\tau,\underline z_\tau)=
\partial_u\Phi(\underline{\overline u}_\tau,\underline{z}_\tau)$, cf.\ 
the first equality in \eq{DPhi-if-quadratic}, and then by
testing \eq{C-N-2nd+b++} by $\underline{\overline u}_\tau$ while assuming a 
uniform monotonicity of $\pl_u\Phi$, cf.\ e.g.\ 
\cite[Step 2 in the proof of Theorem 5.1.2]{MieRou15RIST}, which 
here does not seem to work due to the energy-conserving discretisation
of the inertial term. The same difficulties apply to weakening of
\eq{qualif-b}. Fortunately, sometimes particular techniques work,
cf.\ \cite{Roub??ECTD}. 
%
\end{remark}
}

\begin{remark}[{\it Numerical implementation.}]
\upshape
Often, the potentials $\Psi$'s are sum of quadratic functions
with degree-1 homogeneous functions, and are then nonsmooth at 0.  
The qualification \eqref{struct-ass} or \eqref{structural-ass} then ensures 
that the minimization problems in \eq{C-R-formula-2nd} or in 
\eq{C-R-system-potential} have, after possibly a Mosco-type transformation, a 
structure of {\it Quadratic-Programming} problems (QP)
if $\Psi$'s have polyhedral epigraph (as e.g.\ in damage after space 
discretisation)
or Second-Order Cone Programming (SOCP) if $\Psi$'s have  epigraphs of 
a so-called ice-cream-cone type (as e.g.\ in plasticity),
cf.\ \cite[Sect.\,3.6.3]{MieRou15RIST} or \cite[Sect.\,5]{RouVal??SDSR}.
Both for QP and for SOCP, efficient algorithms and even 
prefabricated software packages
do exist, cf.\ e.g.\ \cite{Dost09OQPA} and \cite{AliGol03SOCP,Stur02IIPM},
respectively.
\end{remark}


\vspace*{1em}
\section{\sf APPLICATION IN CONTINUUM MECHANICS OF SOLIDS}\label{sec-cont-mech}

\def\dev{\mathrm{dev}}
\def\DEV{\R_\mathrm{dev}^{d\times d}}
\def\SYM{\R_\mathrm{sym}^{d\times d}}
\def\wD{w_{\mbox{\tiny\rm D}}}

The approaches from Sections~\ref{sec-quadratic} and \ref{sec-nonconvex}
can be combined. We will illustrate it by considering 
$z=(\pi,\zeta)\in Z:=Z_1\times Z_2$ so that one can thus 
consider, in view of Remarks~\ref{rem-Psi} and \ref{rem-more-dissip},  
the dissipation potentials $\Psi_1:Z_2\times(U\times Z_1)\to\R\cup\{\infty\}$
and $\Psi_2:Z_2\to\R\cup\{\infty\}$ and the system\hspace*{-.1em}
\begin{subequations}\label{plast-dam}\begin{align}\label{plast-dam-a}
\mathscr{T}'\,\DT u+&\pl_{\DT u}\Psi_1(\zeta;\DT u,\DT\pi)
+\pl^{}_u\Phi(u,\pi,\zeta)\ni f,
\ \ \ \ \ \ u|_{t=0}^{}=u_0,\ \ \ \DT u|_{t=0}^{}=v_0,
\\\label{plast-dam-b}
&\pl_{\DT\pi}\Psi_1(\zeta;\DT u,\DT\pi)+\pl^{}_\pi\Phi(u,\pi,\zeta)\ni g,\ \ \ \ \ \ \pi|_{t=0}^{}=\pi_0,
\\\label{plast-dam-c}
&\ \ \ \ \ \ \ \,\pl\Psi_2(\DT\zeta)+\pl^{}_\zeta\Phi(u,\pi,\zeta)\ni h,\ \ \ \ \ \ \zeta|_{t=0}^{}=\pi_0,
\end{align}\end{subequations}
If $\Phi
$ is quadratic in terms of $(u,\pi)$ and also in $\zeta$
separately,  one can devise two fractional steps first for $(u,\pi)$ as in
Sections~\ref{sec-quadratic} and second for $\zeta$ as in this
Section \ref{sec-nonconvex}. In a more general case when $\Phi$ is not 
quadratic in terms of $\zeta$, the quotient 
$\mathfrak{D}_{\zeta}^{}\Phi(u,\pi,\zeta,\tilde\zeta)$ should
be used in place of $\pl^{}_\zeta\Phi(u,\pi,\zeta)$, cf.\ 
Section~\ref{sec-nonquadratic}. It results to the formula
\begin{subequations}\label{disc-hyper-system}
\begin{align}\label{disc-hyper-system-a}
&\frac{u_\tau^k{-}u_\tau^{k-1}\hspace*{-.5em}}\tau=\vhalf,&&
u_\tau^0=u_0,
\\\nonumber
&\mathscr{T}'\frac{v_\tau^k{-}v_\tau^{k-1}\hspace*{-.5em}}\tau
+\pl_{\DT u}\Psi_1\Big(\zeta_\tau^{k-1};\vhalf,\frac{\pi_\tau^k{-}\pi_\tau^{k-1}\hspace*{-.5em}}\tau\hspace*{.5em}\Big)
\\\label{disc-hyper-system-b}
&\hspace*{13.5em}+\pl^{}_u\Phi\big(\uhalf,\pihalf,\zeta_\tau^{k-1}\big)
\!\ni\! f_\tau^k,\!&&v_\tau^0=v_0,
\\\label{disc-hyper-system-c}
&\pl_{\DT\pi}\Psi_1\Big(\zeta_\tau^{k-1};\vhalf,\frac{\pi_\tau^k{-}\pi_\tau^{k-1}\hspace*{-.5em}}\tau\hspace*{.5em}\Big)
+\pl^{}_\pi\Phi\big(\uhalf,\pihalf,\zeta_\tau^{k-1}\big)\ni g_\tau^k,
&&\pi_\tau^0=\pi_0,
\\\label{disc-hyper-system-d}
&\pl\Psi_2\Big(\frac{\zeta_\tau^k{-}\zeta_\tau^{k-1}\hspace*{-.5em}}\tau\hspace*{.5em}\Big)
+\begin{cases}
\pl^{}_\zeta\Phi\big(u_\tau^k,\pi_\tau^k,\zeta_\tau^{k-1/2}\big)\ni h_\tau^k\!\!
&\text{if }\Phi_0(u,\pi,\cdot)\text{ quadratic},\\
\mathfrak{D}_{\zeta}^{}\Phi\big(u_\tau^k,\pi_\tau^k,\zeta_\tau^k,\zeta_\tau^{k-1}\big)
\!\ni\!h_\tau^k\!&\text{in general cases},\end{cases}&&\zeta_\tau^0=\zeta_0.
\end{align}\end{subequations}

Hereafter, we illustrate it on a model for 
a damageable elasto-plastic
body at small strains
occupying a bounded Lipschitz domain $\Omega\subset\R^d$, $d=2$ or $3$,
possibly (cf.\ Remark~\ref{rem-interface}) also in a surface variant.
We will present a relatively general model of a linearized single-threshold 
plasticity with hardening in visco-elastic solid in Kelvin-Voigt rheology 
accompanied with damage allowed possibly for healing. The plastic threshold 
(so-called yield stress $\sigma_{\rm y}$)
determines $S$ as a ball with the radius $\sigma_{\rm y}$. The healing is 
an important phenomenon in some applications (in particular in geophysics)
and particularly exploits combination with plastic slip so that healing can
be realized in the permanently (plastically) deformed configuration, forgetting
the original configuration, cf.\ also \cite[Remark 5.2.24]{MieRou15RIST}. 
For readers' convenience, let us summarize 
the basic notation used in what follows:


\vspace{.5em}

\begin{center}
\fbox{
\begin{minipage}[t]{0.39\linewidth}\small\smallskip
$u$ displacements \\
$\zeta$ damage scalar variable \\
$\pi$ plastic strain  \\
$e(u)$ small strain tensor \\
$\zeta_\flat$ delamination scalar variable \\
$\pi_\flat$ surface plastic slip  \\
$\bbC$ elastic-moduli tensor \\
$\bbD$ viscous-moduli tensor \\
$\bbH$ kinematic-hardening-moduli tensor
\smallskip \end{minipage}
\begin{minipage}[t]{0.47\linewidth}\small\smallskip
$\varrho$ mass density\\
$\gamma$ damage (or delamination) coefficient\\
$S\subset\DEV$ elasticity domain (containing 0)\\
$\alpha$ pseudopotential of damage dissipation,\\ 
$\kappa_1$ plastic strain (or slip) gradient coefficient\\
$\kappa_2$ damage gradient coefficient\\
$K$ elastic modulus of the adhesive\\
$\kappa_0$ hardening of plastic slip
\smallskip \end{minipage}
}\end{center}

\vspace{-.9em}

\begin{center}
{\small\sl Table\,1.\ }
{\small
Summary of the basic notation used through Sections~\ref{sec-cont-mech}
and \ref{sec-comput}. 
}
\end{center}

\medskip

We denote by $\vec{n}$ the outward unit normal to $\partial \Omega$.
We further suppose that the boundary of $\Omega$ splits as
\[
\partial\Omega :=\Gamma= \GD\cup \GN\,,
\]
with $\GD$ and $\GN$ open subsets in the relative topology of
$\partial\Omega$, disjoint one from each other and, up to $(d{-}1)$-dimensional
zero measure, covering $\partial\Omega$. Later, the 
Dirichlet or the Neumann boundary conditions will be prescribed on 
$\GD$ and $\GN$, respectively.
Considering $T>0$ a fixed time horizon, we set
\begin{displaymath}
I:=[0,T], \qquad
Q:=(0,T){\times}\Omega, \qquad \Sigma:=I{\times}\Gamma,
\qquad \Sdir\!:=I{\times}\GD, \qquad \Snew\!:=I{\times}\GN.
\end{displaymath}
Further, $\SYM$ and $\DEV$ will denote 
the set of symmetric or symmetric trace-free (=\,deviatoric) 
$(d{\times}d)$-matrices, respectively. 

In the bulk model, the {\it state} is formed by the triple 
$q:=(u,\pi,\zeta)$. The governing equation/inclusions read as:
\begin{subequations}\label{plast-dam+}
\begin{align}\label{plast-dam1}
&
\varrho\,\DDT{u}-\mathrm{div}\,\sigma_\mathrm{el}
=f_0\ \ \ \text{ with }\ \sigma_\mathrm{el}=\gamma(\zeta)(\mathbb D\DT e_\mathrm{el}+
\mathbb Ce_\mathrm{el}),
&&\!\!\!\!\!\!\!\!\text{\sf(momentum equilibrium)}
\\[.2em]\label{plast-dam12}
&
\partial\delta_{S}^*(\DT{\pi})
\ni\mathrm{dev}\,\sigma_\mathrm{el}
-\mathbb H\pi
+\kappa_1\Delta\pi\ \ \ \text{ with }\ e_\mathrm{el}=e(u){-}\pi,
&&\text{\sf(plastic flow rule)}
\\[-.2em]\label{plast-dam13}
&\partial
\alpha(\DT\zeta)\ni-\frac12
\gamma'(\zeta)\mathbb C
e_\mathrm{el}{:}e_\mathrm{el}
+\kappa_2
\Delta\zeta,
&&\!\text{\sf(damage flow rule)}
\end{align}\end{subequations}
with $\delta_{S}$ the indicator function to a convex set $S$ and $\delta_{S}^*$
its convex conjugate and with ``$\dev$'' denoting the deviatoric part
of a tensor, i.e.\ $\dev\,\sigma:=\sigma-\mathrm{tr}\,\sigma/d$. 
Here, $[\mathbb Ce]_{ij}$ means $\sum_{k,l=1}^d\mathbb C_{ijkl}e_{kl}$.

Of course, \eq{plast-dam+} is to be completed by appropriate boundary 
conditions, e.g.
\begin{subequations}\label{plast-dam-BC}
\begin{align}\label{plast-dam-BC1}
&&&u=\wD&&\text{on }\GD,&&&&
\\\label{plast-dam-BC+}
&&&
\sigma_\mathrm{el}{\cdot}\vec{n}
=f_1&&\text{on }\GN,
\\\label{plast-dam-BC3}&&&\kappa_1\nabla\pi\vec{n}=0\ \ \text{ and }\ \ 
\kappa_2\nabla\zeta{\cdot}\vec{n}=0&&\text{on }\Gamma
\end{align}\end{subequations}
with $\vec{n}$ denoting the unit outward normal to $\Omega$. 
We will consider an initial-value problem 
for \eqref{plast-dam+}--\eqref{plast-dam-BC} by asking for 
\begin{align}\label{plast-dam-IC}
u(0)=u_0,\ \ \ \ \ \DT u(0)=v_0,\ \ \ \ \ \pi(0)=\pi_0,\ \ \text{ and }\ \ \zeta(0)=\zeta_0.
\end{align}
The abstract spaces and the energy functionals used in \eqref{plast-dam+} are 
now:
\begin{subequations}\label{elast-visco-plast}\begin{align}
&U=H^1(\Omega;\R^d),\ \ \ \ \ Z_1
=\begin{cases}L^2(\Omega;\R^{d\times d}_{\rm dev})&\text{if }\kappa_1=0,\\
H^1(\Omega;\R^{d\times d}_{\rm dev})&\text{if }\kappa_1>0,\end{cases}
\ \ \ \ \ Z_2=H^1(\Omega)\cap L^\infty(\Omega),
\\&
\Phi(u,\pi,\zeta)=\int_\Omega\frac12\gamma(\zeta)\bbC
(e(u){-}\pi){:}(e(u){-}\pi)+\frac12\bbH\pi{:}\pi
+\frac{\kappa_1}2|\nabla\pi|^2+\frac{\kappa_2}2|\nabla\zeta|^2\,\d x
,
\\&\Psi_1(\zeta;\DT u,\DT\pi)=
\int_\Omega\frac12\gamma(\zeta)\bbD(e(\DT u){-}\DT\pi){:}(e(\DT u){-}\DT\pi)
+
\delta_S^*(\DT\pi)\,\d x,\qquad
\Psi_2(\DT\zeta)=\int_\Omega
\alpha(\DT\zeta)\,\d x,
\label{elast-visco-plast-Psi3}
\\&\mathscr{T}(\DT u)=\int_\Omega\frac\varrho2|\DT u|^2\,\d x,\qquad\quad
\big\langle f(t),u\big\rangle=\int_\Omega\!f_0(t,\cdot){\cdot}u\,\d x+
\int_{\GN}\!\!
f_1(t,\cdot){\cdot}u\,\d S,
\end{align}\end{subequations}
%
with $\bbC$ and $\bbD$ positive-definite 4th-order tensors, $\bbH$
a positive-definite 4th-order tensor, $\alpha:\R\to[0,\infty]$ convex with 
$\alpha(0)=0$, and $\kappa_1\ge0$ and $\kappa_2>0$ given coefficients.
Furthermore, $\gamma:\R\to\R^+$ is positive (i.e.\ allowing only for an 
incomplete damage) continuously differentiable nondecreasing with $\gamma'=0$
on $(-\infty,0]\cup[1,\infty)$, which ensures that the values of $\zeta$ ranges
the interval $[0,1]$ if $\zeta_0$ do so.  

Written in the classical formulation, the differential quotient used in 
\eqref{disc-hyper-system-d} can be now 
taken as
\begin{align}\label{quotient}
\mathfrak{D}_{\zeta}^{}\Phi\big(u,\pi,\zeta,\tilde\zeta\big):=
\begin{cases}
\displaystyle{
\frac12\frac{\gamma(\zeta){-}\gamma(\tilde\zeta)}{\zeta-\tilde\zeta}
\bbC(e(u){-}\pi){:}(e(u){-}\pi)}&
\\[-.4em]\hspace{7.2em}
\displaystyle{-\frac12\kappa_2\Delta\zeta-\frac12\kappa_2\Delta\tilde\zeta}
&\text{on }\{x\!\in\!\Omega;\ 
\zeta(x)\ne\tilde\zeta(x)\},
\\\displaystyle{\frac12\gamma'(\zeta)\bbC
(e(u){-}\pi){:}(e(u){-}\pi)-\kappa_2\Delta\zeta}&\text{on }\{x\!\in\!\Omega;\ 
\zeta(x)=\tilde\zeta(x)\}.
\end{cases}
\end{align}
Note that it obviously satisfies (\ref{ass-for-Phi}a,b). Moreover, the
assumption \eqref{ass-for-Phi-cont} holds in the variant of the 
(strong$\times$weak$\times$weak$\times$weak,weak*)-continuity
since $\gamma(\cdot)$ is assumed continuously
differentiable with $\gamma'$ bounded and $\kappa_2>0$ so that we can use 
Rellich compact-embedding theorem for $\zeta$ and $\tilde\zeta$.
This needs to prove strong convergence of $\overline u_\tau$ mentioned in 
Remark~\ref{rem-conv-modified}. This may be quite technical. Here,
it holds for a special case that $\bbD=\chi\bbC$ for some relaxation time 
$\chi>0$, see \cite[Step 3 in the proof of Prop.\,4]{Roub??ECTD}. Also 
\eqref{ass-for-Phi-pot-z} holds with the potential 
$\mathfrak{G}_{u,\pi,\tilde\zeta}$ given by
\begin{align}&\nonumber
\mathfrak{G}_{u,\pi,\tilde\zeta}(\zeta)=
\int_\Omega\!\bigg(\frac{\kappa_2}4|\nabla\zeta|^2
+\frac12\bbC(e(u(x)){-}\pi(x)){:}(e(u(x)){-}\pi(x))\int_0^{\zeta(x)}
\!\!\!\!\!\!\varGamma_{\tilde\zeta(x)}(z)\,\d z\!\bigg)\d x
\\&\hspace*{14em}\text{with }\ \ \ \varGamma_{\!\tilde z}^{_{}}(z)=
\begin{cases}
(\gamma(z){-}\gamma(\tilde z))/(z-\tilde z)
&\text{if }z\ne\tilde z,\\\qquad\quad\gamma'(z)&\text{if }z=\tilde z.\ \ \
\end{cases}
\end{align}

For the time discretisation \eqref{disc-hyper-system},
one considers the structure \eq{plast-dam} and
use the Crank-Nicolson scheme for displacement with plasticity together
and then backward Euler discretization
for damage. It results to the formula \eq{disc-hyper-system}.
Note that coupling of $\DT u$ and $\DT\pi$ in the dissipation
potential $\Psi_1$ in \eqref{elast-visco-plast}
allows for considering the viscous dissipation acting more physically 
on the elastic strain $e(u){-}\pi$ instead on the total strain $e(u)$ 
and the plastic strain $\pi$ separately.


\TINY{The sophisticated construction of a 
mutual recovery sequence by M.\,Thomas \cite{Thom10PhD} devised 
for the backward Euler formula 
unfortunately does not work if 
because it does not ensure 
$\langle\overline \zeta_\tau(t){-}\tilde \zeta_\tau,\underline \zeta_\tau(t)\rangle_{\Phi,\overline u_\tau(t),\overline\pi_\tau(t)}
\to\langle \zeta(t){-}\tilde \zeta,\zeta(t)\rangle_{\Phi,u(t),\pi(t)}$
needed in \eqref{binomial-MRS}. 
Then the arguments of the proof of
\cite[Proposition 4.2.30]{MieRou15RIST} used for 
$\tilde \zeta_\tau=\max(0,\tilde \zeta-\|\overline \zeta_\tau(t){-}\zeta(t)\|_{C(\bar\Omega)}^{})$
applies ....................
to perform the limit passage from \eqref{disc-semi-stab+}
towards \eqref{ES-semistability} with $\Phi_3$ in place of $\Phi_2$.
Here, as $\Phi_2$ in \eq{elast-visco-plast-Psi3} is degree-1 homogeneous,
we need another semi-stability also for $\pi$-variable; for
the fully quasistatic case with $\mathscr{T}=0$ and $\Psi_1=0$ 
cf.\ \cite[Def.\,3.1]{RouVal??SDSR}. The modification for
the Crank-Nicolson formula is here as in \eq{binomial-MRS}
because $\Phi$ is quadratic in terms of $\pi$. 
}


\begin{remark}[{\it Ambrosio-Tortorelli's approximation of cracks.}]
\upshape
The standard case of rate-independent first-order damage gradient can 
be applied in the energy-conserving Yanenko scheme when a cohesive damage which
does not need the constraint $\zeta\ge0$ is considered. This situation occurs 
if $\gamma'(0)=0$. In particularly, it holds for 
$\gamma(\zeta)=\epsilon+\zeta^2$ with $\epsilon>0$ which
occurs in particular in the so-called Ambrosio-Tortorelli 
approximation \cite{AmbTor90AFDJ} of the fracture, rigorously devised in the 
static scalar case 
and then suggested also for the dynamical situation in 
\cite[Sect.\,5.2.5]{MieRou15RIST} without guaranteeing any convergence 
towards the brittle fracture in the visco-elastic bulk, however.
Keeping $\zeta$ valued in $[0,1]$ can then be made by 
considering $\alpha([0,\infty))=\infty$ so that no healing is allowed.
The analysis of the energy-conserving scheme from Sect.\,\ref{sec-nonconvex}
even coupled with diffusion is in \cite{Roub??ECTD}.
\end{remark}

\begin{remark}[{\it Polycrystalic shape-memory alloys.}]
\upshape
Another noteworthy component-wise quadratic model for phase-transformations 
in shape-memory materials which considers 
$e_{\rm el}=e(u)-\lambda\pi_{\rm tr}$ in \eqref{elast-visco-plast} with
$0\le\lambda\le1$ a volume fraction between so-called austenite and martensite,
and $\pi_{\rm tr}$ a transformation strain that is subjected to a constraint on 
its magnitude.  Cf.\ the polycrystalic models in 
\cite{AuReSt07TDMD,SadBha07MICM} possibly also in combination with 
plasticity like already used in \eqref{elast-visco-plast} with 
a decoupled dissipation similarly like in \eqref{plast-dam}.
A suitable transformation allows more coupled dissipation, 
cf.\ \cite{FrBeSe14MMCM,SFBBS12TMNT},
and then, after penalization of the constraints, a one-step formula 
as used in Sect.~\ref{sec-quadratic} but modified as in Sect.~\ref{sec-nonquadratic}.

\end{remark}

\TINY{
\begin{remark}[{\it Special case: elasto-plasto-dynamics.}]\label{rem-elasto-plasto}
\upshape
Considering $\Psi_1\equiv0$ and $\Psi_2$ 1-homogeneous coercive 
on a Banach space $Z_1\subset Z$ while $\Phi$ 
quadratic would cover the elasto-plasto-dynamics at small strains 
which, without involving
any viscosity, is a hyperbolic problem in a medium involving nonlinear 
processes. It is generally very difficult class of problems but here
complies with the structure \eqref{struct-ass}. Typically, if a 
(e.g.\ kinematic) hardening or a gradient plasticity is considered, then
$u:\Omega\to\R^d$ is a displacement and $z:\Omega\to\R^{d\times d}$ is
a plastic tensor, $U=H^1(\Omega;\R^d)$, $H=L^2(\Omega;\R^d)$, and 
$Z=L^2(\Omega;\R^{d\times d})$ or $Z=H^1(\Omega;\R^{d\times d})$ respectively, 
while 
$Z_1=L^1(\Omega;\R^{d\times d})$ with $\Omega\subset\R^d$ a domain occupied
by the elasto-plastic body; note that $Z_1$ is not reflexive.
However, the loading qualification must be strengthened to $f\in L^1(I;H)$ 
and the last two estimates in \eqref{C-N-2nd-est1}
and the estimates in \eqref{C-N-2nd-est2}
must be weakened to
\begin{align}\label{C-N-2nd-est+}
\|v_\tau\|^{}_{L^\infty(I;H)
}\le C,\ \ \ \ 
\|z_\tau\|^{}_{L^\infty(I;Z)\,\cap\,W^{1,1}(I;Z_1)}\le C,\ \ \ \ 
\|\DT v_\tau\|^{}_{
L^1(I;H)}\le C,\ \ \ \ 
\|\DT u_\tau\|^{}_{L^\infty(I;H)
}\le C.
\end{align}
The convergence of \eqref{weak-form-disc} towards \eqref{weak-form}
with $\Psi_1$ omitted still holds and with $\int_0^T\Psi_2(\DT z)\,\d t$
understood in the sense of measures, cf.\ \eq{variation}. Here it is 
important that the weak convergences $u_\tau(T)\to u(T)$ in $U$, 
$v_\tau(T)\to v(T)$ in $H$, and $z_\tau(T)\to z(T)$ in $Z$ still hold due to 
$L^\infty(I)$-estimates valued in these spaces and due to estimates on the time 
derivatives; in the later case, the Helly selection principle is to be used. 
This allows still for a limit passage in the right-hand side of 
\eqref{weak-form-disc} by weak* lower semicontinuity.
The elasto-plasto-dynamics at small strains has been first rigorously 
analyzed by Gr\"oger \cite{Grog78TDVE}.
If neither hardening nor gradient plasticity is considered, then 
additional analytical difficulties related with the so-called perfect
plasticity occurs. 
............... \COMMENT{seems OK - $\partial_u\Phi(...)\to$
and $\partial_z\Phi(...)\to$ weakly* TO CHECK}
Such problems in one-dimension can also be handled by the theory
of (scalar) hysteresis operators \cite{Krej97HCDH}.   
\end{remark}
}


\TINY{
\begin{remark}[{\it More general dissipation.}]\label{rem-more-dissip+}
\upshape
Making $\Psi_2$ dependent on $u$ as in Remark~\ref{rem-more-dissip} would 
allow for a friction dependent on the normal traction stress
in a normal-compliance regularization, assuming that the specimen is
surely pressed towards the contact boundary $\GC$ so that the sign
of the traction stress does not change and one can consider a bilateral
contact without violating the qualification \eqref{structural-ass-a}.
Moreover, 
as in \eqref{plast-dam},??????????????????????????? the dissipation 
functionals $\Psi_1$ and $\Psi_2$ in \eqref{delamination} can be coupled, 
which here would allow for considering also a term 
$\int_{\GC}\frac12D|u{-}\pi|^2\,\d S$\COMMENT{???} to model a viscosity in the adhesive.
\end{remark}
}

\begin{remark}[{\it Surface plasticity and damage.}]\label{rem-interface}
\upshape
\COLOR{Considering a contact interface $\GC$ as
a part of the boundary of $\Omega$, a useful surface
 analog of the bulk model exploits as} the internal variables 
a surface {\it plastic slip} $\pi_\flat:\GC\to\R^{d-1}$ as in 
\cite{RoKrZe13DACM,RoMaPa13QMMD}
and a surface damage (also called a {\it delamination} parameter)
$\zeta_\flat:\GC\to[0,1]\subset\R$ as invented by Fr\'emond \cite{Frem85DLS}.
A combination with the plastic slip was devised in 
\cite{RoKrZe13DACM,RoMaPa13QMMD} to model
a mode-mixity dependent delamination, reflecting the phenomenon
that the Mode II (shear) needs/dissipates usually considerably more
energy to delaminate than Mode I (opening).
Advantageously, compactness of the trace operator $u\mapsto u|_{\GC}^{}$ 
simplifies some analytical aspects. 
%
%
%
%
%
%
%
\COLOR{In the case of a so-called normal-compliance adhesive contact 
(with the compliance described by a function $p(\cdot)$),
}
the classical formulation of such problem consists in the equilibrium of forces 
on the domain $\Omega$ with the boundary condition on $\GN=\Gamma{\setminus}\GC$
and several complementarity problems on the contact boundary $\GC$.
In the Kelvin-Voigt rheology, the model looks as:
\begin{subequations}\label{eq5:delam-class-form}\begin{align}
\label{eq5:delam-class-form1}
&
\DT u=v,\qquad \varrho\DT v-\mathrm{div}\,\sigma=0,
\qquad \sigma=\bbC e(u)+\bbD e(v),
&&\text{in }
\Omega
,
\\\label{eq5:delam-class-form3}  
&\sigma\nu=
f_1(t)
&&\text{on }
\Gamma,
\\\label{e5:class-form-d-small}
   &\left.\begin{array}{ll}
&\hspace{-1.7em}
\sigma_{\rm n}+\gamma(\zeta_\flat)K u_{\rm n}+p(u_{\rm n})=0,
\qquad\sigma_{\rm t}=\gamma(\zeta_\flat)
K(u_{\rm t}{-}\pi_\flat),
\\[.3em]
&\hspace{-1.7em}
\pl\delta_S^*(\DT\pi_\flat)\ni\sigma_{\rm t}-\kappa_0\pi_\flat+\divS(\kappa_1\nablaS\pi_\flat),
\\[.3em]
&\hspace{-1.7em}
\pl\alpha(\DT{\zeta}_\flat)\ni-\frac12\gamma'(\zeta_\flat)K|u_{\rm t}{-}\pi_\flat|^2
+\divS(\kappa_2\nablaS\zeta_\flat),
   \end{array}\right\}\!\!\!
   &&\text{on }\GC,
\end{align}\end{subequations}
where 
$\eta$ is the driving ``force'' for the plastic-slip evolution,
$\sigma\nu:=(\bbC e(u)\COLOR{{+}\bbD e(v)})\big|_{\Gamma}\nu$ 
is the traction stress on $\Gamma=\GC$ or 
$\GN$.
Moreover, its normal and tangential components on $\GC$ are denoted with 
$\sigma_{\rm n}(u)=(\sigma\nu){\cdot}\nu$ and $\sigma_{\rm t}(u)=
\sigma\nu-((\sigma\nu){\cdot}\nu)\nu$, respectively, so that we
have the decomposition $\sigma\nu=\sigma_{\rm n}\nu+\sigma_{\rm t}$.
%
In \eqref{e5:class-form-d-small}, $\divS:={\rm trace}(\nablaS)$ denotes the 
$(d{-}1)$-dimensional ``surface divergence'' and $\nablaS$ a ``surface
gradient'', i.e.\ the tangential derivative defined as
$\nablaS v=\nabla v-(\nabla v{\cdot}\nu)\nu$ for
$v$ defined on $\GC$. Actually, we assume here that $\GC$ is flat,
otherwise also a curvature term like $(\divS\nu)(\kappa\nablaS\pi_\flat\nu)$
should contribute to the corresponding driving force $\eta$. 
This problem has a structure \eqref{plast-dam} if 
the stored-energy, the kinetic-energy, and the dissipation functionals and
the loading are set now as:
\begin{subequations}\label{delamination}
\begin{align}\nonumber
&\Phi(u,\pi_\flat,\zeta_\flat)=
\int_\Omega\frac12\bbC e(u){:}e(u)\,\d x
+\!\int_{\GC}
\COLOR{\frac12\gamma(\zeta_\flat)K}\big(u_{\rm n}^2+
|u_{\rm t}{-}\pi_\flat|^2\big)+\widehat p(u_{\rm n})
\\[-.4em]&\label{delamination-Phi}\hspace{17em}
+\frac{\kappa_0}2|\pi_\flat|^2
+\frac{\kappa_1}2|\nablaS\pi_\flat|^2+\frac{\kappa_2}2|\nablaS\zeta_\flat|^2
\,\d S,
\\\label{delamination-Psi}
&
\Psi_1(\DT u,\DT\pi_\flat)=\int_\Omega\frac12\bbD e(\DT u){:}e(\DT u)\,\d x
+\int_{\GC}\!\!
\delta_S^*(\DT\pi_\flat)\,\d S,\qquad\Psi_2(\DT\zeta_\flat)
=\int_{\GC}\!\!\alpha(\DT\zeta_\flat)\,\d S,
\\\label{delamination-T-f}
&\mathscr{T}(\DT u)=\int_\Omega\frac\varrho2|\DT u|^2\,\d x,
\qquad
\qquad\big\langle f(t),u\big\rangle=
\int_\Gamma\!\!
f_1(t,\cdot){\cdot}u\,\d S,
\intertext{with $\widehat p$ a primitive function to $p$. The underlying 
Banach spaces are taken as}
\label{delamination-U-Z}
&U=H^1(\Omega;\R^d)
,\qquad Z_1=
\begin{cases}L^2(\GC;\R^{d-1})&\text{if }\kappa_1=0,\\
H^1(\GC;\R^{d-1})&\text{if }\kappa_1>0,\end{cases}
\qquad Z_2=L^\infty(\GC).
\end{align}\end{subequations}
If $[\zeta_\flat]_0\in[0,1]$ on $\GC$ (as usually assumed), 
then $\zeta_\flat$ remains valued 
in $[0,1]$ during the whole evolution \COLOR{provided $\gamma(\cdot)$ is 
qualified as before. Now we can easily afford the complete surface
damage (delamination), i.e.\ $\gamma(\zeta)=0$ if $\zeta\le0$.}
%
\end{remark}

\begin{remark}[{\it Stress formulation, Maxwell/Jeffreys rheology.}]\label{rem-Max}
\upshape
A modification of the original 1st-order system \eqref{eq5:delam-class-form1} 
by elimination of $u$-variable leads to a stress/velocity formulation
\begin{align}\label{KV}
\bbC^{-1}\DT\sigma=e(v)+\bbC^{-1}\bbD e(\DT v)\ \ \ \ \ \ \&\ \ \ \ \ \ 
\varrho\DT v-{\rm div}\,\sigma=0.
\end{align}
This allows for a combination with the so-called mixed FEM,
cf.\ \cite{ByJoTs02NFMF} for the case $\bbD=0$.
It further allows for a straightforward modification 
by adding the terms $\bbD_{\rm Max}^{-1}\sigma$ and $\bbD_{\rm Max}^{-1}\bbD e(v)$:
\begin{align*}
\bbC^{-1}\DT\sigma
+\bbD_{\rm Max}^{-1}\sigma=(\mathbb I{+}\bbD_{\rm Max}^{-1}\bbD)e(v)+\bbC^{-1}\bbD e(\DT v)
\ \ \ \ \ \ \&\ \ \ \ \ \ 
\varrho\DT v-{\rm div}\,\sigma=0
\end{align*}
with $\bbD_{\rm Max}$ another viscous-moduli tensor. For $\bbD=0$ and 
$\bbD_{\rm Max}>0$ we thus obtain the Maxwell rheology while for 
$\bbD>0$ and $\bbD_{\rm Max}>0$ such rheological model is called the 
Jeffreys material as a Norton-Hoff (Stokes) and a Maxwell rheologies in parallel.
For $\bbD_{\rm Max}\to\infty$ we obtain the Kelvin-Voigt model \eqref{KV}.
Both are fluids in the sense that they cannot 
permanently withstand stress under a bounded strain response.
Creep effects can thus be modeled. After being discretised by Crank-Nicolson 
scheme, it looks as:
\begin{align*}
&
\varrho\frac{v_\tau^k{-}v_\tau^{k-1}\hspace*{-.5em}}\tau
-{\rm div}
\,\sigma_\tau^{k-1/2}=
0,
\\&
\bbC^{-1}\frac{\sigma_\tau^k{-}\sigma_\tau^{k-1}\hspace*{-.5em}}\tau
+\bbD_{\rm Max}^{-1}\sigma_\tau^{k-1/2}
=(\mathbb I{+}\bbD_{\rm Max}^{-1}\bbD)e
\big(v_\tau^{k-1/2}\big)+\bbC^{-1}\bbD e\Big(\frac{v_\tau^k{-}v_\tau^{k-1}\hspace*{-.5em}}\tau\hspace*{.5em}\Big).
\end{align*}
To reveal the energetics, we express $e(v_\tau^{k-1/2})=\bbC^{-1}(\frac{\sigma_\tau^k{-}\sigma_\tau^{k-1}\!\!}\tau
-\bbD e(\frac{v_\tau^k{-}v_\tau^{k-1}}\tau))
+\bbD_{\rm Max}^{-1}(\sigma_\tau^{k-1/2}\!-\bbD e(v_\tau^{k-1/2})$ and 
test it by $
\sigma_\tau^{k-1/2}\!-\bbD e(v_\tau^{k-1/2})$. 
After using several binomial formulas similarly as before, we obtain the
discrete energy balance (as an equality):
 \begin{align*}\nonumber
\frac1\tau\bigg(\int_\Omega\frac12\bbC^{-1}(\sigma_\tau^k\!-\bbD e(v_\tau^k)){:}(\sigma_\tau^k\!-\bbD e(v_\tau^k))
-
\frac12\bbC^{-1}(\sigma_\tau^{k-1}\!-\bbD e(v_\tau^{k-1})){:}(\sigma_\tau^{k-1}\!-\bbD e(v_\tau^{k-1}))\,\d x\bigg)
\\\nonumber\qquad+\int_\Omega\bbD e(v_\tau^{k-1/2}){:}e(v_\tau^{k-1/2})
+\bbD_{\rm Max}^{-1}(\sigma_\tau^{k-1/2}\!-\bbD e(v_\tau^{k-1/2})){:}(\sigma_\tau^{k-1/2}\!-\bbD e(v_\tau^{k-1/2}))\,\d x
\\\hspace{15em}=\int_\Omega\!f_\tau^k\cdot v_\tau^k\,\d x
+
\int_\Gamma\!\!g_\tau^k\cdot v_\tau^k\,\d S.
\end{align*}
\end{remark}


\vspace*{1em}
\section{\sf ILLUSTRATIVE COMPUTATIONAL SIMULATIONS}\label{sec-comput}

We illustrate the above sections on 
a Mode-II frictional-type or adhesive 
contact problem 
\COLOR{from Remark~\ref{rem-interface}.}
All these inelastic process occurs only on (a part)
of the boundary, let us denote it by $\GC$. The resting part of the
boundary will be denoted by $\GN$ and part of it, denoted by $\GNN$, will 
be loaded by time-dependent force.

The time-discretisation scheme \eqref{disc-hyper-system} has been used. 
In particular, it complies with the energy conservation at all time intervals 
when the delamination $\zeta_\flat$ does not evolve. The space-discretisation was performed by FEM, namely by using 
the Q1-elements for $u$ and P0-elements for both $\pi_\flat$ and $\zeta_\flat$ on 
the boundary $\GC$. For the special case where $\bbD{=}\chi\bbC$, with $\chi$ a given relaxation time, damping matrix may be defined as proportional to stiffness, that is a special case for the device known as Rayleigh damping \cite{cook2001concepts}. Here we exploited that we consider only such particular 
cases which do not require gradient of $\pi_\flat$ 
i.e.\ both $\kappa=0$ 
can be used for the two experiments below because they consider only either
the plastic slip or the delamination, but not their combination, 
cf.\ also Remark~\ref{ram-dam+plast}.

The isotropic \COLOR{visco}elastic material of the bulk is aluminum with the Young modulus $E=70$GPa, Poisson ratio $\nu=0.35$, \COLOR{i.e.\ 
$\bbC_{ijkl}
\doteq[60\delta_{ij}\delta_{kl}+26(\delta_{ik}\delta_{jl}{+}\delta_{il}\delta_{jk})]{\rm GPa}$,
and mass density $\varrho=2700\,$kg/m$^3$. Elastic plain strain 
is considered \COLOR{for 2-dimensional computational experiments. Moreover, 
we consider a Mode-II contact in the sense that $u_{\rm n}=0$ on $\GC$, which 
means sending the compliance slope $p(\cdot)$ in \eqref{e5:class-form-d-small} 
to infinity}. Rather formally, we consider a very small viscosity 
$\bbD=\chi\bbC$ with $\chi=2\,$ns implemented simply as in \cite{PaMaRo14BVE}
to comply with our convergence-analysis arguments, although it does not have 
visible effects in our simulations.}
\COLOR{We also do not consider hardening, i.e.\ $\kappa_0=0$ in 
\eq{delamination-Phi}; in fact, the coercivity of 
$\Phi$ in terms of $\pi_\flat$ is then ensured through the coercivity 
of $u_{\rm t}$ by cooperation of Korn's inequality and inertia
if $\inf\gamma(\cdot)>0$ or through the coercivity
of $\alpha$ in  \eq{delamination-Psi} with $\kappa_1>0$, which is
applicable in Sections~\ref{sec-friction} and \ref{sec-delam}.}

\COLOR{We perform our illustrative computational experiments for $d=2$, 
using} a rectangular 2-dimensional domain $\Omega$ in Fig\COLOR{s}.\,\ref{fig4-dam-plast-geom} and 
\ref{fig4-dam-plast-geom+}. For the spatial discretisation, $\Omega$ was 
divided into 2$\times$40=80 squares 
for the coarsest spacial discretisation, and then also refined 4 and 9 times
for 320 and 720 square elements.

\COMMENT{THE SECOND EXAMPLE NOW NEEDS $\nablaS\zeta_\flat$ IF 
IT WOULD BE APPROXIMATED BY NONLINEAR $\gamma(.)$ TO ALLOW FOR 
ENERGY-CONSERVING DISCRETISATION :-(( SO I ERRASED DETAILS ABOUT 
DISCRETISATION}
\TINY{
Let us 
remark that the following
experiments allow for ignoring gradients by setting 
$\kappa$'s 
zero, which is what we used when implementing P0-FEM for 
$\pi_\flat$ and $\zeta_\flat$ on $\GC$. 
As to the former case, the mutual recovery sequence can be
constructed as in the proof of Proposition~\ref{prop-RI-quadratic}
because $\pi_\flat\mapsto\Phi(u,\pi_\flat,1)$ from \eq{delamination-Phi} is quadratic 
without no constraints and $\Phi_2$  from \eq{delamination-Psi} is degree-1 
homogeneous. As to the later case, we can exploit that 
$\zeta_\flat\mapsto\Phi(u,0,\zeta_\flat)$ from \eq{delamination-Phi} is affine and
use the standard construction of the mutual recovery sequence 
$[\tilde\zeta_\flat]_\tau=[\overline\zeta_\flat]_\tau(t)\tilde\zeta_\flat/\zeta_\flat(t)$ which 
converges to $\tilde\zeta_\flat$ weakly* in $L^\infty(\GC)$ while complying with the 
constraints $0\le[\tilde\zeta_\flat]_\tau\le[\overline\zeta_\flat]_\tau(t)$ coming 
from $\Phi$ 
and $\Psi_3$\COMMENT{NOTATION TO CHECK} provided $0\le\tilde\zeta_\flat\le\zeta_\flat(t)$, and this attribute is 
here inherited for the modified sequence needed for \eq{disc-semi-stab+} with 
$\Psi_3$\COMMENT{NOTATION TO CHECK} in place of $\Phi_2$, namely $[\tilde\zeta_\flat]_\tau{+}[\underline\zeta_\flat]_\tau(t)$
again converges weakly* in $L^\infty(\GC)$ and here this mode
of convergence suffices because of the mentioned affinity of $\Phi(u,0,\cdot)$
and because of the compactness of the trace operator $u\mapsto u|_{\GC}$.}



\subsection{Frictional-contact experiment without adhesion}\label{sec-friction}

As a particular situation, we send $\alpha$ to infinity which
practically means that $\DT\zeta_\flat=0$ and $\zeta_\flat$ stays constant
during the whole evolution, say $\zeta_\flat=[\zeta_\flat]_0\equiv1$. 
In other words,
we just omit $\zeta_\flat$ here. Further modeling ansatz is to consider $K$ 
large on $\GC$, namely $K=75\,$GPa/m. 
This causes that the difference $u_{\rm t}{-}\pi_\flat$ is small
during the whole evolution and we can recognize the model as 
a regularized dry-friction model with $\sigma_\mathrm{y}$ playing a role of
the friction coefficient, cf.\ also \cite[Remark 5.2.17]{MieRou15RIST}.
More specifically, it is a so-called given friction (also called 
Tresca's friction) which neglects influence of possible variation of the 
normal force on the tangential friction.
\begin{figure}[th]
\begin{center}
\psfrag{GN}{\footnotesize $\GNNN$}
\psfrag{GD}{\footnotesize $\GNN$}
\psfrag{GC}{\footnotesize $\GC$}
\psfrag{elastic}{\footnotesize elastic body}
\psfrag{obstacle}{\footnotesize rigid obstacle}
\psfrag{adhesive}{\footnotesize adhesive}
\psfrag{LC}{$L_c$}
\psfrag{L}{\footnotesize $L=\ $250\,mm}
\psfrag{H}{\footnotesize $H=$}
\psfrag{12.5}{\scriptsize 12.5\,mm}
\psfrag{loading}{\footnotesize loading}
\hspace*{-.5em}\vspace*{-.1em}{\includegraphics[width=.95\textwidth]{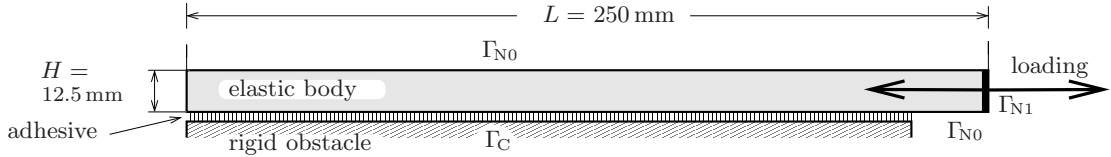}}
\end{center}
\vspace*{-1em}
\caption{\sl Geometry of a 2-dimensional rectangular-shaped specimens 
subjected to a cyclically loading $f_1=f_1(t)$ the 
right-hand side $\GNN$.
}
\label{fig4-dam-plast-geom}
\end{figure}
To show influence of the friction, two values of the yield stress determining
 the sliding resistance (i.e.\ the prescribed Tresca friction) has been chosen,
namely $\sigma_\mathrm{y}=3\,$MPa and 6\,MPa, referred respectively as ``small'' 
and ``large'' friction in Figures~\ref{fig-loading-f+}--\ref{fig-forc-disp}.
We consider a 2-dimensional rectangular specimen, cf.\ Fig.\,\ref{fig4-dam-plast-geom}. The adhesive part \COLOR{$\GC$} is considered to be that of 0.9 
of the total length, \COLOR{cf.\ again Fig.\,\ref{fig4-dam-plast-geom}.} 

\begin{figure}[th]
\begin{center}
\TINY{\hspace*{-1.5em}\vspace*{-.1em}\includegraphics[width=.36\textwidth,height=.3\textwidth]{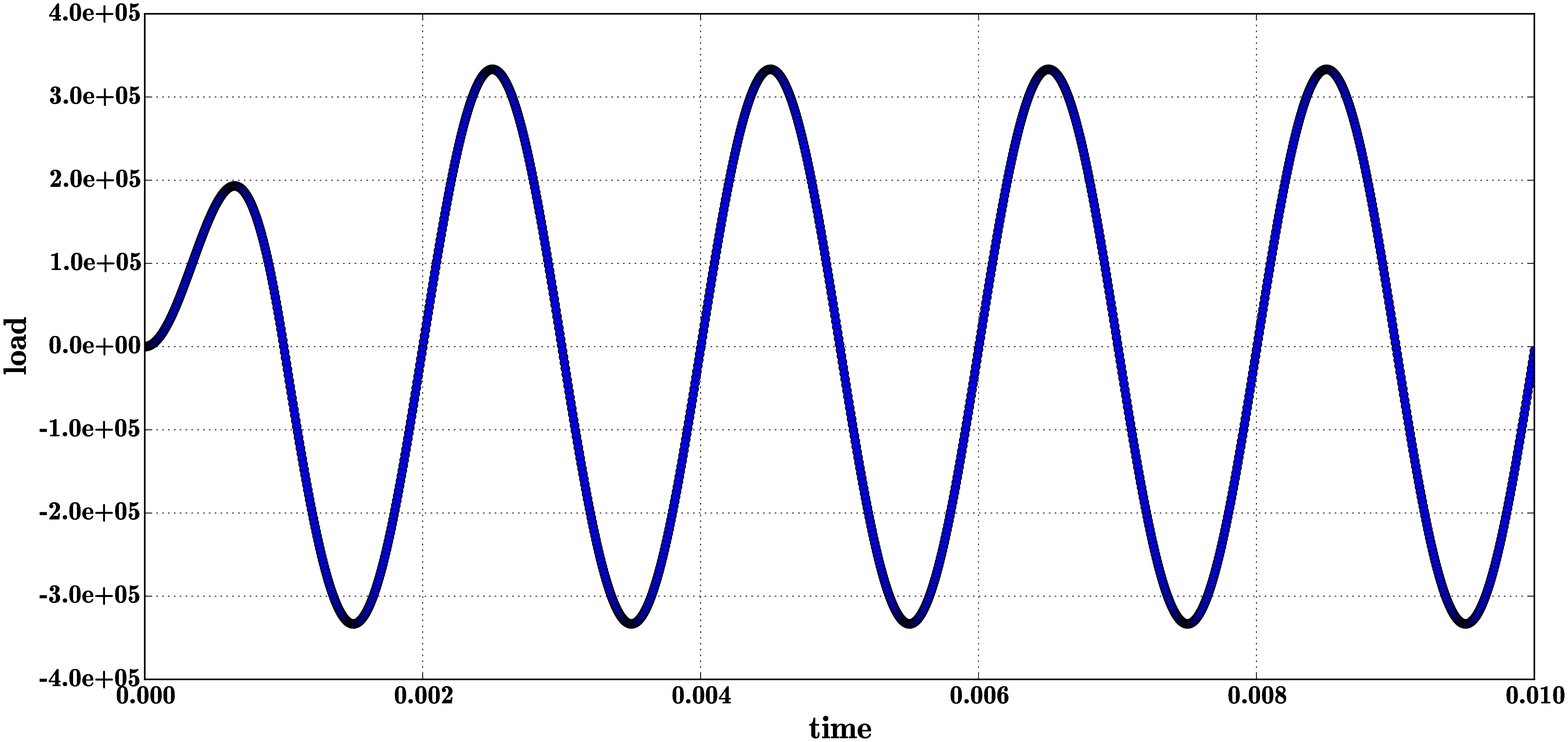}
\hspace*{-1.em}\includegraphics[width=.36\textwidth,height=.3\textwidth]{elastoplastNRG.eps}
\hspace*{-1.em}\includegraphics[width=.36\textwidth,height=.3\textwidth]{kineticNRG.eps}\hspace*{-2.5em}
\\}
\hspace*{-.5em}\vspace*{-.1em}\includegraphics[width=.36\textwidth,height=.25\textwidth]{cyclic_loading.eps}
\hspace*{-1.em}\includegraphics[width=.36\textwidth,height=.25\textwidth]{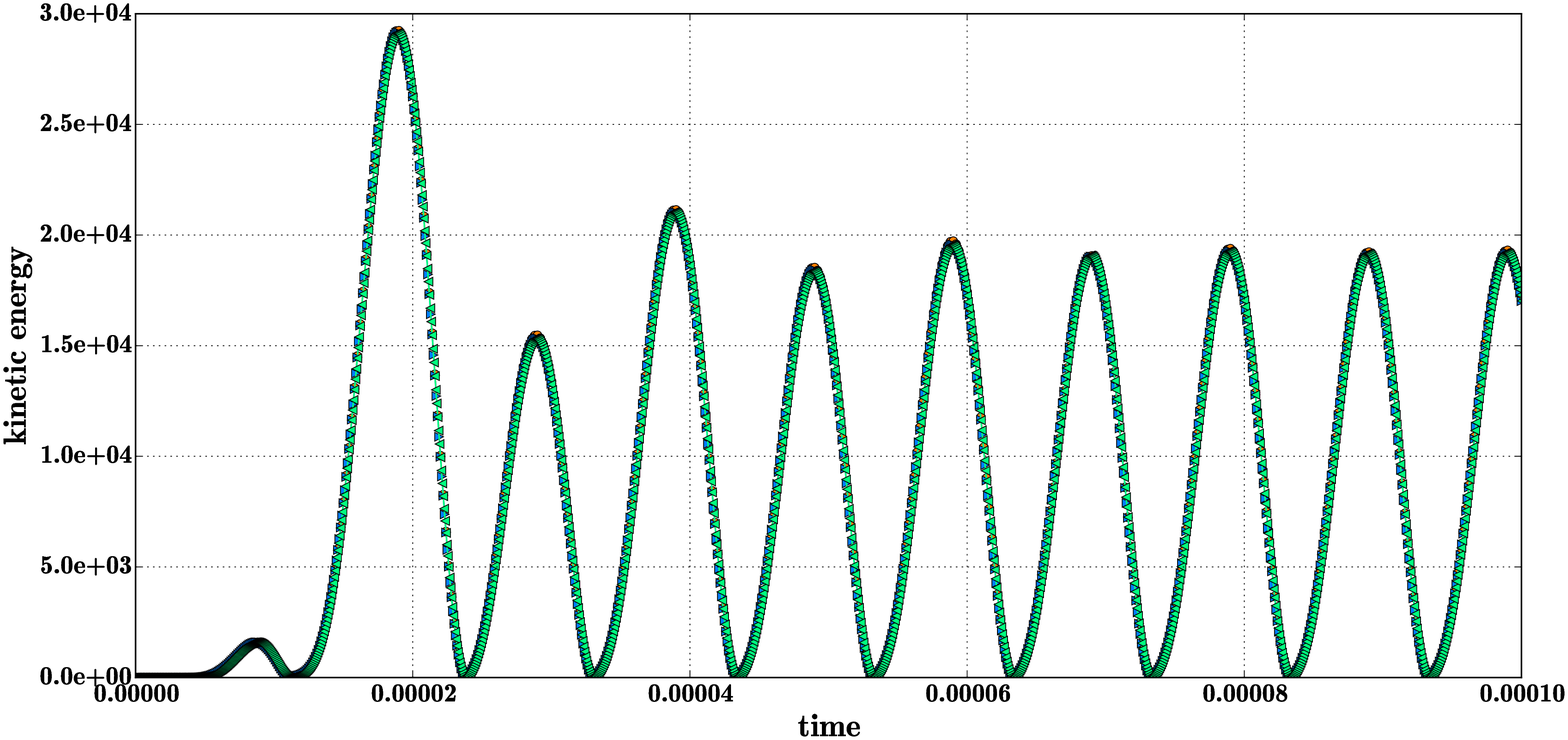}
\hspace*{-1.em}\includegraphics[width=.36\textwidth,height=.25\textwidth]{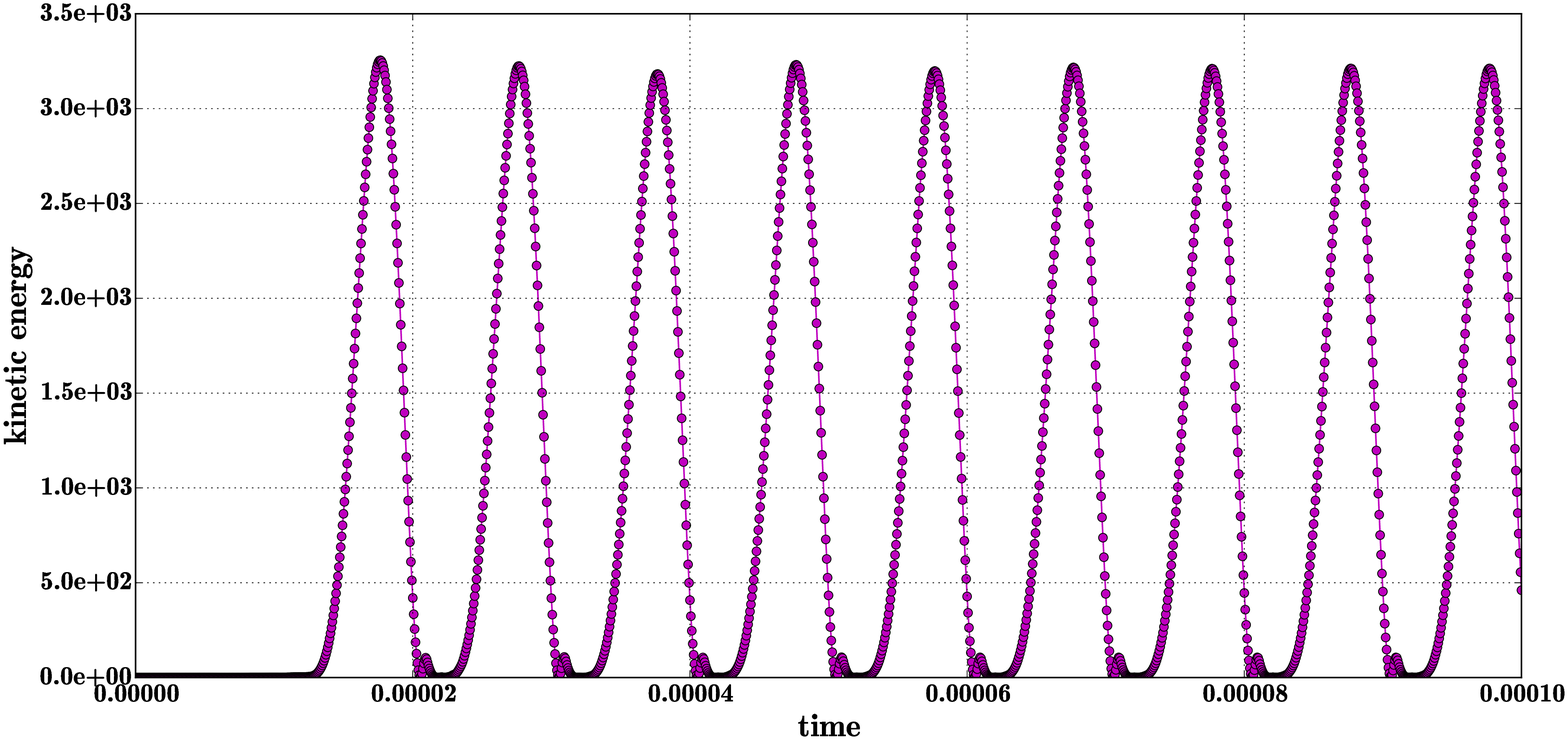}\hspace*{-1.5em}
\\[-.5em]{\bf given load}\hspace{8em}{\bf small friction}\hspace{8em}{\bf large friction}
\end{center}
\vspace*{-1em}
\caption{\sl Computational results of the frictional-contact experiment from Fig.\,\ref{fig4-dam-plast-geom}:
\newline Left: the prescribed cycling force loading $f_1=f_1(t)$.
\newline Middle/Right: 
%
the kinetic energy 
for all of the three meshes used for calculations
(differences practically invisible) for two friction coefficients.
}
\label{fig-loading-f+}
\end{figure}
%
The one-step Crank-Nicolson scheme from Sect.\,\ref{sec-quadratic} is used.
The coarsest time discretisation used 500 time steps, for total time of the 
experiment that of 
1 ms, and then we compare it
also for a twice and three-times finer time discretisation when refining
simultaneously the space discretisation.

\begin{figure}[th]
\begin{center}
\TINY{\hspace*{-.5em}\vspace*{-.1em}{\includegraphics[width=.5\textwidth,height=.3\textwidth]{disp.eps}}\hspace*{-2em}
\includegraphics[width=.6\textwidth,height=.3\textwidth]{velc.eps}\hspace*{-.5em}
\\}
\hspace*{-.5em}\vspace*{-.1em}{\includegraphics[width=.5\textwidth,height=.27\textwidth]{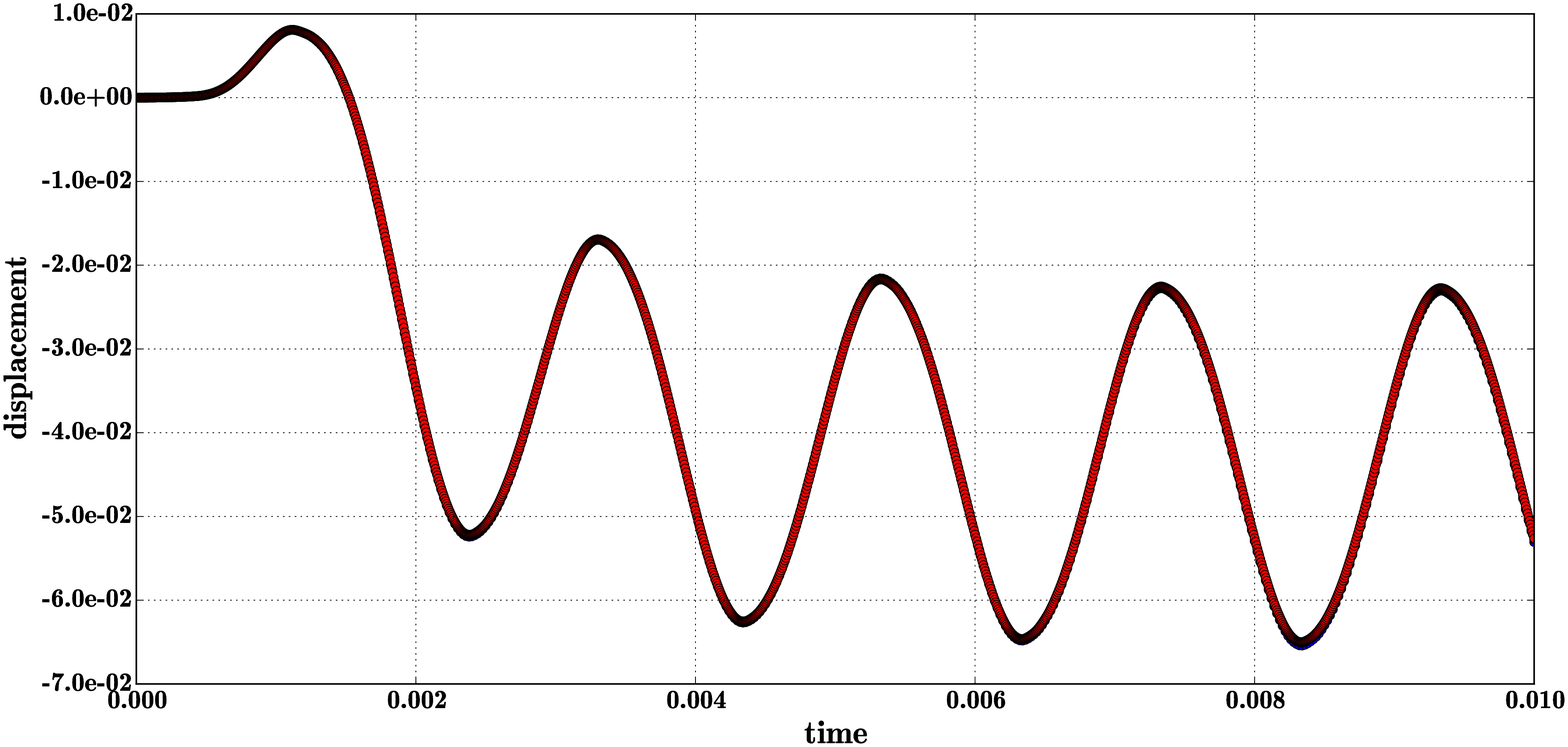}}\hspace*{-2em}
\includegraphics[width=.6\textwidth,height=.27\textwidth]{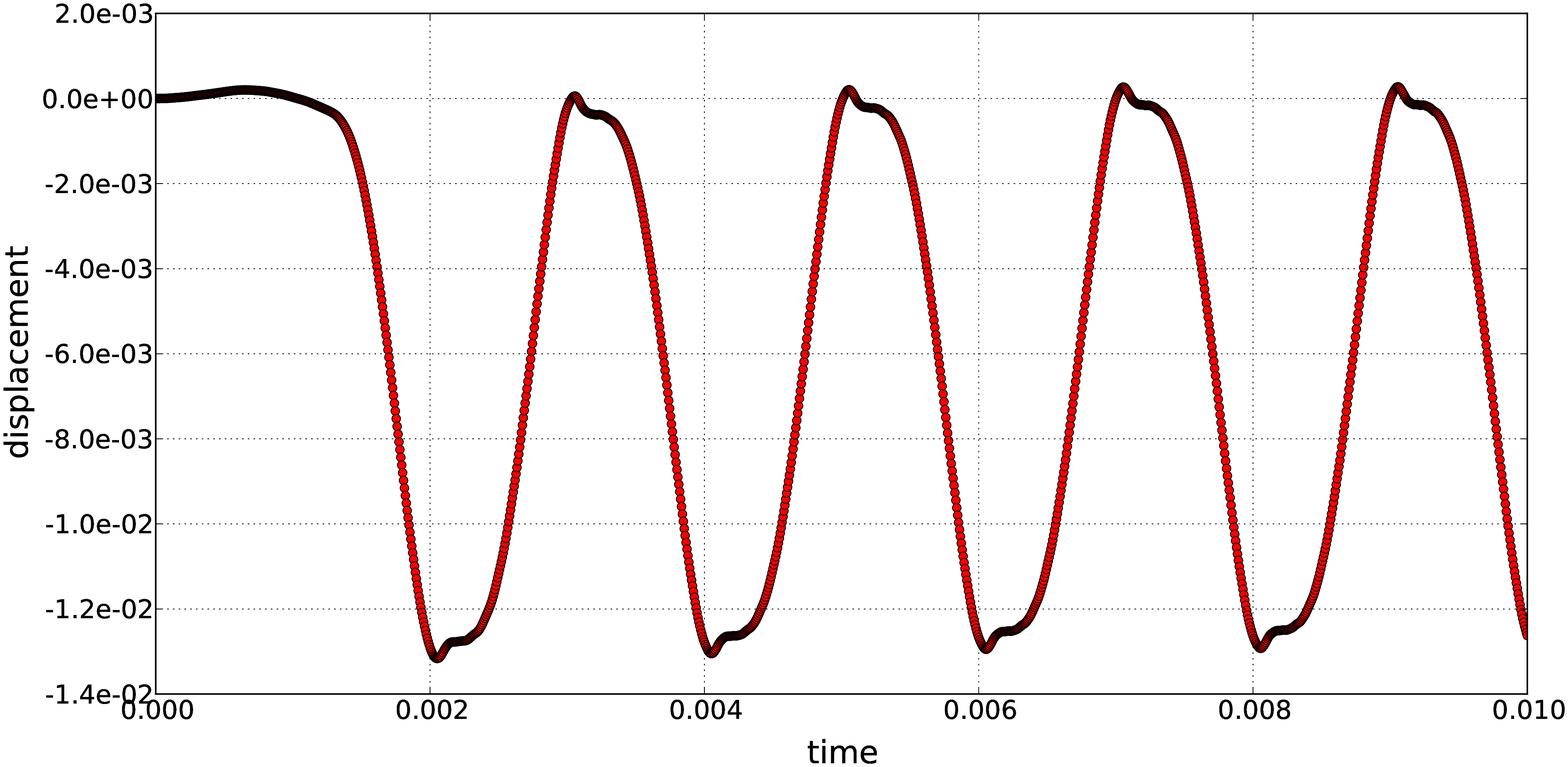}\hspace*{-.5em}
\\\hspace*{-.5em}\vspace*{-.1em}{\includegraphics[width=.5\textwidth,height=.27\textwidth]{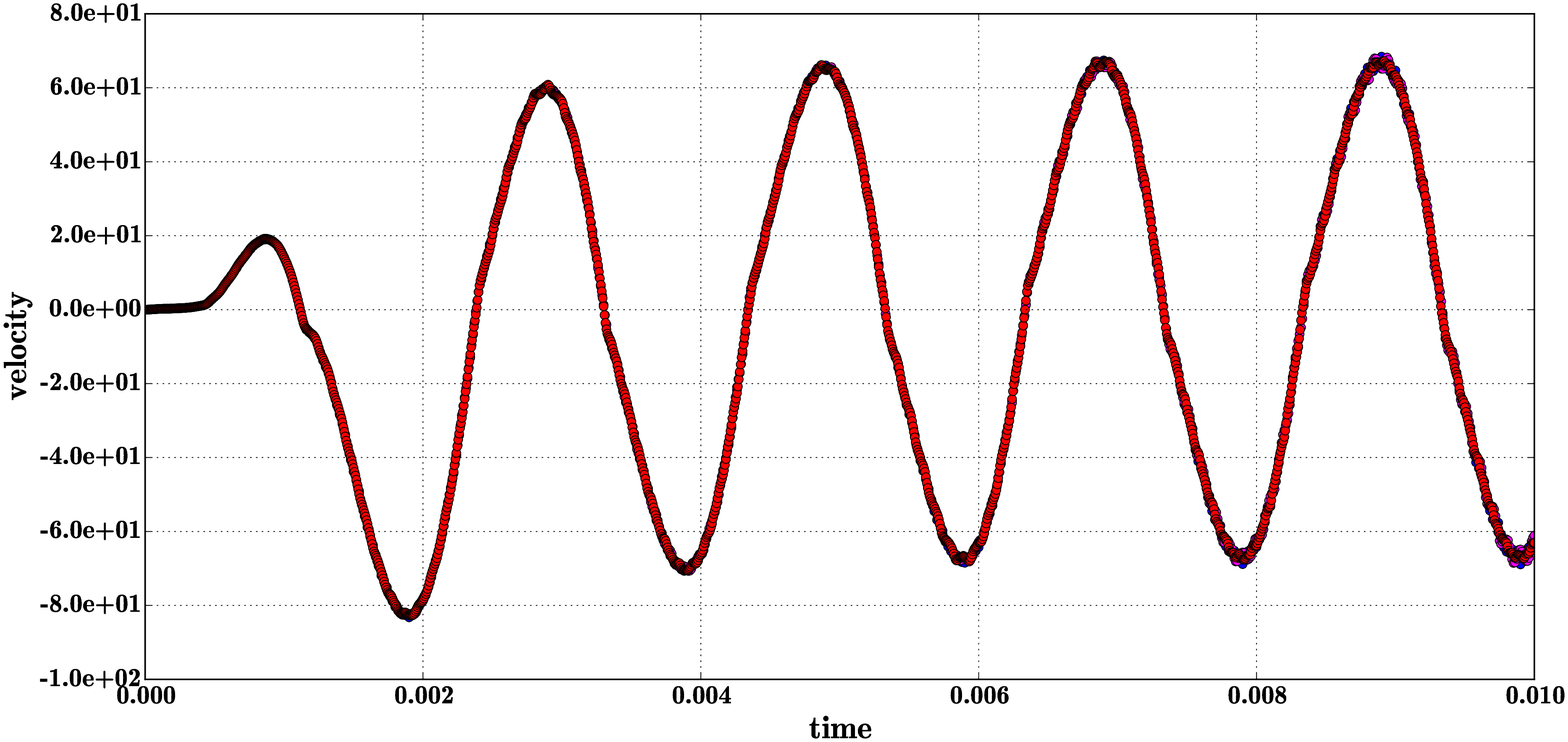}}\hspace*{-2em}
\includegraphics[width=.6\textwidth,height=.27\textwidth]{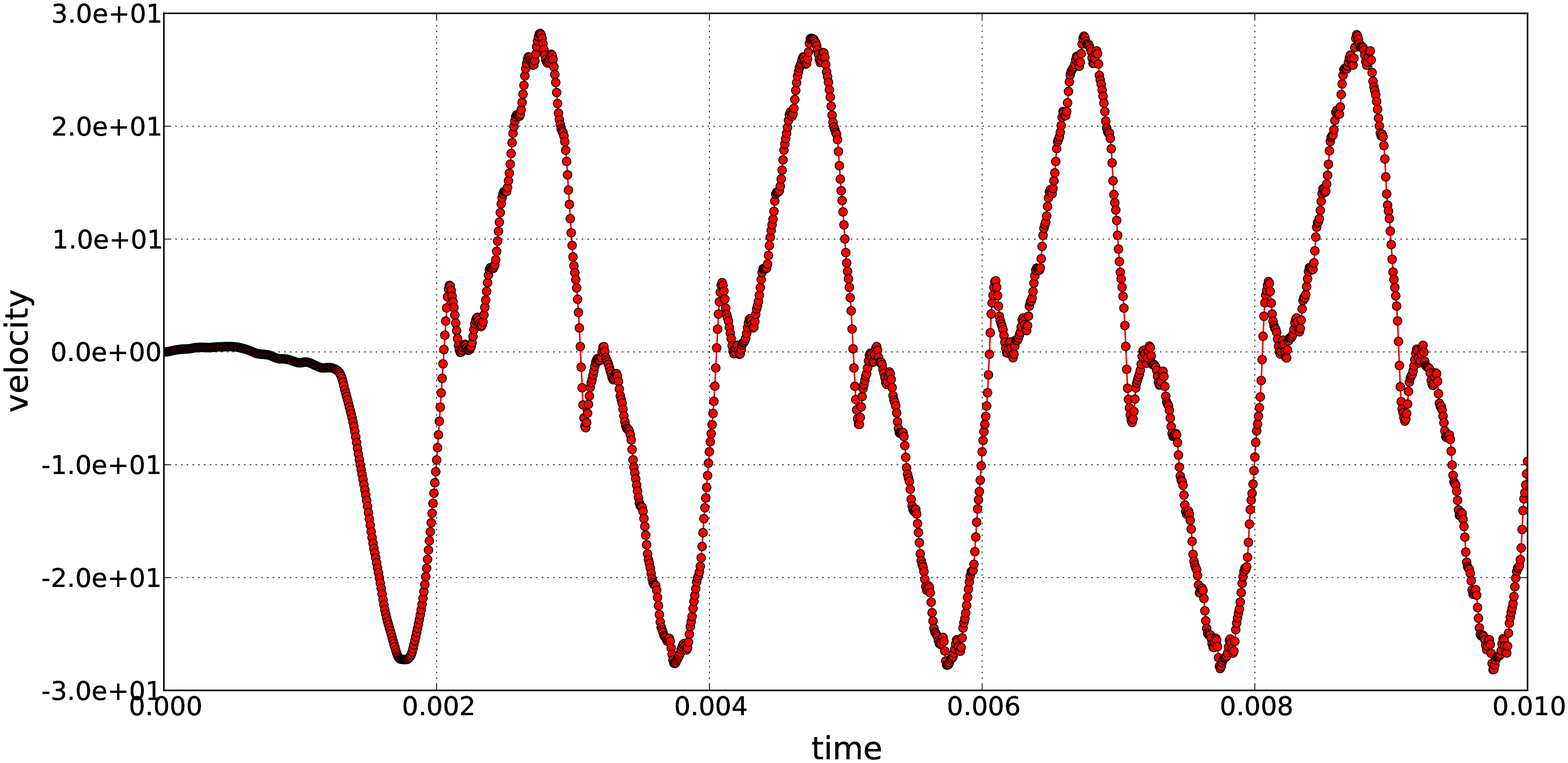}\hspace*{-.5em}
\\[-.5em]{\bf small friction}\hspace{11em}{\bf large friction}
\end{center}
\vspace*{-1em}
\caption{\sl The displacement (upper diagrammes) at the right-hand end of the 
bar 
\COLOR{from Fig.\,\ref{fig4-dam-plast-geom}} and its time derivative, 
i.e.\ the velocity (lower diagrammes).
\TINY{Higher-frequency vibrations of the visco-elastic bar
are clearly visible, being superposed on the lower-frequency loading 
from Fig.\,\ref{fig-loading-f+}-left.} \COLOR{The oscillatory behaviour 
clearly demonstrates the role of inertia preventing 
the immediate blow up after the sliding resistence due to friction is first 
overcome. In addition, high-frequency vibrations in eigen modes of the
bar are visible under large friction (i.e.\ in the right-hand collum).}}
\label{fig-displ-velc}
\end{figure}
\begin{figure}[th]
\begin{center}
\TINY{\hspace*{-5em}\vspace*{-.1em}{\includegraphics[width=.5\textwidth,height=.3\textwidth]{f_u.eps}}
\\}
\hspace*{-2em}\includegraphics[width=.5\textwidth,height=.25\textwidth]{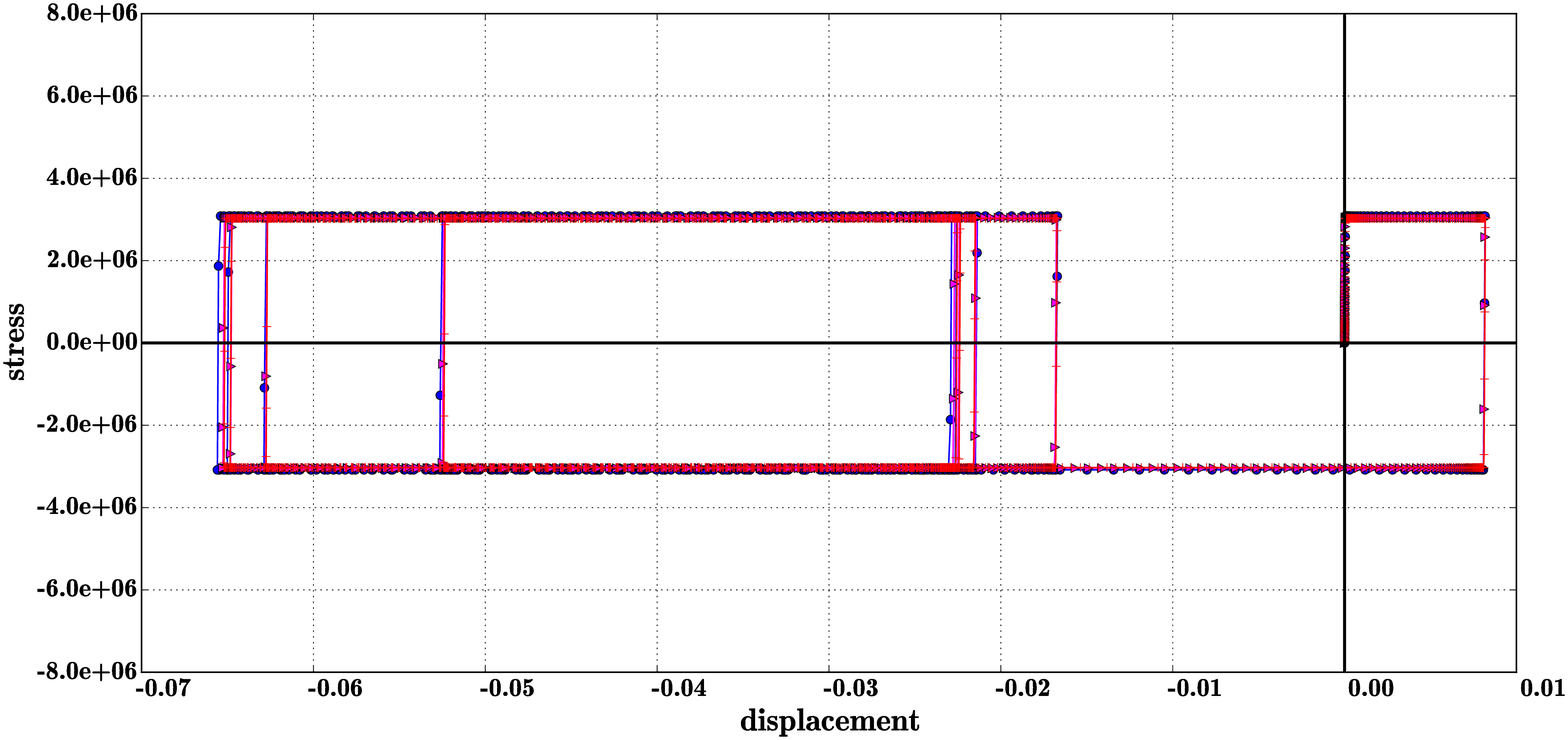}
\includegraphics[width=.5\textwidth,height=.25\textwidth]{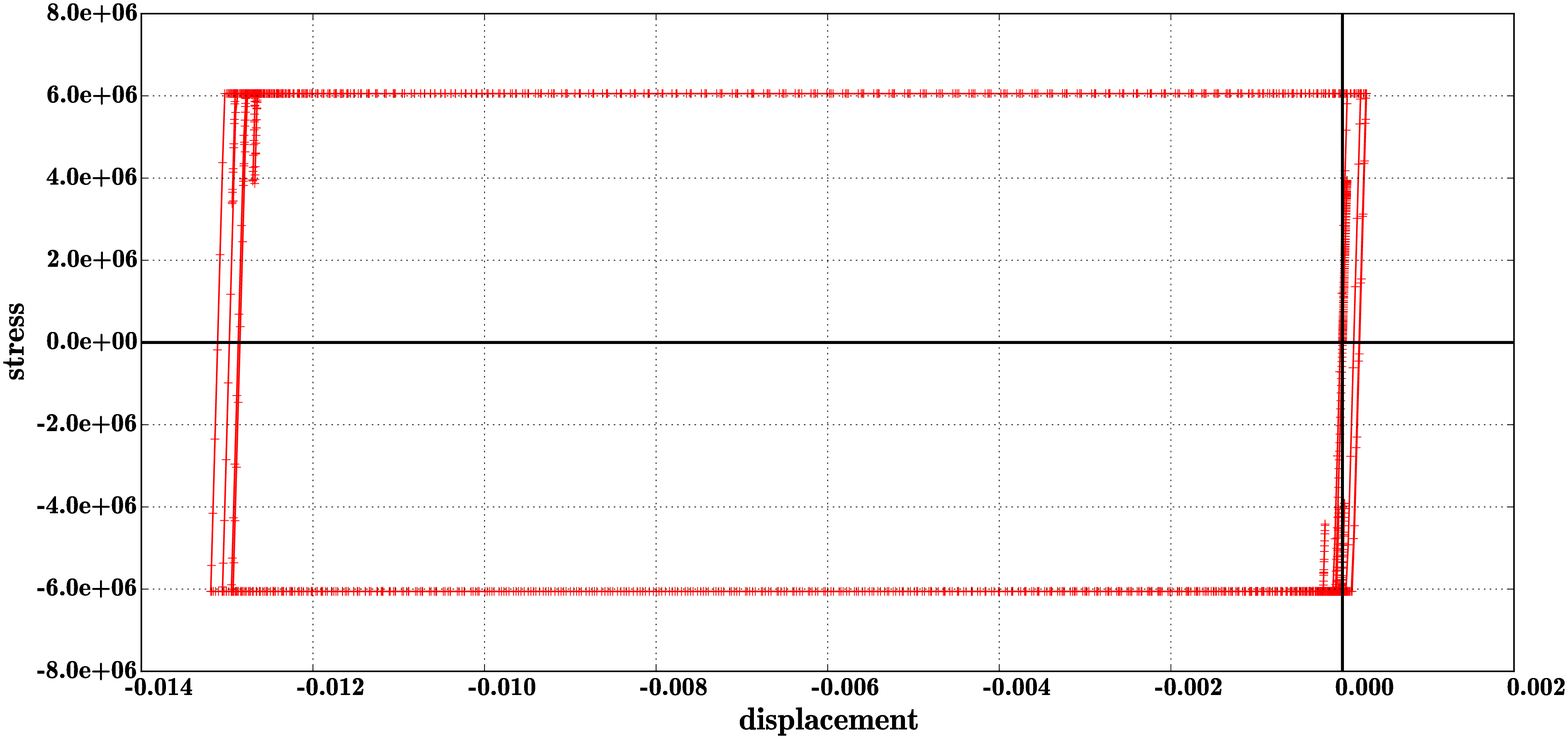}\hspace*{-2em}
\\[-.5em]{\bf small friction}\hspace{15em}{\bf large friction}
\\[-5em]
\end{center}
\vspace*{-1em}
\caption{\sl Typical clock-wise (so-called stop-operator) hysteretic loop 
in the graph displacement versus stress\TINY{\COLOR{ slightly disturbed by higher-frequency 
vibrations visible on Fig.\,\ref{fig-displ-velc}-right}}.}
\label{fig-forc-disp}
\end{figure}\newpage
On the right-hand side $\GNN$ of the specimen we assume a cyclic loading in 
time as depicted in the left plot of Fig.\,\ref{fig-loading-f+}-left. 
Numerical convergence is documented by the plots in 
Fig.\,\ref{fig-loading-f+}-middle/right. 
The computed response on the free right-hand side of the bar is depicted in plots of Fig.\,\ref{fig-displ-velc} both for displacement (left) and velocity (right). In order to check algorithmic performance reults of all the three meshes considered are shown on the same plots. Finally a typical hysteresis loop in the graph displacement versus stress for a material point of the adhesive interface part is shown in Fig.\,\ref{fig-forc-disp}, once again results obtained for all the three considered meshes are given there.
\COLOR{The presence of inertia is essential otherwise the coercivity under the 
mere force loading without any hardening would be lost and instead of an 
oscillating transient response on Figure~\ref{fig-displ-velc}, we would see 
an instantaneous blow-up at time when the friction threshold is reached.} 
%

\subsection{Delamination experiment}\label{sec-delam}
\begin{figure}[th]
\begin{center}
\psfrag{GN}{\footnotesize $\GNNN$}
\psfrag{GD}{\footnotesize $\GNN$}
\psfrag{GC}{\footnotesize $\GC$}
\psfrag{elastic}{\footnotesize elastic body}
\psfrag{obstacle}{\footnotesize rigid obstacle}
\psfrag{adhesive}{\footnotesize adhesive}
\psfrag{LC}{$L_c$}
\psfrag{L}{\footnotesize $L=\ $250\,mm}
\psfrag{H}{\footnotesize $H=$}
\psfrag{12.5}{\scriptsize 12.5\,mm}
\psfrag{loading}{\footnotesize loading}
\hspace*{-.5em}\vspace*{-.1em}{\includegraphics[width=.95\textwidth]{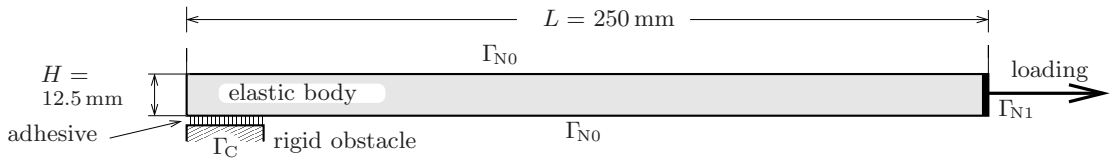}}
\end{center}
\vspace*{-1em}
\caption{\sl Geometry of a 2-dimensional rectangular-shaped specimens 
subjected to the monotonically increasing loading $f_1=f_1(t)$ on the 
right-hand side $\GNN$. 
}
\label{fig4-dam-plast-geom+}
\end{figure}

The geometry is similar to the previous example except that we now consider 
the adhesive contact only on a smaller part of the bottom side (namely 1/10 
of the total length of the specimen, cf.\ Fig.\,\ref{fig4-dam-plast-geom+}, 
and a monotonically increasing loading on the right-hand side \COLOR{ up to a 
complete-rupture time (about 0.7\,s) after which it drops to zero}. 
\COLOR{The  material of the bar as well as of the adhesive is the same 
as in the previous example from Sect.\,\ref{sec-friction},}
except of the fact that in this case we assume Mode-II fracture toughness to have the value $a_2=187.5$ J/m$^2$, while we send $\sigma_y$ to infinity which practically means that 
\COLOR{$\DT\pi_\flat=0$ 
and thus, putting also $\pi_\flat|_{t=0}^{}=0$, the plasticity has been supressed in this experiment.}

This experiment wants to show emission of an elastic wave during sudden 
rupture of the adhesive contact. The fractional-step (with two steps) 
Crank-Nicolson scheme from Sect.\,\ref{sec-nonconvex} is used. In this case 
the coarsest time discretisation used 5000 time steps, for total time of the 
experiment that of 5 sec, and then we compare it
also for a twice and three-times finer time discretisation when refining also the space discretisation. 

\begin{figure}[th]
\begin{center}
\includegraphics[width=.55\textwidth,height=.29\textwidth]{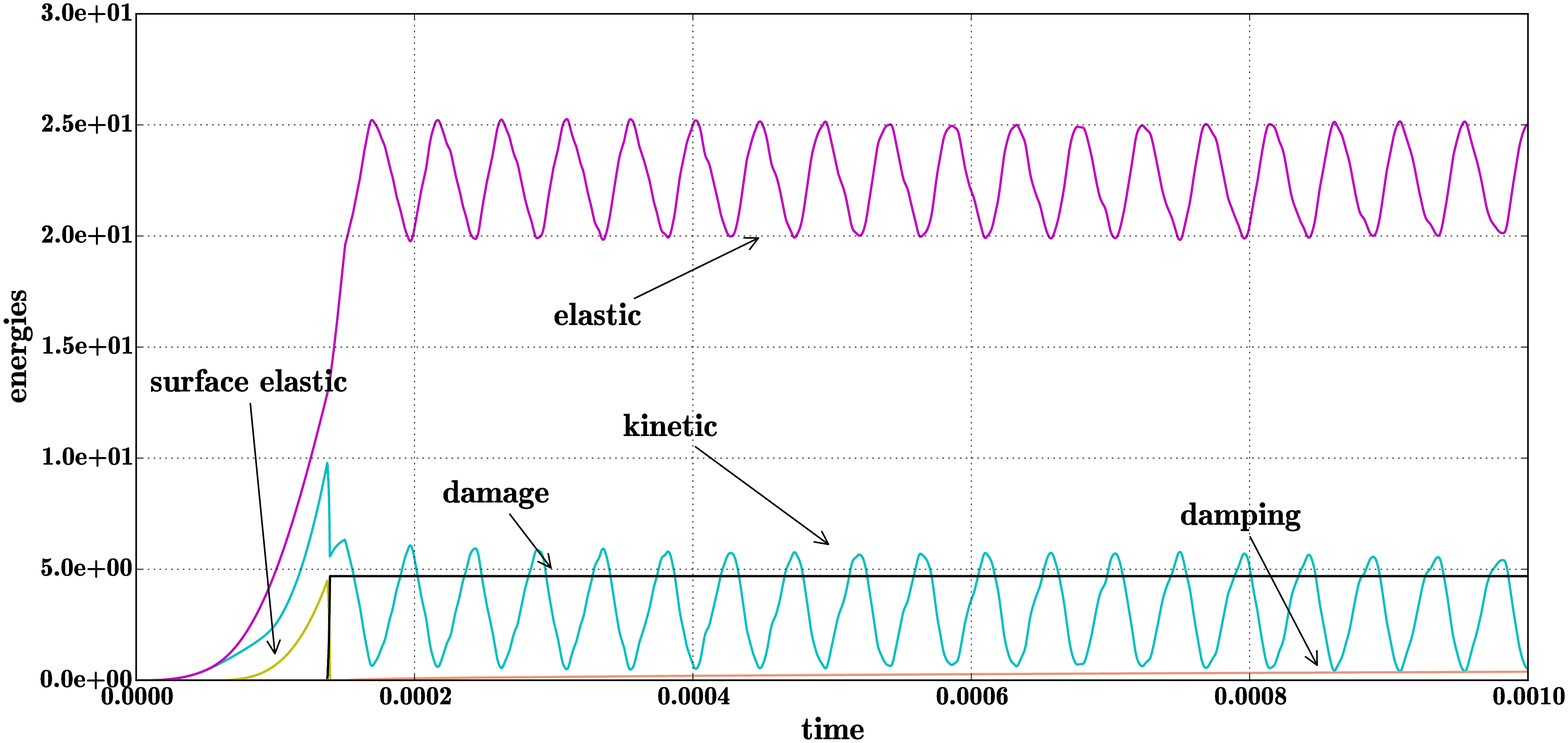}
\includegraphics[width=.43\textwidth,height=.29\textwidth]{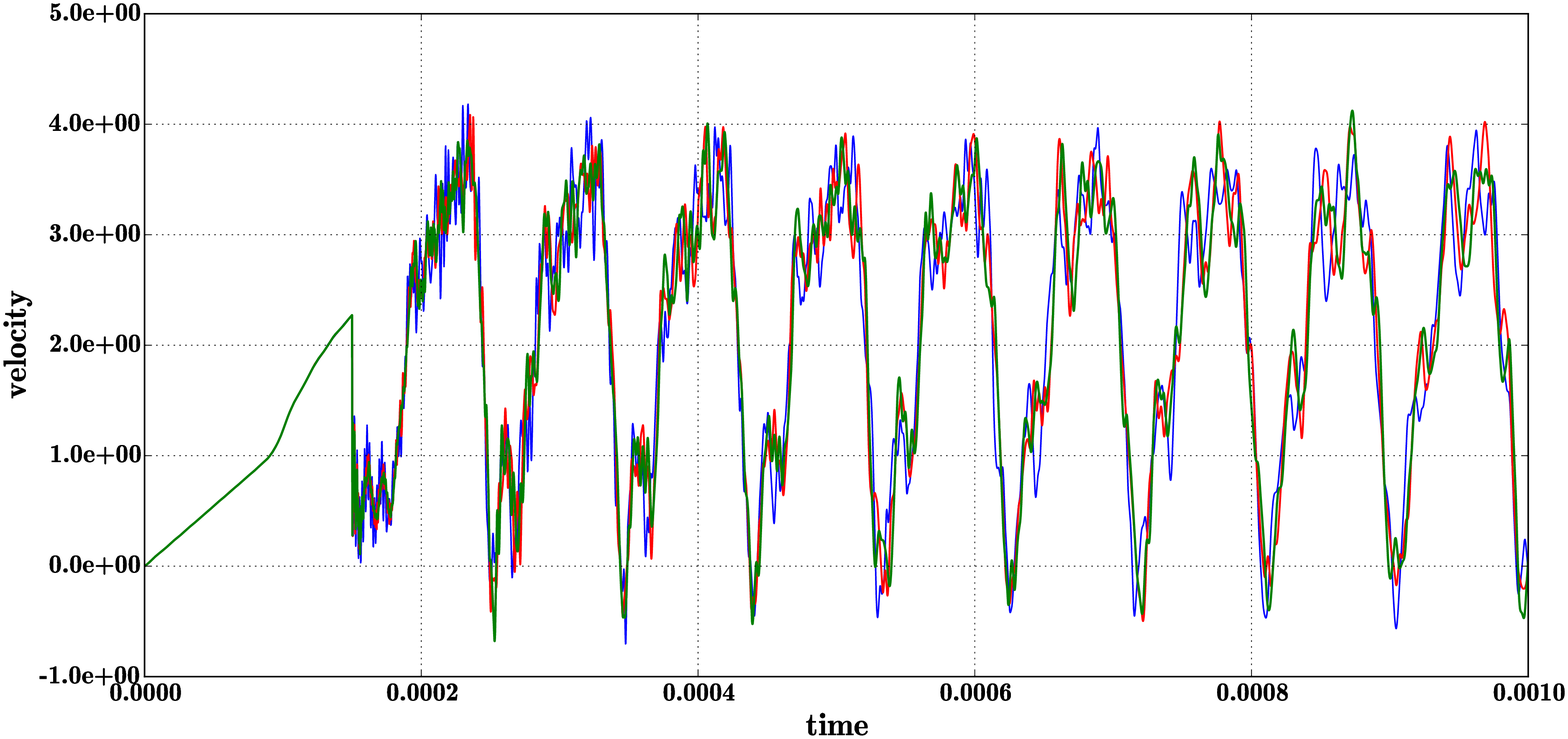}
\end{center}
\vspace*{-1em}
\caption{\sl Results of the delamination experiment from Fig.\,\ref{fig4-dam-plast-geom+}:\COMMENT{THE $t$-SCALE BETTER 0-5s}
\newline Left: time evolution of 
energies of the system showing in particular \COLOR{non-attenuated} 
oscillations of kinetic and 
elastic energies after rupture \COLOR{at time 
$t\doteq 0.14\,{\rm ms}$.}
\newline 
Right: computed velocity of the right-hand side of the domain. \COLOR{The high-frequency vibrations are superposed to the oscillation of the bar on the lowest eigen-frequence.}}
\label{dam-plots_2}
\end{figure}

\begin{figure}[th]
\begin{center}
\includegraphics[width=.55\textwidth,height=.29\textwidth]{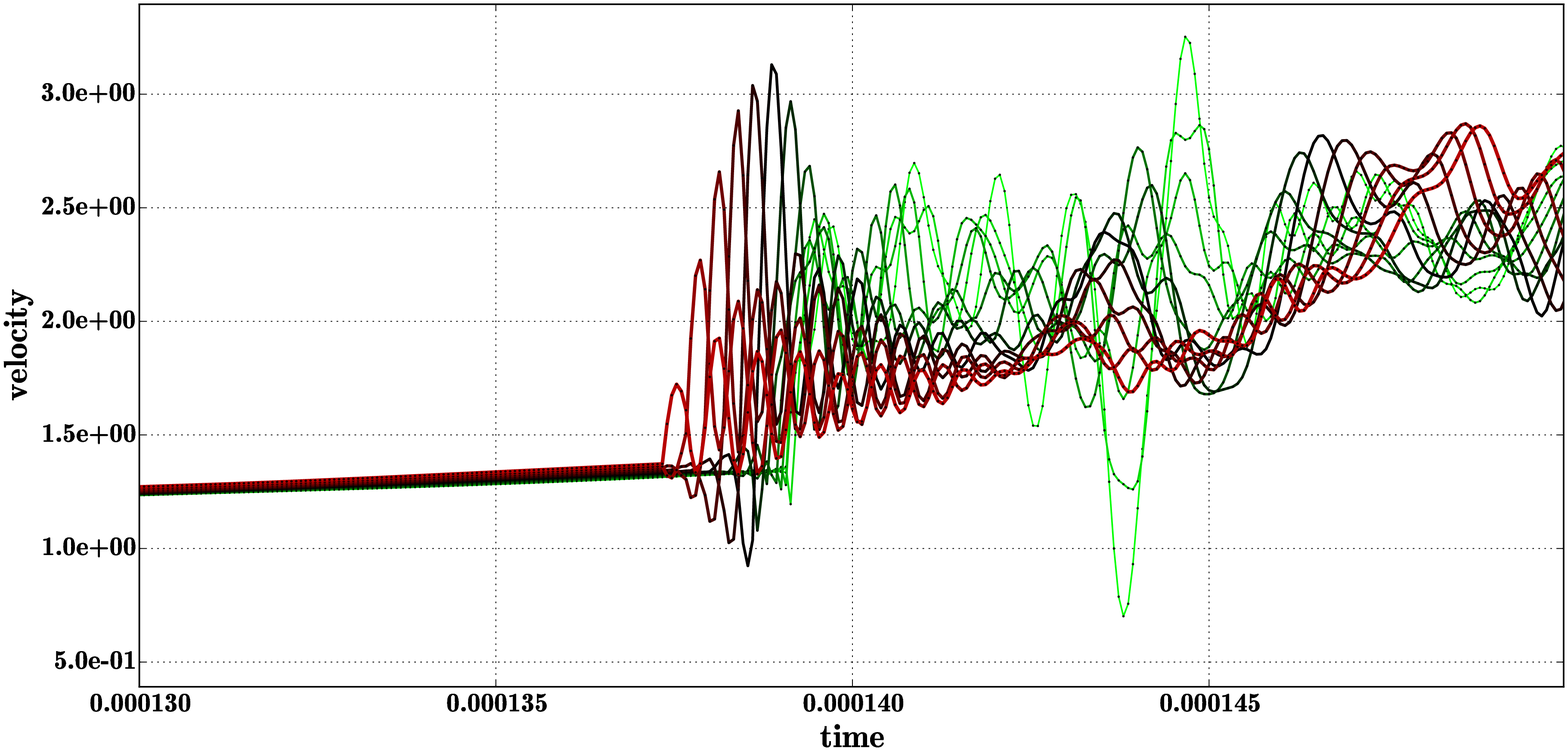}
\includegraphics[width=.43\textwidth,height=.29\textwidth]{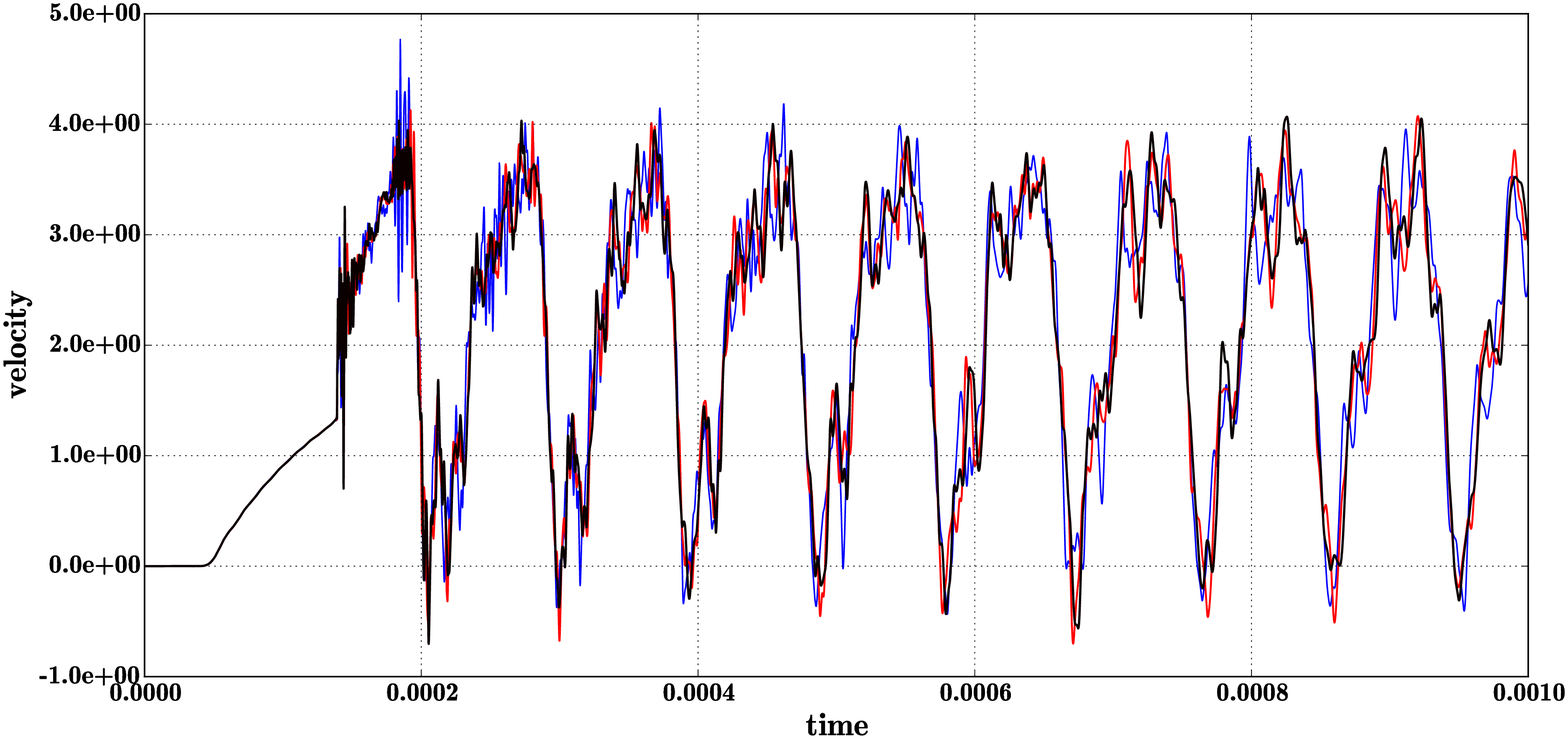}
\end{center}
\vspace*{-1em}
\caption{\sl  Results of the delamination experiment from Fig.\,\ref{fig4-dam-plast-geom+}:
\newline Left: A detail of damage evolution 
\COLOR{on the time interval 
$[0.13\,{\rm ms},0.15\,{\rm ms}]$ around the rupture 
time} tracking the velocity
\COLOR{at particular mesh } 
points on the adhesive.
\newline
Right: velocity like in Fig.~\ref{dam-plots_2}-right but on the left-hand side
of the specimen from Fig.~\ref{fig4-dam-plast-geom+}.
}
\label{dam-plots_3}
\end{figure}

Furthermore, distinguishing energies of the system may be seen in the plot of Fig.\,\ref{dam-plots_2}-left, where more specific elastic, kinetic together with the adhesive stored energy, due to damage dissipation and the damping are shown. Moreover, velocity of the right-hand side of the domain, as computed using all the three meshes, is depicted in Fig.\,\ref{dam-plots_2}-right.
Finally, velocity of the left-most bottom point is depicted in Fig.\,\ref{dam-plots_3}-left. Result are given for all the meshes, while in plot of  Fig.\,\ref{dam-plots_3}-right, results computed using the finer mesh are given for sequentially point on the adhesive. It might observed there the sequence of damage evolution and furth motion after that.


\TINY{
Note also that, not having the gradient of damage, the stored energy is affine
in $\zeta_\flat$. Then both options in \eq{disc-hyper-system-d}, 
i.e.\ the Crank-Nicolson formula and the backward-Euler formula, coincide 
with each other. In particular, 
the energy conservation is granted in 
our time/space discrete scheme.
}

\COLOR{
\begin{remark}[{\it Handling the constraint $\zeta_\flat\ge0$.}]
\upshape
\MARGINOTE{HERE THE PAPER F.Armero, E.Pet\"ocz: Formulation and analysis of 
conserving algorithm for frictionless dynammic contact/impact problems 
STILL TO CHECK ONCE AGAIN}
A more conventional handling of the constraint $\zeta_\flat\ge0$ is
to implement it as a variational inequality and then one can 
consider $A(\cdot)$ linear to model a standard adhesive delamination. 
We can obtain it only as a limit from our model\COMMENT{ if ......}
Noteworthy, this variational-inequality limit problem is not compatible
with our energy-conserving discretisation, although the discretisation of 
all the approximate problems exactly conserves energy. 
In fact, for simplicity, the calculations presented in this 
Section~\ref{sec-delam} have been performed by this variational-inequality
limit problem, the violation of energy conservation only
during the short time interval of rupture having been practically invisible 
below 0.1\%.
\end{remark}
}

\begin{remark}[{\it Combination of damage and plasticity.}]\label{ram-dam+plast}
\upshape
One can combine both inelastic processes like in the bulk model mentioned
in Section~\ref{sec-cont-mech}.
One would then deal with the model 
from \cite{RoKrZe13DACM,RoMaPa13QMMD} originally considered in quasistatic 
variant and devised to distinguish
Mode I and Mode II delamination. Here we considered intentionally only
Mode II otherwise one should rather consider a unilateral contact which 
would then be incompatible with our assumption of $\Phi$ being 
quadratic in terms of $u$. Here $\kappa>0$ 
is needed to facilitate the analysis. 
\end{remark}

\COMMENT{SOMEWHERE STILL VERIFICATION OF THE continuity \eqref{ass-for-Phi}
for $\gamma(\cdot)$ smooth with $\gamma'(\cdot)$ bounded}

\baselineskip=12pt
{\small
\noindent{\it Acknowledgments:} This research was supported from the grants 
14-15264S ``Experimentally justified multiscale modelling of 
shape memory alloys'' 
and 16-03823S ``Homogenization and multi-scale computational modelling of 
flow and nonlinear interactions in porous smart structures''
of the Czech Science Foundation, and from the 
institutional support RVO:\,61388998 (\v CR).}

\bibliographystyle{plain}
\bibliography{trcreta}

\end{sloppypar}
\end{document}